\documentclass[a4paper,11pt]{article} 

\usepackage[french]{babel}  

\usepackage{amssymb,amsmath,amsthm}
\usepackage{nccmath}
\usepackage{graphicx,floatrow}
\usepackage[margin=2.3cm]{geometry} 

\def\ZZ{\mathbb{Z}} 
\def\QQ{\mathbb{Q}}
\def\RR{\mathbb{R}}
\def\CC{\mathbb{C}}
\def\NN{\mathbb{N}}

\def\cL{\mathcal{L}}
\def\cA{\mathcal{A}}

\def\dt{\, \mathrm{d}t}
\def\ds{\, \mathrm{d}s}
\def\u{\mathsf{u}} 
\def\v{\mathsf{v}}
\def\w{\mathsf{w}}

\renewcommand{\bold}{\boldsymbol}

\def\({\langle} 
\def\){\rangle} 

\def\Cyl{\operatorname{Cyl}}

\def\epsilon{\varepsilon}
\def\phi{\varphi}

\title{Solutions born\'ees et auto-similaires d'\'equations int\'egrales}
\author{Jean-Fran{\c c}ois Bertazzon et Vincent Delecroix}
\author{
	Jean-Fran\c cois Bertazzon
	\footnote{
		Lyc\'ee Notre-Dame de Sion,
		231 Rue Paradis, 13006 Marseille, France.		
		{\it E-mail addresse:} {\tt jeffbertazzon@gmail.com}
		}
	\quad et \quad  
	Vincent Delecroix
	\footnote{
		LaBRI, UMR 5800
		B\^atiment A30
		351, cours de la Lib\'eration
		33405 Talence cedex, France.		
		{\it E-mail addresse:} {\tt vincent.delecroix@labri.fr}
	}
}

\date{}

\begin{document}

\renewcommand{\arraystretch}{1.5}

\makeatletter
\renewcommand\section{\@startsection{section}{1}{\z@}{6pt\@plus0pt}{6pt}{\scshape\centering \large}}
\renewcommand\subsection{\@startsection {subsection}{2}{\z@}{6pt \@plus 0ex  \@minus 0ex}{-6pt \@plus 0pt}{\reset@font\itshape \ \ \ \ \ \ \ \ }}
\makeatother

\renewcommand{\proofname}{\textup{Preuve}}
\newcommand{\Floor}[1]{{\left\lfloor #1 \right\rfloor}}

\newtheorem{definition}{D\'efinition}
\newtheorem{lemme}{Lemme}
\newtheorem{theoreme}{Th\'eor\`eme}
\newtheorem{corollaire}[lemme]{Corollaire}
\newtheorem{proposition}[lemme]{Proposition}

\theoremstyle{definition}
\newtheorem{remarque}{Remarque}

\maketitle

\renewcommand{\abstractname}{\textup{R\'esum\'e}}
\begin{abstract}
Nous construisons des solutions born\'ees aux \'equations int\'egrales
\[
\int_0^{\lambda x} f(t)\, \dt= f(x) - f(0)
\]
o\`u $\lambda \geq 2$ est un entier. Cette construction s'appuie sur une m\'ethode
originale de limite de sommes de Birkhoff it\'er\'ees. 
\end{abstract}

\section{Introduction} \label{se:intro}

Nous nous int\'eressons aux \'equations int\'egrales :
\begin{equation}
\int_0^{\lambda x} f(t)\,\dt = \delta \big(f(x) - f(0)\big) \qquad \text{pour tout $x \geq 0$.}
\tag{$E_{\lambda,\delta}$}
\label{eq:E_lambda_delta}
\end{equation}
Remarquons tout d'abord que le param\`etre $\delta$ est accessoire : si $f$ est solution de~\eqref{eq:E_lambda_delta} alors $g(x) = f(|\delta| x)$ est solution de $(E_{\lambda,1})$ ou $(E_{\lambda,-1})$. Cependant, il s'av\`ere naturel de consid\'erer des param\`etres $\lambda$ et $\delta$ entiers dans notre construction. On pr\'ef\'erera parfois la formulation \'equivalente
\begin{equation*}
f'(x) = \frac{\lambda}{\delta} f(\lambda x) \qquad \text{pour tout $x \geq 0$.}
\tag{$E'_{\lambda,\delta}$}
\label{eq:Ebis_lambda_delta}
\end{equation*}
L'\'equation ci-dessus est un cas particulier de \emph{l'\'equation du pantographe} dont la forme g\'en\'erale est 
\[
f'(x) = a f(\lambda x) + b f(x) \quad \text{avec $a,b \in \RR$ et $\lambda \in \RR_+$.}
\]
Nous renvoyons \`a l'introduction de~\cite{BogachevaDerfelbMolchanovcOckendond} pour une bibliographie
r\'ecente sur le sujet.

\bigskip

Si $0<\lambda\leq1$, l'espace des solutions de l'\'equation \eqref{eq:E_lambda_delta} est de dimension $1$ et les solutions sont analytiques.
Dans le cas qui nous int\'eresse, $\lambda>1$, le th\'eor\`eme suivant montre que l'espace des solutions est de dimension infinie.
\begin{theoreme} \label{thm:prolong}
Soit $f \in C^\infty([1,\lambda])$ telle que, pour tout entier $n$, on ait $f^{(n)}(1) =  f^{(n)}(\lambda)=0$.
Alors $f$ se prolonge de mani\`ere unique en une fonction $C^\infty([0,\infty[)$ solution de \eqref{eq:E_lambda_delta}.
\end{theoreme}

Nous montrons que les solutions ne peuvent pas \^etre p\'eriodiques et s'annulent n\'ecessairement infiniment souvent si $\delta > 0$.
\begin{theoreme} \label{thm:proprietes}
Soit $f$ une solution de l'\'equation~\eqref{eq:E_lambda_delta} avec $\lambda > 1$.  
\begin{itemize}
\item Si $\delta \not= 0$ et $f$ est p\'eriodique, alors $f=0$,
\item Si $\delta > 0$, alors pour tout $t > \dfrac{\delta}{\lambda (\lambda -
1)}$ la fonction $f$ s'annule au moins une fois dans l'intervalle
$\left[\dfrac{t}{\delta}, \dfrac{t}{\delta} \lambda^3 \right]$.
\end{itemize}
\end{theoreme}
L'objectif principal de cet article est la construction de solutions born\'ees et auto-similaires de ces \'equations. Donnons nous deux fonctions $f_a: [0,\lambda] \to \RR$ et $f_b:[0,\lambda] \to \RR$ v\'erifiant $f_a(0) = f_a(\lambda) = f_b(0) = f_b(\lambda)$. Pour un mot infini
$\w = \w_0 \w_1 \ldots \in \{a,b\}^\NN$, nous noterons $f_\w: \RR_+ \to \RR$ la fonction d\'efinie par
\[
f_\w(n \lambda + x) = f_{\w_n}(x) \quad \text{pour $x \in [0,\lambda]$ et $n \in \NN$.}
\]
On appelera $f$ la \emph{concat\'enation de $f_a$ et $f_b$ le long de $\w$}.
Notons que cette d\'efinition fait aussi sens pour des mots $\w$ finis, auquel
cas la fonction obtenue est d\'efinie sur le  domaine $[0, |\w|\lambda]$ o\`u $|\w|$ d\'esigne
la longueur du mot $\w$. 
Pour les mots de longueur $1$, nous retrouvons bien les fonctions de d\'eparts $f_a$ et $f_b$.

Nos solutions sont construites sur des mots $\w$ substitutifs.
Une \emph{substitution} $\sigma$ sur $\{a,b\}$ est un morphisme du mono\"ide libre $\{a,b\}^*$. 
Notons qu'une substitution est enti\`erement d\'etermin\'ee par les images des g\'en\'erateurs $\sigma(a)$ et $\sigma(b)$. 
Notons \'egalement que l'action d'une substitution s'\'etend en une action sur les mots infinis (nous renvoyons \`a la partie \ref{subsec:mots_infinis} pour plus de d\'etails).

Prenons la substitution $\sigma: a \mapsto ab, b \mapsto ba$ dite \emph{de Prouet-Thue-Morse}. 
Les premi\`eres it\'erations de cette substitution sur la lettre $a$ donnent : 
\[
\sigma^1(a) = ab,\quad \sigma^2(a) = abba, \quad \sigma^3(a) = abbabaab, \quad \sigma^4(a) = abbabaabbaababba  \ldots
\]
On obtient une suite de pr\'efixes embo\^it\'es qui convergent (pour la topologie produit) vers un mot infini
\[
\w = abbabaabbaababbabaababbaabbabaabbaababbaabbabaababbabaabb\ldots
\]
appel\'e mot \emph{de Prouet-Thue-Morse}. C'est l'unique mot (infini) $\w$ commen\c{c}ant par $a$ tel que $\sigma(\w) = \w$. Autrement dit c'est un point fixe de $\sigma$. Pour tout mot fini $\u$, nous noterons $|\u|$ sa longueur et $|\u|_\alpha$ le nombre de lettre $\alpha$ dans ce mot.
Une substitution $\sigma$ sur $\{a,b\}$ est dite $\lambda$-\emph{uniforme} si $|\sigma(a)| = |\sigma(b)| = \lambda$.

\begin{theoreme}  \label{thm:sols_subs}
Soit $f_a: [0,\lambda] \to \RR$ et $f_b:[0,\lambda] \to \RR$ et $\w \in \{a,b\}^\NN$ tels que $f = f_\w$ soit solution de \eqref{eq:E_lambda_delta}.
Si $f$ n'est pas identiquement nulle alors il existe une unique substitution $\lambda$-uniforme $\sigma$ telle que $\sigma(\w) = \w$.
De plus $|\sigma(a)|_a = |\sigma(b)|_a$ et pour tout $x,y \in [0,\lambda]$ on a
\[
\int_{\lambda x}^{\lambda y} f_{\sigma(a)}(t) \dt = \delta \Big(f_a(y) - f_a(x)\Big)
\quad \text{et} \quad
\int_{\lambda x}^{\lambda y} f_{\sigma(b)}(t) \dt = \delta \Big(f_b(y) - f_b(x)\Big).
\]
\end{theoreme}

\bigskip

Notre m\'ethode pour construire des solutions de \eqref{eq:E_lambda_delta} repose sur la convergence de sommes de Birkhoff.
Repartons du mot de Prouet-Thue-Morse $\w = abbabaabb\ldots$ et 
consid\'erons  $\phi: \{a,b\} \rightarrow \RR$ la fonction d\'efinie par $\phi = \chi_a - \chi_b$ o\`u $\chi_\alpha$ d\'esigne la fonction qui vaut un pour $\alpha$ et $0$ sinon. Il faut ici voir la suite $n \mapsto \phi(\w_n)$ comme une fonction d\'efinie sur les entiers. Une version discr\`ete de l'int\'egrale consiste \`a prendre des \emph{sommes de Birkhoff} de cette fonction : on pose $S^{(1)}_0(\phi,\w)=0$ et 
\[
S^{(1)}_n(\phi,\w) = \sum_{k=0}^{n-1} \phi(\w_k) \quad \text{pour tout entier $n\geq 1$.}
\]
\begin{figure}[H]
\begin{center}
\includegraphics[width=0.95\textwidth]{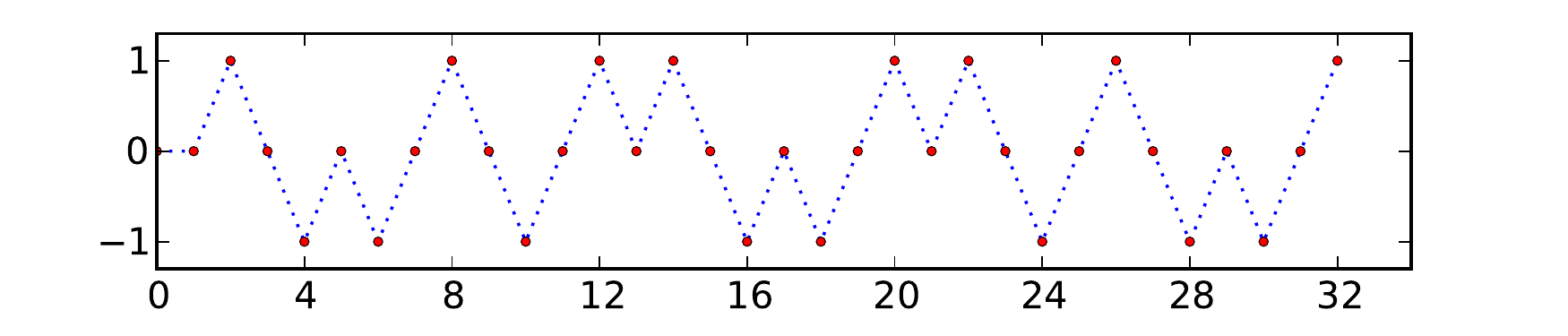} \\
\includegraphics[width=0.95\textwidth]{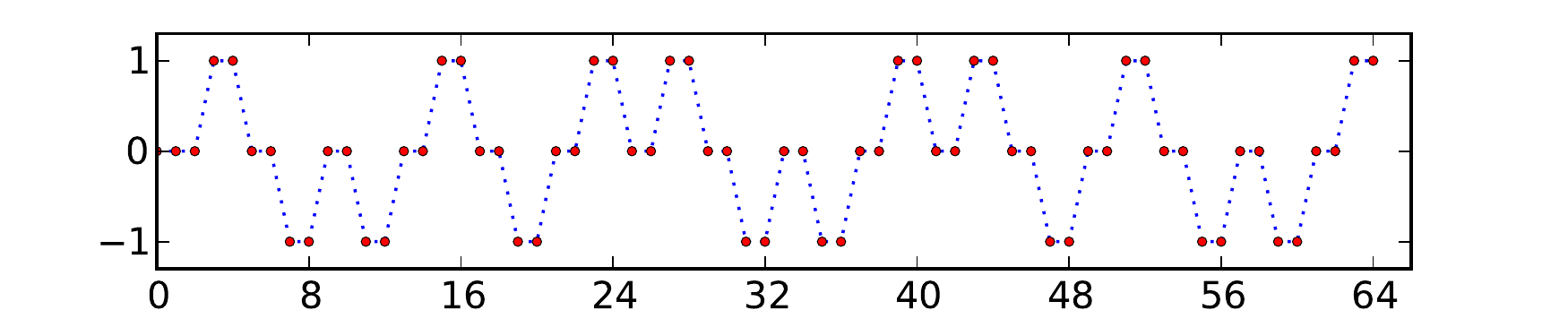} \\
\includegraphics[width=0.95\textwidth]{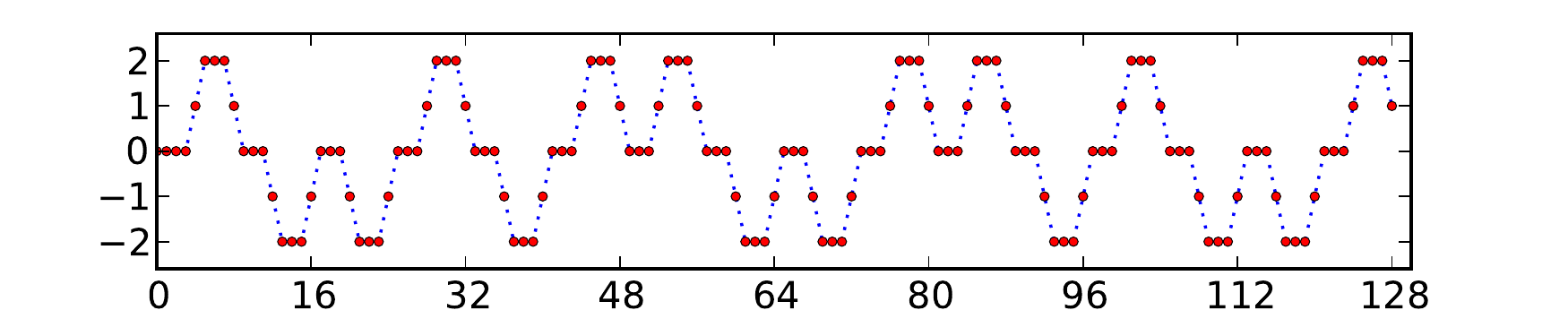} \\
\end{center}
\caption{Les sommes $S^{(k)}(\phi,\w)$ pour $k=1,2,3$ pour le mot de Prouet-Thue-Morse $\w$ et $\phi = \chi_a - \chi_b$.}
\label{fig:ptm}
\end{figure}

\begin{figure}[H]
\begin{center}
\includegraphics[width=0.95\textwidth]{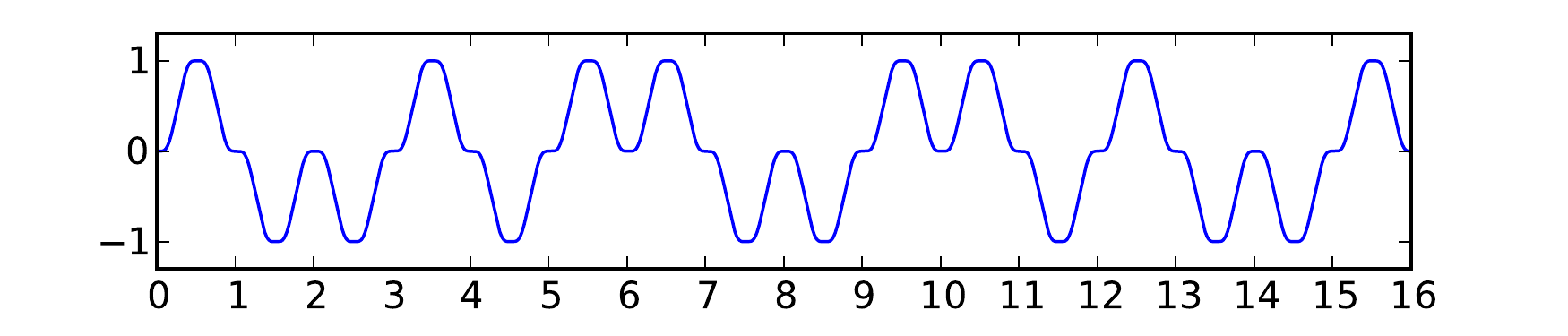}
\end{center}
\caption{La fonction limite $f_{\w}$ pour le mot de Prouet-Thue-Morse $\w$ obtenue avec notre construction pour la fonction $\phi = \chi_a - \chi_b$ (fonction de Fabius~\cite{Fabius66}.).}
\label{fig:f_ptm}
\end{figure}

Les quantit\'es $S^{(1)}_n(\phi,\w)$ peuvent \`a nouveau \^etre consid\'er\'ees comme une fonction $n \mapsto S^{(1)}_n(\phi,\w)$. 
On construit alors par r\'ecurrence une suite de fonctions discr\`etes en d\'efinissant $S^{(\ell+1)}(\phi,\w)$ 
	comme la somme de Birkhoff de $S^{(\ell)}(\phi,\w)$. 
Nous montrons que cette suite de fonctions discr\`etes proprement renormalis\'ee converge vers une solution de $(E_{2,1})$ 
	(voir les figures~\ref{fig:ptm} et~\ref{fig:f_ptm}).

Plus g\'en\'eralement, pour tout mot infini $\w$ sur $\{a,b\}$ et toute fonction $\phi:\{a,b\} \rightarrow \RR$, nous pouvons d\'efinir ses sommes de Birkhoff $S^{(\ell)}_n(\phi,\w)$. 
Notons que si $\w = \w_0 \w_1 \ldots \w_{m-1}$ est un mot fini, la somme de Birkhoff $S^{(i)}_n(\phi,\w)$ a bien un sens tant que $n$ est inf\'erieur ou \'egal \`a la longueur $m$ du mot $\w$.
Pour d\'eterminer le param\`etre $\delta$ de l'\'equation fonctionnelle, nous introduisons la d\'efinition suivante.
\begin{lemme} \label{le:def_delta}
Soit $\sigma$ une substitution $\lambda$-uniforme sur $\{a,b\}$. Soit $\phi: \{a,b\} \rightarrow \RR$ telle que $\phi(a) \not= \phi(b)$, alors pour tout $i \geq 1$, la quantit\'e
\[
\frac{S^{(i)}_\lambda(\phi, \sigma(a)) - S^{(i)}_\lambda(\phi, \sigma(b))}{\phi(a) - \phi(b)}
\]
est ind\'ependente de $\phi$ (telle que $\phi(a) \not= \phi(b)$). On notera $\delta_i(\sigma)$ cette quantit\'e et 
\[ 
\quad \delta(\sigma) = \Big(\delta_1(\sigma),  \ldots,\delta_\lambda(\sigma)\Big).
\]
\end{lemme}
Nous renvoyons le lecteur \`a la section~\ref{subsec:sommes_iterees} pour la preuve de ce lemme.
Remarquons cependant que $\delta_1(\sigma) = 0$ si et seulement si $\sigma(a)$ et $\sigma(b)$ contiennent le m\^eme nombre de $a$ (et \`a fortiori le m\^eme nombre de $b$).

\begin{figure}[H]  
\begin{center} \includegraphics[width=0.95\textwidth]{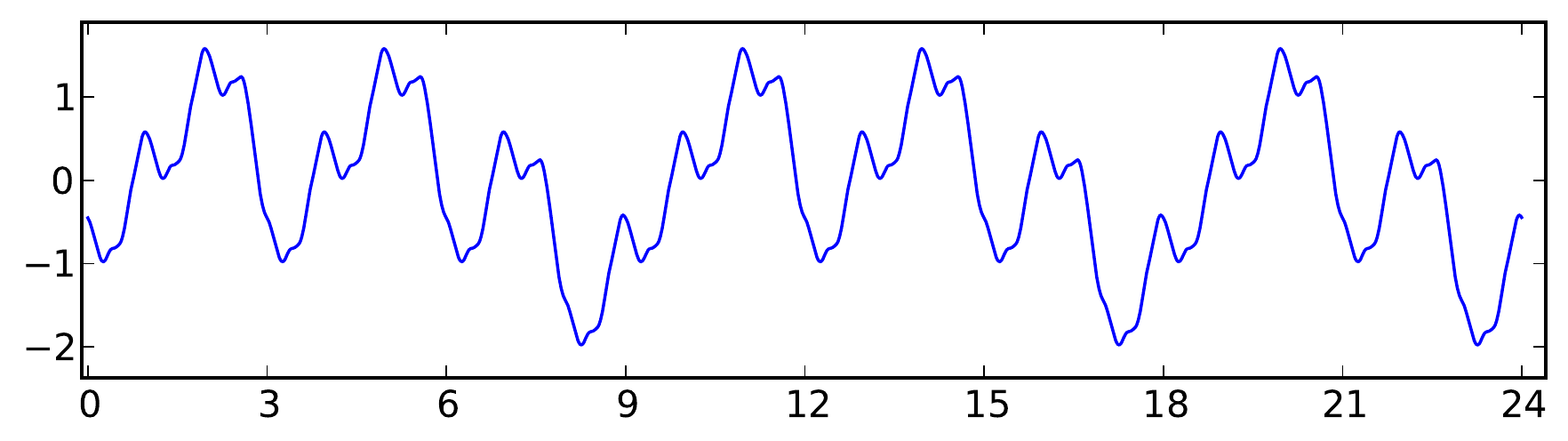} \end{center}
\caption{Fonction $f_{\w}$ pour le point fixe commen\c{c}ant par $a$ de la substitution 
	$(a \mapsto aab, b \mapsto aba)$ obtenue avec notre construction pour la fonction  $\phi = \chi_a - 2 \chi_b$.}
\label{fig:exemple_f_sigma1}
\end{figure}

\begin{figure}[H]
\begin{center} \includegraphics[width=0.95\textwidth]{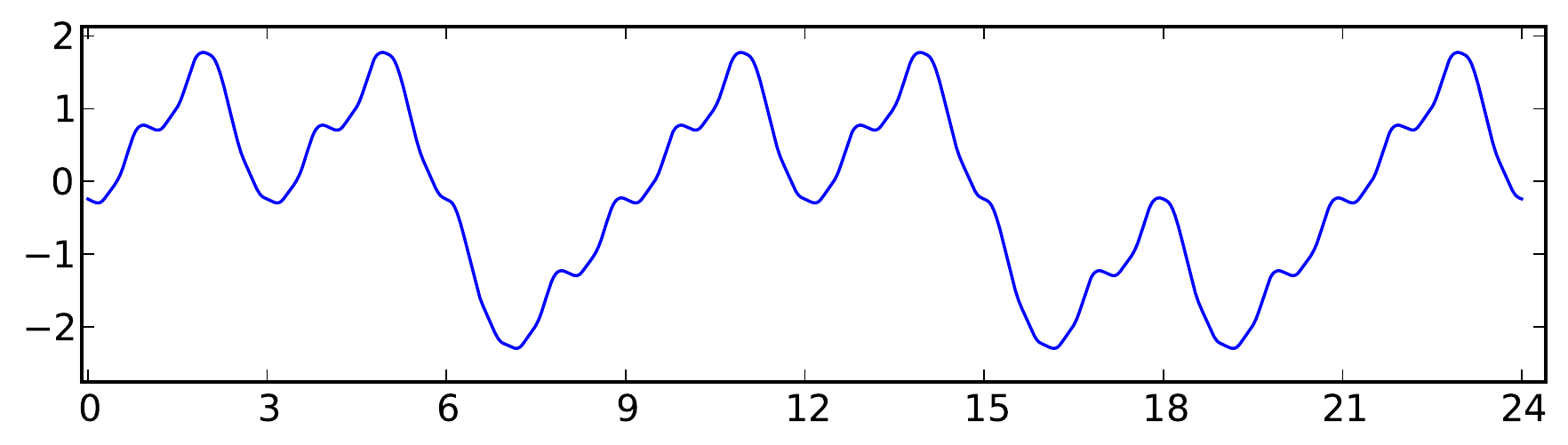} \end{center}
\caption{Fonction $f_{\w}$ pour le point fixe commen\c{c}ant par $a$ de la substitution $(a \mapsto aab,b \mapsto baa)$ 
	obtenue avec notre construction avec la fonction $\phi = \chi_a - 2 \chi_b$.}
\label{fig:exemple_f_sigma2}
\end{figure}

\begin{figure}[H]
\begin{center} \includegraphics[width=0.95\textwidth]{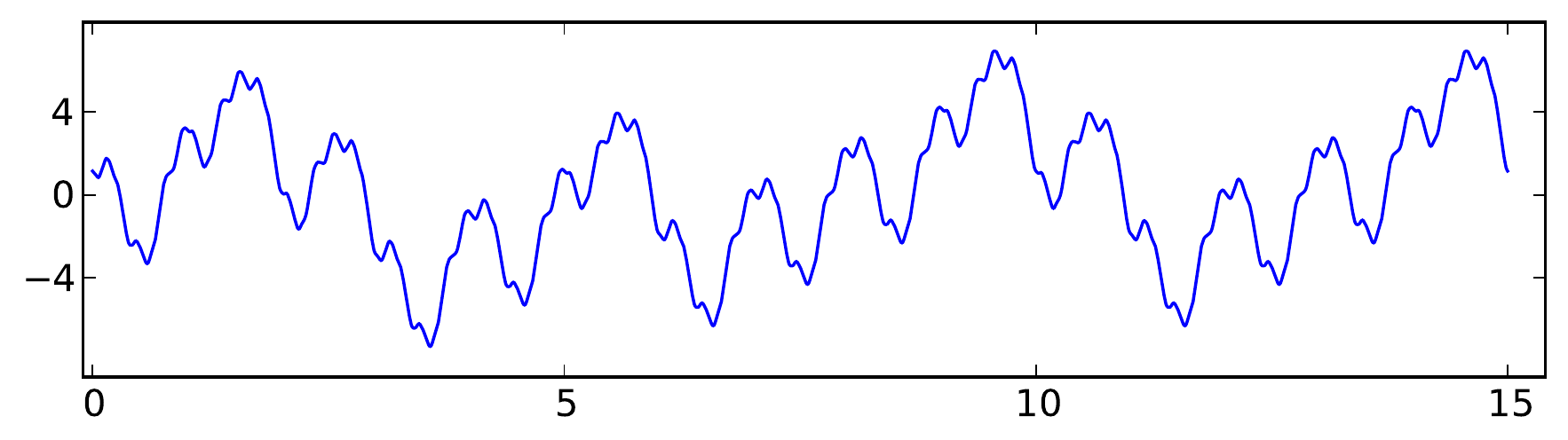} \end{center}
\caption{Fonction $f_{\w}$ pour le point fixe commen\c{c}ant par $a$ de la substitution $(a \mapsto abbaa, b \mapsto baaab)$ 
	obtenue avec notre construction avec la fonction $\phi = 2\chi_a - 3 \chi_b$.}
\label{fig:exemple_f_sigma3}
\end{figure}

Notre r\'esultat principal est le suivant.
\begin{theoreme} \label{thm:construction}
Soit $\sigma$ une substitution $\lambda$-uniforme sur $\{a,b\}$ telle que $\delta_1(\sigma) = 0$, $\delta = \delta_2(\sigma) \not= 0$ et admettant un point fixe $\w = \w_0 \w_1 \w_2 \cdots$. Alors il existe deux fonctions continues $f_a,f_b: [0,\lambda] \rightarrow \RR$ non-nulles avec $f_a(0) = f_a(\lambda) = f_b(0) = f_b(\lambda) = 0$ telles que $f = f_{\w}$ est une solution de l'\'equation~\eqref{eq:E_lambda_delta}.
\end{theoreme}
Ce th\'eor\`eme est une g\'en\'eralisation d'un r\'esultat du premier auteur~\cite{jeff2}. Dans ce dernier, la fonction limite $f_{\w}$ est construite dans le cas particulier de la substitution de Prouet-Thue-Morse pour laquelle $\delta = (0,1)$. Dans le cas g\'en\'eral, un des points d\'elicats est d'obtenir la valeur en $0$ de la fonction limite. On trouvera plusieurs autres exemples dans les figures~\ref{fig:exemple_f_sigma1}, \ref{fig:exemple_f_sigma2} et~\ref{fig:exemple_f_sigma3} pour les substitutions $(a \mapsto aab, b \mapsto aba)$, $(a \mapsto aab, b \mapsto baa)$ et $(a \mapsto abbaa, b \mapsto baaab)$ qui v\'erifient
\[
\delta(aab, aba) = (0,1,0),
\qquad
\delta(aab, baa) = (0,2,1)
\quad \text{et} \quad
\delta(abbaa, baaab) = (0, -1, 2, 3, 1),
\]
o\`u nous avons abr\'eg\'e $\delta(a \mapsto \w_a, b \mapsto \w_b)$ par $\delta(\w_a,\w_b)$. 
\'Etant donn\'ee une substitution $\lambda$-uniforme telle que $\delta_1(\sigma)=0$, 
nous verrons dans la section~\ref{se:preuve} que la fonction $\phi$ la mieux adapt\'ee \`a la construction est 
$|\sigma(a)| _b \chi_a - |\sigma(a)| _a \chi_b$.

\begin{remarque}
Il se peut qu'une substitution sur $\{a,b\}$ aient deux points fixes (un commen\c{c}ant par $a$ et un commen\c{c}ant par $b$). Dans ce cas, les fonctions $f_\w$ obtenues \`a partir de ces deux points fixes sont les concat\'enations des m\^emes fonctions $f_a$ et $f_b$, seul l'ordre de la concat\'enation diff\`ere.
\end{remarque}

\bigskip

Les premiers travaux consid\'erant de telles fonctions remontent \`a J.~Fabius~\cite{Fabius66}. 
Dans son travail, il construit une fonction qui s'\'ecrit comme concat\'enation de $f_a$ et $f_b = -f_a$ le long du mot de Prouet-Thue-Morse et qui est solution de~($E_{2,1}$). La d\'efinition de Fabius de la fonction est donn\'ee comme fonction de r\'epartition d'une somme infinie de variables al\'eatoires. Nous obtenons la m\^eme fonction par une m\'ethode bien diff\'erente dans le Th\'eor\`eme~\ref{thm:construction} en prenant la substitution de Prouet-Thue-Morse. Cette fonction est \'egalement celle consid\'er\'ee dans~\cite{jeff2}. Pour tout param\`etre $\lambda>1$ (y compris irrationnel), T.~Yoneda~\cite{Yoneda2006} et~\cite{Yoneda2007} obtient une solution de~\ref{eq:E_lambda_delta}. De plus il prouve que cette solution est born\'ee si et seulement si $\lambda \geq 2$. Pour le cas $\lambda = 2$, il retrouve la fonction de Fabius. Pour les cas $\lambda > 2$, cette solution est distincte de celles que l'on peut obtenir avec notre approche.

\bigskip

Pour d\'emontrer notre r\'esultat principal (le th\'eor\`eme~\ref{thm:construction}), nous commen\c{c}ons par introduire dans la section~\ref{se:grounil} un groupe $G$ construit comme produit semi-direct de l'espace des suites $\RR^\NN$ par $\RR$. Dans ce groupe, les sommes de Birkhoff $S^{(1)}_n$, $S^{(2)}_n$, \ldots apparaissent naturellement en regardant des produits dans $G$ (voir le lemme~\ref{lem:connexion}).

\'Etant donn\'e une substitution $\sigma$ de longueur constante, nous construisons dans la section~\ref{se:endo} un endomorphismes $\cL_\sigma$ de $G$ qui g\'en\'eralise l'ab\'elianis\'e d'une substitution : la puissance $k$-\`eme de l'endomorphisme $\cL$ d\'ecrit les sommes de Birkhoff au temps $n = \lambda^k$. Ce sont alors les puissances de $\cL$ qui vont nous permettre d'analyser le comportement des sommes. Cette analyse est le point central de cet article et est contenu dans les propositions~\ref{prop:R_et_c} et~\ref{prop:Phi}. Une fois le comportement des puissances de $\cL_\sigma$ ma\^itris\'e, nous d\'emontrons la convergence des sommes de Birkhoff dans la section~\ref{se:preuve}. \`A noter que dans cette section, les fonction $f_a$, $f_b$ et $f$ sont explicitement d\'ecrites comme des limites de somme de Birkhoff (voir le th\'eor\`eme~\ref{thm:thm_principal_precis}).

\section{Preuve des th\'eor\`emes~\ref{thm:prolong}, \ref{thm:proprietes} et~\ref{thm:sols_subs}}

Dans cette section nous donnons des preuves \'el\'ementaires des th\'eor\`emes~\ref{thm:prolong}, \ref{thm:proprietes} et~\ref{thm:sols_subs}.

\begin{proof}[Preuve du th\'eor\`eme \ref{thm:prolong}]
Soit $f$ une fonction comme dans l'\'enonc\'e.
En utilisant~\eqref{eq:Ebis_lambda_delta}, on peut prolonger $f$ sur $[\lambda,\lambda^2]$, puis $[\lambda^2,\lambda^3]$ et ainsi construire de proche en proche une fonction $f: [1,+\infty[ \rightarrow \RR$. La condition impos\'ee sur les d\'eriv\'ees garantit que la fonction est bien d\'efinie aux points $\lambda^n$ avec $n > 1$. Pour prolonger la fonction sur $[1/\lambda,1]$, on utilise la relation
\[
f(x) = f(1) - \frac{1}{\delta}\ \int_{\lambda x}^\lambda f(s) \ds.
\]
Par construction, $f$ est de classe $C^\infty$ sur $[1/\lambda,\lambda]$. Plus g\'en\'eralement, si $f$ est d\'efinie sur $[\lambda^{-n}, \lambda]$, on la prolonge sur $[\lambda^{-n-1},\lambda^{-n}]$ en posant
\[
f(x) = f(\lambda^{-n}) - \frac{1}{\delta} \int_{\lambda x}^{\lambda^{-n+1}} f(s) \ds.
\]
Aux points $\lambda^{-n}$, la fonction ainsi prolong\'ee est bien lisse (les formules qui d\'efinissent $f$ co\"incident \`a gauche et \`a droite des points $\lambda^{-n}$). Par contre, il n'est pas imm\'ediat de voir que la fonction ainsi construite se prolonge en $0$. Introduisons pour $n \geq 0$ les quantit\'es
\[
u_n = \max_{x \in [\lambda^{-n},\lambda^{-n+1}]} |f(x)| 
\quad \text{et} \quad
v_n = \max_{x,y \in [\lambda^{-n},\lambda^{-n+1}]} |f(x) - f(y)|.
\]
Nous allons montrer que $v_n$ est le terme g\'en\'eral d'une s\'erie convergente.
Par d\'efinition et l'\'equation fonctionnelle~\eqref{eq:E_lambda_delta} ces suites v\'erifient pour $n \geq 1$ :
\begin{align*}
v_n & \leq |\delta|^{-1} \int_{\lambda^{-n+1}}^{\lambda^{-n+2}} |f(s)| ds \leq |\delta|^{-1} \lambda^{-n+1} (\lambda - 1) u_{n-1} \\
u_n & \leq u_{n-1} +  v_n \\
u_{n-1} & \leq u_n +  v_{n-1}. 
\end{align*}
En particulier
\[
|u_n - u_{n-1}| \leq |\delta|^{-1} \lambda^{-n+2} \max (u_{n-1}, u_{n-2}).
\]
On peut alors utiliser le r\'esultat du lemme~\ref{lem:convergence_suite} ci-dessous.
On remarque enfin que cette construction d\'efinie de mani\`ere unique $f$ sur chaque intervalle 
	$[\lambda^{k},\lambda^{k+1}]$ pour $k\in\mathbb Z$. 
\end{proof}

\begin{lemme} \label{lem:convergence_suite}
Soit $(u_n)_{n \in \NN}$ une suite de r\'eels telle qu'il existe $\theta \in [0,1[$ tel que
\[
|u_{n+1}-u_n| \leq \theta^n \max_{0 \leq k \leq n} |u_k|.
\]
Alors la suite $u_n$ converge.
\end{lemme}

\begin{proof}[Preuve du lemme \ref{lem:convergence_suite}]
Si nous montrons que $(u_n)$ est  born\'ee alors $|u_{n+1}-u_n|$ est le terme g\'en\'eral d'une s\'erie convergente,
ce qui montrera que $(u_n)$ est une suite de Cauchy.

Pour obtenir une borne sur $u_n$ posons
\[
u'_0 = |u_0| \quad \text{et} \quad u'_{n+1} = (1 + \theta^n) u'_n.
\]
Alors
\[
|u_n| \leq u'_n \leq u'_0 \prod_{k=0}^n (1 + \theta^k) \leq u'_0 \prod_{k=0}^\infty (1 + \theta^k).
\]
Ce qui montre que $(u'_n)$ est born\'ee et conclut la preuve.
\end{proof}

\begin{proof}[Preuve du th\'eor\`eme~\ref{thm:proprietes}]
Nous commen\c{c}ons par d\'emontrer, l'annulation de $f$. Supposons que $\lambda > 1$ et $\delta > 0$. Quitte \`a remplacer $f(x)$ par $f(\delta x)$, on peut supposer $\delta=1$.

L'\'equation $(E_{\lambda,1})$ permet d'\'ecrire pour tous r\'eels $x,y$
\[
f(y) - f(x) = \int_{\lambda x}^{\lambda y} f(t) \dt.
\]

Soit $t \geq 1/(\lambda (\lambda-1))$ tel que $f(t) \not= 0$. Quitte \`a changer $f$ par $-f$ on peut supposer que $f(t) > 0$.
Supposons que $f$ ne s'annule pas sur $[t,(\lambda^3+\epsilon) t]$ pour un $\epsilon > 0$. L'\'equation $(E'_{\lambda,1})$ montre que $f$ est strictement croissante sur $[t,(\lambda^2+\epsilon/\lambda) t]$. On obtient alors que pour tout r\'eel $x \in [\lambda t, (\lambda + \epsilon/\lambda^2) t]$:
\begin{align*}
f(x) 
&> f(x) - f(t)                           && \text{car $f(t) > 0$} \\
&= \int_{\lambda t}^{\lambda x} f(s) \ds && \text{par l'\'equation fonctionnelle} \\
& \geq \int_x^{\lambda x} f(s) \ds       && \text{car $f \geq 0$ sur $[\lambda t,x] \subset [t, (\lambda^2+\epsilon/\lambda)t]$} \\
& > (\lambda x-x) f(x)                   &&  \text{car $f$ croissante sur $[x,\lambda x] \subset [t,(\lambda^3+\epsilon)t]$}.
\end{align*}
On obtient une contradiction car $x \geq \lambda t \geq 1/(\lambda-1)$.
Ainsi $f$ s'annule sur $[t, (\lambda^3+\epsilon) t]$. Par continuit\'e de $f$ et $\epsilon$ \'etant arbitraire, $f$ s'annule sur $[t,\lambda^3 t]$.

\bigskip

Passons maintenant \`a la seconde partie du th\'eor\`eme.
Soient $f$ une solution p\'eriodique et $T$ sa plus petite p\'eriode :
\[
T = \inf \{ t > 0 \ ; \ \ \forall x > 0,\, f(x) = f(x+t) \}.
\]

Alors pour tout entier $n$ :  
\[
n \int _0 ^{T} f(t)\dt = \int _0 ^{nT} f(t) \dt = \delta \left(f \left(\frac{nT}{\lambda}\right) - f(0) \right).
\]
Comme $f$ est born\'ee, en passant \`a la limite sur $n$, on obtient que $\int_0^T f(t) \dt = 0$.

Maintenant, si $x$ est un r\'eel positif on a
\begin{align*}
f(x) 
&= \frac{1}{\delta} \int _0 ^{\lambda x } f(t) \dt + f(0) 
	=\frac{1}{\delta} \int _T^{\lambda x + T } f(t) \dt  + \frac{1}{\delta} \int _0  ^{T } f(s) \ds  + f(0) \\
& =\frac{1}{\delta} \int _0  ^{\lambda x+T } f(s) \ds  +   f(0)
	= f\left(x+\frac{T}{\lambda} \right).
\end{align*}
Comme $\lambda > 1$ ceci contredit le fait que $T$ est la plus petite p\'eriode sauf si $T=0$.
\end{proof}

\begin{proof}[Preuve du th\'eor\`eme~\ref{thm:sols_subs}]
On sait par le th\'eor\`eme~\ref{thm:proprietes} que la solution $f$ ne peut pas \^etre p\'eriodique. Ainsi $f_a \not= f_b$. De plus, comme $f_a(0) = f_b(0)$, on a \'egalement $f'_a \not= f'_b$.

L'\'equation fonctionnelle peut se r\'e\'ecrire $f'(t) = \frac{\lambda}{\delta} f(\lambda t)$.
Pour tout $k \geq 0$ entier et tout $x \in [0,\lambda]$ nous avons
\[
f(\lambda^2 k + \lambda x) = \frac{\delta}{\lambda} f'(\lambda k + x)
= \left\{ \begin{array}{ll}
\dfrac{\delta}{\lambda} f'_a(x) & \text{si $\w_k = a$} \\
\dfrac{\delta}{\lambda} f'_b(x) & \text{si $\w_k = b$}.
\end{array} \right.
\]
Ainsi les valeurs de $f$ sur le segment $[\lambda^2 k, \lambda^2 k + \lambda^2]$ ne d\'ependent que de $\w_k$. Comme $f_a$ est distinct de $f_b$, et en consid\'erant des positions $k$ telles que $\w_k=a$ et $\w_k=b$, on d\'eduit qu'il existe deux mots finis $\sigma_a$ et $\sigma_b$ uniquement d\'etermin\'es de longueur $\lambda$ et tels que
\[
\dfrac{\delta}{\lambda} f'_a(x) = f_{\sigma_a}(\lambda x)
\qquad \text{et} \qquad
\dfrac{\delta}{\lambda} f'_b(x) = f_{\sigma_b}(\lambda x).
\]
On conclut alors que le mot $\w$ est point fixe de la substitution $a \mapsto \sigma_a$ et $b \mapsto \sigma_b$.
\end{proof}

\section{Un groupe nilpotent $G$ pour les sommes de Birkhoff it\'er\'ees} \label{se:grounil}

Nous introduisons le groupe $G$, produit semi-direct de $\CC$ par $\CC^\NN$. Il peut \^etre vu comme limite projective de groupes nilpotents $G_\ell$ d'indice de nilpotence $\ell$. Nous explicitons par la suite le lien entre les sommes de Birkhoff et ce groupe.

\subsection{D\'efinitions.} \label{subsee:grnil}

Soit $V = \CC^{\NN^*}$ l'ensemble des suites complexes indic\'ees \`a partir de $1$. C'est un $\CC$-espace vectoriel.
On d\'efinit l'op\'erateur de d\'ecalage $X: V \rightarrow V$ comme
\[
X \cdot (s_1, s_2, \ldots) = (0, s_1, s_2, \ldots).
\]
On note aussi $V_\ell = X^\ell \cdot V$ l'ensemble des suites qui sont nulles \`a partir du rang $\ell+1$.

L'espace vectoriel $V$ est muni de la topologie produit. Pour cette toplogie, $X$ est continue.
On dira qu'une application d'un $\CC$-espace vectoriel $W$ de dimension finie dans $V$ est \emph{polynomiale} si toutes ses coordonn\'ees sont des polyn\^omes. Une application polynomiale est continue.

Soit $A = X+I$. Cet op\'erateur agit sur les suites par : 
\[
A \cdot (s_1, s_2, s_3, \ldots) = (s_1, s_2 + s_1, s_3 + s_2, \ldots).
\]
Si on identifie $V$ aux s\'eries formelles en $X$ alors $A$ est simplement la multiplication par $1+X$.
Les puissances de $A$ sont directement donn\'ees par le bin\^ome de Newton $A^k = \sum \limits _{i=0}^k \binom{k}{i} X^i$. 
On peut \'ecrire en coordon\'ees
\[
A^k \cdot s = \left( \sum_{i=0}^k \binom{k}{i} s_{n+1-i} \right)_{n \geq 1}
\]
o\`u par convention $s_i=0$ si $i \leq 0$. Cette d\'efinition s'\'etend \`a tout nombre complexe $\xi$ en posant
\[
A^\xi = \sum_{i \geq 0} \binom{\xi}{i} X^i
\]
o\`u les coefficients binomiaux sont d\'efinis pour tout entier $k \geq 0$ par
\[
\binom{\xi}{k} = \frac{\xi(\xi-1)\cdots(\xi-k+1)}{k(k-1)\cdots 2 \cdot 1}.
\]

On note $G$ le produit semi-direct $\CC \rtimes V$ avec la multiplication
\[
\(\xi, s\) \cdot \(\zeta, t \) = \left\(\xi+\zeta, \ A^\zeta\, s + t \right\).
\]
Les sous-groupes $V_\ell$ sont distingu\'es dans $G$ et on note $G_\ell = G / V_{\ell+1}$.
On a plus pr\'ecis\'ement pour un \'el\'ement $v \in V$ et $\(z,s\) \in G$
\[
\(z,s\)^{-1} v \(z,s\) = A^z v.
\]
Le groupe $G$ est une limite projective des $G_\ell$.

Le \emph{commutateur} de deux \'el\'ements $\(\xi,s\)$ et $\(\zeta,t\)$ de $G$ est
\[
[\(\xi,s\),\(\zeta,t\)] = \(\xi,s\)^{-1} \(\zeta,t\)^{-1} \(\xi,s\) \(\zeta,t\).
\]
Nous aurons souvent besoin des calculs suivants.
\begin{lemme} \label{lem:group_law}
Soient $\xi,\zeta$ deux nombres complexes et $s,t$ deux suites de $V = \CC^{\NN^*}$. Alors
\[
\( \xi ,s \)^{-1} = \( -\xi, -A^{-\xi}  s \)
\qquad \text{et} \qquad
\left[ \(\xi,s\),\, \(\zeta,t\) \right] =
\(0, (A^\zeta - I)s - (A^\xi - I) t\).
\]
En particulier, $\(\xi,s\)$ et $\(\zeta,t\)$ commutent si et seulement si $(A^\zeta - I) s = (A^\xi - I) t$.
\end{lemme}

\begin{proof}[Preuve du lemme \ref{lem:group_law}]
Soient $\( \xi, s\)$ et $\(\zeta, t\)$ deux \'el\'ements de $G$. Alors
\[
\(\xi,s\)\cdot \(-\xi, -A^{-\xi}  s ) =\( 0, A^{-\xi} s - A^{-\xi} s \) = \(0,0\). 
\]

Pour le commutateur, un calcul direct nous donne :
\[
{[} \(\xi,s\),\(\zeta,t\) {]} 
= \(-\xi-\zeta, -A^{-\xi-\zeta}s - A^{-\zeta} t \) \(\xi+\zeta, A^\zeta s + t\) 
= \(0, A^\zeta s + t - s - A^\xi t\). \qedhere
\]
\end{proof}

\begin{remarque}
Le groupe $G_\ell$ est isomorphe au groupe multiplicatif des matrices complexes triangulaires inf\'erieures de la forme :
\[
\begin{pmatrix}
1 & 0 & \ldots \\
z & 1 & 0 & \ldots \\
\binom{z}{2} & z & 1 & 0 & \ldots \\
\vdots & \vdots & \ddots & \ddots \\
\binom{z}{\ell-1} & \binom{z}{\ell-2} & \ldots & z & 1 & 0 \\
s_\ell & s_{\ell-1} & \ldots & s_2 & s_1 & 1 \\
\end{pmatrix}
\]
En particulier, pour $\ell=2$, il s'agit du groupe de Heisenberg.

\end{remarque}

\subsection{Mots finis, mots infinis et substitutions.} \label{subsec:mots_infinis}
Soit $\cA$ un ensemble fini que l'on appelle \emph{alphabet}.
Dans la suite, on prendra tr\`es souvent $\mathcal A = \{a,b\}$.
On note $\cA^* = \cA^0 \cup \cA^1 \cup \cA^2 \cup \ldots$ l'ensemble des \emph{mots finis} sur $\cA$. Il s'agit d'un mono\"ide pour la concat\'enation et on note $\epsilon$ le mot vide. Un \emph{facteur} d'un mot $\u = \u_0 \u_1 \ldots \u_{n-1}$ est un mot $\v$ tel que $\v = \u_i \u_{i+1} \ldots \u_{i+k-1}$. On consid\'erera \'egalement $\mathcal A^\NN$ l'ensemble des \emph{mots infinis} et $\cA^\ZZ$ l'ensemble des \emph{mots bi-infinis}. On munit $\cA^\NN$ et $\cA^\ZZ$ de la topologie produit ce qui en fait des espaces compacts.

Une \emph{substitution} sur $\cA$ est un morphisme de $\cA^*$ vu comme mono\"ide.
On dira qu'une substitution est \emph{positive} si chaque image de lettre contient toutes les lettres de l'alphabet. Par exemple la substitution de Prouet-Thue-Morse $a \mapsto ab, b \mapsto ba$ est positive. Le \emph{langage} d'une substitution positive est le plus petit sous-ensemble de $\cA^*$ contenant $\cA$, stable par facteur et par $\sigma$. C'est aussi l'ensemble des facteurs des mots $\sigma^n(\alpha)$ o\`u $\alpha$ est une lettre de l'alphabet et $n$ un entier. Par exemple, le langage de la substitution de Prouet-Thue-Morse est
\[
\{ \ \epsilon, \ a,  \  b,  \ aa,  \ ab, \  ba,  \ bb,  \ aab,  \ aba,  \ abb,  \ baa,  \ bab,  \ bba,  \ \ldots \   \}.
\]

\`A une substitution positive $\sigma$, on associe \'egalement le sous-ensemble de $\cA^\NN$ form\'e des mots infinis dont tous les facteurs sont dans le langage de $\sigma$. C'est un ensemble compact et non-vide que l'on note $K_{\sigma,\NN}$. De la m\^eme fa\c{c}on, on associe un sous-espace de $\cA^\ZZ$ que l'on notera $K_{\sigma,\ZZ}$. La substitution $\sigma$ agit continuement sur $K_{\sigma,\NN}$ et $K_{\sigma,\ZZ}$. Un point p\'eriodique pour $\sigma$ de p\'eriode $p$ est un \'el\'ement de $K_{\sigma,\NN}$ ou $K_{\sigma,\ZZ}$ tel que $\sigma^p(\w) = \w$. Par exemple, le mot de Prouet-Thue-Morse commen\c{c}ant par
\[
\w = abbabaabbaababbabaababbaabbabaabbaababbaabbabaababbabaabbaababbabaaba \ldots
\]
est un point fixe de $(a \mapsto ab, b \mapsto ba)$ (i.e. un mot de p\'eriode 1).

\medskip

Dans la suite, on note indiff\'eremment $K$ l'ensemble des mots infinis ou bi-infinis sur $\cA$.
Le d\'ecalage sur $K$ est l'application qui consiste \`a d\'ecaler l'origine $T: (\u_i)_i \mapsto (\u_{i+1})_i$.
Si l'on munit $X$ de la topologie produit, alors $X$ est compact et $T$ est une application continue.
Sur $\cA^\ZZ$ c'est un hom\'eomorphisme alors que sur $\cA^\NN$ c'est une application dont chaque fibre a le m\^eme cardinal que $\cA$.
L'\emph{orbite} d'un point $\u \in K$ est la suite de points $\u$, $T\u$, $T^2\u$, \ldots.

Si $\sigma$ est une substitution positive, alors $K_{\sigma,\NN}$  et $K_{\sigma,\ZZ}$ sont invariants par $T$: ce sont des \emph{sous-d\'ecalages} de $\cA^\NN$ et $\cA^\ZZ$ respectivement. Ils ont la propri\'et\'e d'\^etre minimaux: pour tout mot $\w \in K_{\sigma,\ZZ}$ son orbite positive $\{T^n \w\}_{n \in \NN}$ est dense dans $K_{\sigma,\ZZ}$. 



\subsection{Sommes de Birkhoff.}  \label{subsec:sommes_iterees} \label{se:birkhoff}
Soit $\cA$ un alphabet. 
On note $K$ le d\'ecalage $\cA^\NN$ ou $\cA^\ZZ$.
Soit $\phi: K \rightarrow \RR$ une application continue. Par exemple, la fonction $\chi_a: K \rightarrow \RR$ d\'efinie par
\[
\chi_a(\u) = \left\{ \begin{array}{ll}
1 & \text{si $\u_0 = a$} \\
0 & \text{sinon.}
\end{array} \right.
\]

Soit $\u$ un \'el\'ement de $K$. Nous associons \`a la paire $(\phi,\u)$ sa \emph{somme de Birkhoff} d\'efinie par
\[
S^{(1)}_n(\phi,\u) = \sum_{k=0}^{n-1} \phi(T^k \u).
\]
Si $\phi$ est la fonction $\chi_a$, alors $S^{(1)}_n(\phi,\u)$ est le nombre de $a$ dans le pr\'efixe de $\u_0 \u_1 \ldots \u_{n-1}$ de $\u$.
On d\'efinit par r\'ecurrence les \emph{sommes de Birkhoff it\'er\'ees} de $\phi$ par
\begin{equation} \label{eq:sums} 
S_n^{(\ell+1)}(\phi,\u) = \sum_{k=0}^{n-1} S^{(\ell)}_k(\phi,\u).
\end{equation}

Soit $\sigma$ une substitution positive et $K_{\sigma,\NN}$ le d\'ecalage associ\'e.
Le th\'eor\`eme suivant est un r\'esultat classique de th\'eorie ergodique.
Nous renvoyons \`a \cite{Py} pour des d\'etails.
\begin{proposition}[Unique ergodicit\'e]
Soit $\sigma$ une substitution positive et $(K,T)$ le d\'ecalage associ\'e.
Alors il existe une unique mesure de probabilit\'e $\mu$ sur $K$ invariante par $T$.
De plus, pour toute fonction $\phi \in C(K)$, pour tout $\ell \geq 1$, on a uniform\'ement en $\u$
\[
\lim_{n \to \infty} \binom{n}{\ell}^{-1} S_n^{(\ell)}(\phi,\u) = \int_K \phi d\mu.
\]
\end{proposition}
L'\'enonc\'e ci-dessus est souvent \'ecrit avec $\ell=1$. La convergence pour $\ell \geq 2$ d\'ecoule directement du cas $\ell=1$ en remarquant que les moyennes au rang $\ell+1$ sont les moyennes de Ces\`aro ce celles au rang $\ell$.

Ce r\'esultat montre que la premi\`ere \'echelle d'approximation de $S^{(\ell)}_n(\phi,\u)$ est $\binom{n}{\ell} \int_K \phi$. Dans cet article, pour certaines substitutions et certaines fonctions, nous montrons que quitte \`a soustraire un polyn\^ome, les sommes de Birkhoff restent born\'ees. Ce r\'esultat est tr\`es sp\'ecifique et ne concerne pas toutes les substitutions. On trouvera une \'etude pr\'ecise de la croissance des sommes de Birkhoff $S^{(1)}_n$ pour une substitution quelconque dans~\cite{adam1}.

\bigskip

Nous faisons maintenant le lien avec l'op\'erateur $A$ et le groupe $G$ introduits dans la section~\ref{subsee:grnil}.
Pour chaque $\ell$ et chaque $n$, on peut voir $S_n^{(\ell)}$ comme un endomorphisme de $C(K)$.
On note $S_n = (S_n^{(1)}, S_n^{(2)}, S_n^{(3)}, \ldots)$ vu comme op\'erateur de $C(K)$ dans $C(K)^{\NN^*}$.
On plonge $C(K)$ dans $C(K)^{\NN^*}$ via $f \mapsto (f,0,0,\ldots)$.
On d\'efinit deux op\'erateurs $U_T$ et $A$ sur $C(K)^{\NN^*}$ comme
\[
U_T \cdot (\phi_1, \phi_2, \ldots) = (\phi_1 \circ T, \phi_2 \circ T, \ldots)
\quad \text{et} \quad
A \cdot (\phi_1, \phi_2, \ldots) = (\phi_1, \phi_2 + \phi_1, \phi_3 + \phi_2, \ldots).
\]
La notation $A$ surcharge la d\'efinition que nous avons donn\'ee sur $V$ mais on a la compatibilit\'e suivante. Pour tout mot $\u$ dans $K$ on a une projection $e_\u: C(K)^{\NN^*} \rightarrow V$ qui consiste \`a \'evaluer toutes les coordonn\'ees en $\u$. On a alors pour tout $\phi$, $e_\u (A \phi) = A (e_\u(\phi))$. Notons que les op\'erateurs $U_T$ et $A$ commutent.

\'Etant donn\'e $\phi = (\phi_1, \phi_2, \ldots) \in C(K)^{\NN^*}$ on peut lui associer ses sommes de Birkhoff en posant
\[
S_0 = I
\quad \text{et} \quad
S_{n+1} = U_T^n + A S_n.
\]
Bien entendu, si $\phi \in C(X)$ alors on retrouve la d\'efinition standard via le plongement $C(K) \to C(K)^{\NN^*}$
\[
S_n(\phi,\u) = (S^{(1)}_n(\phi,\u),\, S^{(2)}_n(\phi,\u),\, \ldots).
\]
\begin{lemme} \label{lem:bs_explicite}
On a
\[
S_n = \sum_{k=0}^{n-1} A^{n-k-1} U_T^k.
\]
En particulier, pour une fonction $\phi \in C(K)$ et un entier $\ell \geq 0$ on a
\[
S_n^{(\ell)}(\phi, \u) = \sum_{k=0}^{n-\ell} \binom{n-k-1}{\ell-1} \phi(T^k \u).
\]
\end{lemme}
Du fait de leur d\'efinition, les sommes de Birkhoff $S_n$ peuvent \^etre construites directement dans le groupe $G$.
\'Etant donn\'e $\phi \in C(K)$, on d\'efinit $\pi_\phi: K \rightarrow G$ par
\[
\pi_\phi(\u) = \(1, (\phi(\u), 0, 0, \ldots) \).
\]
\begin{lemme} \label{lem:connexion}
Soit $\u \in \cA^\NN$ et $\phi \in C(K)$ alors
\[
\pi_\phi(\u) \pi_\phi(T \u) \ldots \pi_\phi(T^{n-1} \u) = \( n, S_n(\phi,\u)\) = \(n, (S^{(1)}_n(\phi,\u), S^{(2)}_n(\phi,\u), \ldots)\).
\]
\end{lemme}

On d\'emontre maintenant le lemme \'enonc\'e dans l'introduction qui dit que les quantit\'es $\delta_i$ sont bien d\'efinies.
\begin{proof}[Preuve du lemme~\ref{le:def_delta}]
Pour tous r\'eels $\alpha,\beta$ on a $S_n^{(i)}(\alpha \phi + \beta, \u) = \alpha S_n^{(i)}(\phi, \u) + \binom{n}{i} \beta$. On en d\'eduit donc que la quantit\'e $(S^{(i)}_\lambda(\phi,\sigma(a)) - S^{(i)}_\lambda(\phi,\sigma(b)))/ (\phi(a)-\phi(b))$ est invariante par toute transformation affine $\phi \mapsto \alpha \phi + \beta$. Comme toutes les fonctions $\phi: \{a,b\} \rightarrow \RR$ telles que $\phi(a) \not= \phi(b)$ sont reli\'ees entre elle par une relation affine, on en d\'eduit le r\'esultat.
\end{proof}

\begin{remarque}
Nous pouvons consid\'erer la double suite $(S^{(\ell)}_n(\phi,\u))_{\ell,n}$ comme un Çtriangle de Pascal g\'en\'eralis\'eÈ,
 pour lequel la premi\`ere colonne, constitu\'ee de 1 dans le triangle de Pascal, est remplac\'ee par $u_n = \phi(T^n \u)$ et dans lequel on applique la m\^eme r\`egle de construction :

\begin{center}
\begin{tabular}{|l||c|c|c|c|c|c|}
\hline
$n \backslash \ell$ & 0 & 1 & 2 & 3 & 4 & 5 \\
\hline \hline
0 & $u_0$ & 0               & 0          & 0  & 0 & 0 \\
\hline
1 & $u_1$ & $u_0$             & 0          & 0 & 0 & 0 \\
\hline
2 & $u_2$ & $u_0 + u_1$       & $u_0$        & 0  & 0 & 0 \\
\hline
3 & $u_3$ & $u_0 + u_1 + u_2$ & $2u_0 + u_1$ & $u_0$ & 0 & 0\\
\hline
4 & $u_4$ & $u_0 + u_1 + u_2 + u_3$ & $3u_0 + 2u_1 + u_2$ & $3u_0 + u_1$ & $u_0$ & 0 \\
\hline
5 & $u_5$ & $u_0 + u_1 + u_2 + u_3 + u_4$ & $4u_0 + 3u_1 + 2u_2 + u_3$ & $6u_0 + 3u_1 + u_2$ & $4u_0+u_1$ & $u_0$\\
\hline
\end{tabular}
\end{center}
\end{remarque}

%

\subsection{Puissances dans $G$.} \label{subse:pui}

Nous allons d\'efinir la puissance $\xi$-i\`eme d'un \'element $\(z,s)$ de $G$ pour tout nombre complexe $\xi$.
Notons d\'ej\`a que pour toute puissance enti\`ere $k$ :
\[
\(1,s\)^k = \left\( \xi, (I+A+\cdots+A^{k-1}) \, s \right\).
\]
Posons pour tout entier $k \geq 0$, 
\[
\displaystyle B(k) = I+A+\cdots+A^{k-1} = \sum_{n \geq 0} \binom{k}{n+1} X^n.
\]
Comme pour les puissances de $A$, on \'etend cette d\'efinition aux nombres complexes $\xi$ non nuls en posant
\[
B(\xi) = \sum_{n \geq 0} \binom{\xi}{n+1} X^n = \frac{A^\xi - 1}{A - 1} = \frac{A^\xi - 1}{X}.
\]
Nous d\'efinissons pour tout complexe $\xi$ la puissance $\xi$-i\`eme d'un \'el\'ement $\(z,s\)$ de $G$ en posant
\[
\(z,s\)^\xi = \left\{ \begin{array}{ll}
\(z\, \xi, B(z \xi) \, B(z)^{-1} \cdot s\) = \left\(z\, \xi, \frac{A^{z \xi} - 1}{A^z - 1} \cdot s \right\) & \text{si $z \not= 0$} \\
\(0, \xi \, s\) & \text{si $z = 0$}
\end{array} \right. .
\]

\begin{lemme} \label{le:puis}
L'application $(\(z,s\), \xi) \mapsto \(z,s\)^\xi$ est continue. De plus, pour tout \'el\'ement $\(z,s\)$ de $G$ et tous nombres complexes $\xi$ et $\zeta$ :
\[
\(z,s\)^{\xi+\zeta} = \(z,s\)^\xi \cdot \(z,s\) ^\zeta.
\]
Soit $\(z,s\)$ un \'el\'ement de $G$ avec $z$ non nul, alors l'ensemble des \'el\'ements de $G$ qui commutent avec $\(z,s\)$ sont les \'el\'ements de la forme $\(z,s\)^\xi$ avec $\xi \in \RR$. 
\end{lemme}

\begin{proof}[Preuve du lemme \ref{le:puis}]
Fixons un \'element $\(z, s\)$ de $G_\ell$.
L'application $\xi \mapsto \(z, s\)^{\xi}$ est continue de $\RR^*$ dans $G$.
Nous allons montrer qu'elle se prolonge en $0$ par la formule donn\'ee dans le lemme.

Il nous suffit de montrer que pour tous les nombres complexes $z$ et $\xi$
\begin{equation} \label{eq:puipui}
B( z \xi) \, B(z)^{-1} \underset{m \to 0}{\longrightarrow} \xi.
\end{equation}
Pour cela, il suffit de remarquer que
\[
\frac{A^{z \xi} - 1}{z} \underset{z \to 0}{\longrightarrow} \xi.
\]

Nous allons maintenant montrer la relation sur les puissances. Soit $\xi,\zeta \in \RR$ et $\(z,s\) \in G$. Nous avons
\begin{align*}
\(z, s\)^\xi \cdot (z, s\)^\zeta
&= \(z\xi, B(z\xi) B(z)^{-1} s \) \, \(z \zeta, B(z \zeta) B(z)^{-1} s \) \\
&= \left\(z (\xi + \zeta), \left(A^{z\zeta} B(z\xi) + B(z \zeta) \right) B(z)^{-1} s \right\).
\end{align*}
Le r\'esultat suit en remarquant que pour tous les nombres complexes $\xi$ et $\zeta$ on a 
\[
B(\xi+\zeta) = B(\xi) + A^\xi B(\zeta).
\]

D'apr\`es le lemme~\ref{lem:group_law}, deux \'el\'ements $(z,s)$ et $(z',t)$ de $G$ commutent si et seulement si $(A^{z'} - 1) s = (A^z - 1) t$. Cette relation peut \'egalement s'\'ecrire $B(z') s = B(z) t$ et si $z \not= 0$ on a $t = B(z') B(z)^{-1} s$. Autrement dit $\(z,s\)^{z'/z} = \(z',t\)$.
\end{proof}

%
%
%
%
%
%
\subsection{Sous-groupes discrets.} \label{subse:ssdiscrets}

Notons $\Gamma$ (respectivement $\Gamma_\ell$) les \'el\'ements de $G$ (resp. $G_\ell$) dont tous les coefficients sont entiers. 
Autrement dit, $\Gamma = \ZZ \rtimes \ZZ^\NN$ (resp. $\Gamma_\ell = \ZZ \rtimes \ZZ^\ell$).

Soit $\bold{b} = \(1,0\)$, et pour $i \geq 1$, $\bold{a}_i  = \(0, X^{i-1} \cdot 1 \)$.
Remarquons que
\begin{align*}
[\bold{a}_i, \bold{b}] &= 
\(0,-X^{i-1}\)\, \(-1,0\)\, \(0,X^{i-1} \cdot 1 \)\, \(1,0\) \\ &=
\(-1,-(1+X)^{-1} X^{i-1} \cdot 1 \)\, \(1, (1+X)X^{i-1} \cdot 1\) \\ &=
\(0, X^i \cdot 1 \) = \bold{a}_{i+1}.
\end{align*}

\begin{lemme} \label{le;puisdef}
Avec les m{\^e}mes notations que ci-dessus. Le groupe $\Gamma$ est isomorphe au groupe engendr\'e par $\bold{b},\bold{a}_1,\bold{a}_2,\ldots$ avec les relations
\begin{enumerate}
\item $\bold{a}_{n+1} = [\bold{a}_n,\bold{b}]$ pour $n$,
\item $[\bold{a}_i,\bold{a}_j] = 1$ pour tout $i,j$.
\end{enumerate}
Pour obtenir le groupe $\Gamma_\ell$, il suffit d'ajouter la relation $\bold{a}_{\ell+1} = 1$.
\end{lemme}

\begin{proof}
Soit $H$ le groupe engendr\'e par $\bold{b},\bold{a}_1,\bold{a}_2,\ldots$ et les relations donn\'ees dans l'\'enonc\'e.
On construit un morphisme $f:\Gamma \rightarrow H$ en envoyant l'\'el\'ement $\(1,0\)$ sur $\bold{b}$ et $\(0,X^{i-1} \cdot 1\)$ sur $\bold{a}_i$.
Il est par construction injectif car $\Gamma$ v\'erifie les relations donn\'ees dans le lemme.

Le morphisme inverse est construit de la mani\`ere suivante.
En utilisant la  relation $\bold{a}_{n+1} = [\bold{a}_n,\bold{b}]$ qui se r\'e\'ecrit $\bold{a}_n \bold{b} = \bold{b} \bold{a}_n \bold{a}_{n+1}$ et la commutation des \'el\'ements $\bold{a}_i$, tout \'el\'ement de $H_\ell$ peut se mettre sous la forme
$\bold{b}^m \bold{a}_1^{s_1} \bold{a}_2^{s_2} \bold{a}_3^{s_3} \ldots \bold{a}_\ell^{s_\ell}$
	o\`u $s_1, s_2, \ldots, s_\ell, m \in \ZZ$.
Cet \'el\'ement est bien s\^ur l'image de $\(m, (s_1,\ldots,s_\ell)\)$ ce qui montre que $f$ est surjectif.
\end{proof}

\begin{remarque}
Le groupe $\Gamma_\ell$ est engendr\'e par $\bold{b}_1$ et $\bold{a}$. Pour $\ell=2$, on obtient assez simplement la liste de relations
\[
[\bold{a}_1,[\bold{a}_1,\bold{b}]] = 1 \quad et \quad  [\bold{b},[\bold{a}_1,\bold{b}]] = 1
\]
qui sont les relations standards pour le groupe d'Heisenberg \`a coefficients entiers.
\end{remarque}

\section{Endomorphismes de $G$ et substitutions} \label{se:endo}

On se donne une substitution $\sigma$ et une fonction $\phi: \{a,b\} \rightarrow \RR$ non constante. Le but de cette section est de d\'efinir un endomorphisme $\cL = \cL_{\phi,\sigma}$ qui permettent de traduire l'action de $\sigma$ au niveau des sommes de Birkhoff. Plus pr\'ecis\'ement, nous construisons tout d'abord dans la section~\ref{subse:projection} un morphisme $\pi_\phi: \{a,b\}^* \rightarrow G$ (semblable \`a celui de la section \ref{subsec:sommes_iterees}) tel que pour tout mot fini $\u$
\[
\pi_\phi(\u) = \(|u|, S_{|u|}(\phi,\u)\).
\]
Nous construisons ensuite un endomorphisme $\cL_{\phi,\sigma}$ de $G$ tel que $\pi_\phi \circ \sigma = \cL_{\phi,\sigma} \circ \pi_\phi$.

\subsection{Projections du mono{\"\i}de libre dans $G$.} \label{subse:projection}
Soit $\phi: \{a,b\} \rightarrow \RR$ une fonction. On consid\`ere implicitement $\phi$ comme une fonction sur les mots infinis en posant $\phi(\u) = \phi(\u_0)$. \`A un mot fini $\u$ du mono\"ide libre g\'en\'er\'e par $a$ et $b$, nous pouvons lui associer ses sommes de Birkhoff it\'er\'ees (voir le lemme~\ref{lem:connexion}):
\[
\pi_\phi(\mathsf{u}) = \(1, \phi(\mathsf{u}_0))\) \(1, \phi(\mathsf{u}_1)\) \ldots \(1, \phi(\mathsf{u}_{|u|-1})\) = \(|\mathsf{u}|, S_{|\mathsf{u}|}(\phi,\mathsf{u}) \).
\]
L'application $\pi_\phi$ ainsi d\'efinie est un morphisme: $\pi_\phi(\u \v) = \pi_\phi(\u) \pi_\phi(\v)$.

En prenant $\phi = \chi_a$ on a par exemple
\[
\begin{array}{|l|l|}
\hline
bbbb & (4, (0, 0, 0, 0)) \\
\hline
bbba & (4, (1, 0, 0, 0)) \\
bbab & (4, (1, 1, 0, 0)) \\
babb & (4, (1, 2, 1, 0)) \\
abbb & (4, (1, 3, 3, 1)) \\
\hline
\end{array} \hspace{.5cm} \begin{array}{|l|l|}
\hline
bbaa & (4, (2, 1, 0, 0)) \\
baba & (4, (2, 2, 1, 0)) \\
baab & (4, (2, 3, 1, 0)) \\
abba & (4, (2, 3, 3, 1)) \\
abab & (4, (2, 4, 3, 1)) \\
aabb & (4, (2, 5, 4, 1)) \\
\hline
\end{array} \hspace{.5cm} \begin{array}{|l|l|}
\hline
baaa & (4, (3, 3, 1, 0)) \\
abaa & (4, (3, 4, 3, 1)) \\
aaba & (4, (3, 5, 4, 1)) \\
aaab & (4, (3, 6, 4, 1)) \\
\hline
aaaa & (4, (4, 6, 4, 1)) \\
\hline
\end{array}
\]

\begin{proposition} \label{prop:injection_monoide_dans_G}
Soit $\phi: \{a,b\} \rightarrow \RR$ telle que $\phi(a) \not= \phi(b)$. Alors l'application $\pi_\phi:\{a,b\}^* \rightarrow G$ est injective. De mani\`ere plus pr\'ecise, $\pi_\phi$ est une injection de l'ensemble des mots de longueurs $\ell$ dans $G_\ell$.
\end{proposition}

\begin{proof}
Si $\u$ et $\v$ n'ont pas la m\^eme longueur alors clairement $\pi_\phi(\u) \not= \pi_\phi(\v)$. On suppose alors que $\u$ et $\v$ ont m\^eme longueur et on proc\`ede par r\'ecurrence sur la longueur.

Le r{\'e}sultat est imm{\'e}diat pour les mots $a$ et $b$ puisque nous avons suppos\'e que $\phi(a) \neq \phi(b)$.

Supposons le r{\'e}sultat vrai pour les mots de longueurs plus petites ou \'egales \`a $n-1$.
Fixons alors deux mots $\u = \u_0 \ldots \u_{n-1}$ et $\v = \v_0 \ldots \v_{n-1}$ de longueur $n$ tels que $\pi_\phi(\u)=\pi_\phi(\v)$.

On a $S^{(n)}_n (\phi,\u) =\phi(\u_0)$ et $S^{(n)}_n(\phi,\v) =\phi(\v_0)$.
Donc n{\'e}cessairement, les deux mots commencent par la m{\^e}me lettre.
Notons $\u' = \u_1 \ldots \u_{n-1}$ et $\v' = \v_1 \ldots \v_{n-1}$.
Par d\'efinition $\pi_\phi(\u) = \pi_\phi(\u_0) \pi_\phi(\u')$ et donc l'\'egalit\'e $\pi_\phi(\u) = \pi_\phi(\v)$ entra\^ine que $\pi_\phi(\u') = \pi_\phi(\v')$. On peut alors appliquer l'hypoth\`ese de r\'ecurrence pour d\'eduire que $\u' = \v'$ et donc que $\u = \v$.
\end{proof}

\subsection{G\'en\'eralit\'es sur les endomorphismes de $G$.}
Dans cette section, nous classifions les endomorphismes des groupes $G_\ell$ et du groupe $G$ introduit dans la section~\ref{subsee:grnil}. Rappelons que $G$ est construit \`a partir de l'espace vectoriel des suites $V = \CC^{\NN^*}$ et de l'op\'erateur $A = 1+X$.

Un \emph{endomorphisme} de $G$ est une application $\cL: G \rightarrow G$ telle que
\[
\cL(\(z,s\) \(z',t\)) = \cL(\(z,s\)) \cdot \cL(\(z',t\)) \quad et \quad \cL(\bold{0})=\bold{0}.
\]

Il est alors clair d'apr\`es le lemme \ref{le;puisdef} que pour tout r\'eel $\xi$, $\cL(\(z,s\)^\xi)=\cL(\(z,s\))^\xi$.
\bigskip

Notons tout d'abord que l'endomorphisme int\'erieur $\cL(x) = \(z,s\)^{-1} x \(z,s\)$ v\'erifie
\[
\cL(\(1, 0\)) = \(1,(I-A) s\)
\quad \text{et} \quad
\cL(\(0,1\)) = \(0, A^z\).
\]

\begin{lemme} \label{lem:endo_Gl}
Soit $\cL: G \rightarrow G$ un endomorphisme de $G$.
Les quantit\'es $\lambda,\mu \in \RR$ et $\beta,\delta \in V$ telles que
\[
\cL( \(1,0\) ) = \(\lambda, \beta \)
\quad \text{et} \quad
\cL( \(0,1\) ) = \(\mu, \delta \)
\]
v\'erifient l'\'equation
\begin{equation} \label{eq:condition_endomorphisme}
(A^\mu - 1)\Big((A^\lambda - 1) \delta - (A^\mu - 1) \beta\Big) = 0.
\end{equation}
Autrement dit, si $\mu=0$ il n'y a pas de condition, sinon il faut que $(A^\lambda - 1) \delta = (A^\mu - 1) \beta$.

R\'eciproquement, \'etant donn\'e $\lambda,\nu \in \RR$ et $\beta,\delta \in V$ v\'erifiant~\eqref{eq:condition_endomorphisme} il existe un unique endomorphisme $\cL:G \to G$ tel que $\cL \(1,0\) = \(\lambda,\beta\)$ et $\cL \(0,1\) = \(\mu,\delta\)$.
\end{lemme}
En passant au quotient, on retrouve le cas $\ell=2$ de~\cite{MR1326950}.



\begin{proof}[Preuve du lemme \ref{lem:endo_Gl}]
On sait que $G$ est engendr\'e par les deux \'el\'ements $\bold{b}=\(1,0\)$ et $\bold{a}_1=\(0,1\)$ et leurs puissances. Donc leurs images par $\cL$ d\'eterminent enti\`erement $\cL$. Nous allons \'etablir quelles relations doivent v\'erifier leurs images. Pour cela, nous commen{\c c}ons par calculer l'image de $\bold{a}_n$ pour $n \geq 2$.

D'apr\`es le lemme~\ref{lem:group_law}, le commutateur d'un \'el\'ement $\(r,w\)$ avec $\(\lambda,\beta\)$ est
\[
[\(\lambda,\beta\),\(r,w\)] = \(0, (A^r - 1) \beta - (A^\lambda-1) w\).
\]
On en d\'eduit $\cL(\bold{a}_1) = \(\mu,\delta\)$, $\cL(\bold{a}_2) = \left\(0, (A^\lambda - 1) \delta - (A^\mu - 1) \beta \right\)$ et
\[
\cL(\bold{a}_n) = \left\(0, (A^\lambda - 1)^{n-2} ((A^\lambda - 1) \delta - (A^\mu - 1) \beta) \right\).
\]
Il est alors clair que les images de $\bold{a}_n$ commutent pour $n \geq 2$. 
La seule obstruction est donc  $\cL([\bold{a}_1,\bold{a}_2]) = 1$ et
un calcul explicite nous donne
\begin{align*}
[\cL(\bold{a}_2),\cL(\bold{a}_1)] = \(0, (A^\mu - 1)((A^\lambda - 1) \delta - (A^\mu - 1)\beta\).
\end{align*}
D'o\`u le r\'esultat.
\qedhere
\end{proof}

\subsection{Endomorphismes associ\'es aux sommes de Birkhoff.}
\label{subse:induit}
Rappelons que nous voulons analyser les sommes de Birkhoff it\'er\'ees d'une fonction $\phi: \{a,b\} \to \RR$. 
Nous avons vu en~\ref{subsec:sommes_iterees} que ces sommes de Birkhoff it\'er\'ees peuvent \^etre vues comme une projection du mono\"ide libre $\{a,b\}^*$ dans le groupe $G$.

\begin{lemme}\label{lem:L_f_sigma}
Soit $\sigma$ une substitution de longueur constante $\lambda > 0$ sur $\{a,b\}$ telle que $\sigma(a) \not= \sigma(b)$.
Soit $\phi: \{a,b\} \rightarrow \RR$ une fonction telle que $\phi(a) \not= \phi(b)$.
Alors il existe un unique morphisme $\cL = \cL_{\phi,\sigma}$ tel que
\begin{equation} \label{eq:mopi}
\cL \circ \pi_\phi = \pi_\phi \circ \sigma .
\end{equation}
De plus, ce morphisme pr\'es\`erve $V$ et plus pr\'ecis\'ement, pour tout $n \geq 1$ on a :
\[
\cL (\bold{a}_n) = \left\( 0, (A^\lambda - I)^{n-1} \delta \right\),
\]
o\`u $\delta = \delta(\sigma) = (\delta_1(\sigma), \delta_2(\sigma), \ldots)$ est le vecteur d\'efini dans le lemme~\ref{le:def_delta}.
\end{lemme}
Remarquons que l'action de $\cL$ sur $V$ ne d\'epend pas de la fonction $\phi$ choisie.

\begin{proof}
D'apr\`es le lemme~\ref{lem:connexion}, il faut et il suffit de construire un morphisme tel que
\[
\cL \big( \(1, \phi(a)\) \big)= \(\lambda, S(\phi,\sigma(a))\)
\quad \text{et} \quad
\cL \big( \(1, \phi(b)\)\big) = \(\lambda, S(\phi,\sigma(b))\).
\]
Pour ce faire, on se ram\`ene au lemme~\ref{lem:endo_Gl}. On a, 
\[
\(1,\phi(a)\)^{-1} \(1,\phi(b)\) = \(0,\phi(b)-\phi(a)\) = \(0,1\)^{\phi(b)-\phi(a)}.
\]
Donc 
\[
\(0,1\) = (\(1,\phi(a)\)^{-1} \(1,\phi(b)\))^{1/(\phi(b)-\phi(a))}.
\]

De la m\^eme fa\c{c}on on peut exprimer $\(1,0\)$ \`a partir de $\(1,\phi(a)\)$ et $\(1,\phi(b)\)$ en utilisant $\(1,0\) = \(1,\phi(a)\) \(0,1\)^{-\phi(a)}$. On peut alors utiliser le lemme~\ref{lem:endo_Gl} en remarquant qu'on a un exemple dans lequel $\mu=0$. En particulier, $\cL (V) \subset V$. L'unicit\'e vient de la construction car les \'el\'ements $\(1,\phi(a)\)$ et $\(1,\phi(b)\)$ et leurs puissances engendrent tout $G$.

La formule explicite sur $V$ s'obtient d'une part en remarquant que 
\[
\cL \Big( \big\(0, \phi(a) - \phi(b) \big\) \Big) = \big\(0, S(\phi, \sigma(a)) - S(\phi,\sigma(b)\big\) = \big\(0, (\phi(a) - \phi(b)) \delta\big\)
\]
et en suivant les calculs du lemme~\ref{lem:endo_Gl}.
\end{proof}

Soit $\lambda > 0$, $\beta,\delta \in V$. On consid\`ere le morphisme $\cL = \cL_{\lambda,\beta,\delta}$ qui v\'erifie
\[
\cL \big( \(1,0\)\big) = \(\lambda,\beta\)
\quad \text{et} \quad
\cL \big( \(0,1\)\big) = \(0,\delta\).
\]
Par la remarque~\ref{lem:L_f_sigma} on a $\cL (\mathbf{a}_n) = \(0, (A^\lambda - 1)^{n-1} \delta \)$. On obtient alors pour un \'el\'ement $\(z,s\)$ de $G$ l'expression suivante
\begin{align*}
\cL \big(\(z, s\) \big)
&= \cL (\mathbf{b}^z \mathbf{a}_1^{s_1} \ldots \mathbf{a}_\ell^{s_\ell} \ldots) \\
&= \(\lambda,\beta\)^z \ \prod_{k=1}^\infty \(0, (A^\lambda - 1)^{k-1} \delta\)^{s_k} \\
&= \left\(\lambda z, \left( \sum_{k=1}^\infty s_k (A^\lambda - 1)^{k-1}\right) \delta + B(\lambda z) B(\lambda)^{-1} \beta \right\).
\end{align*}
Notons que la somme est bien valide car dans chaque coordonn\'ee, la somme est finie. Pour le passage \`a la derni\`ere ligne, on a utilis\'e les formules pour les puissances du lemme~\ref{le:puis}. Afin d'avoir une expression plus explicite, on introduit les polyn\^omes $q_{i,n}(\lambda)$ d\'efinis par :
\begin{equation} \label{eq:q_def}
q_{i,n}(\lambda) = \sum_{\substack{n_1 + \ldots + n_i = n\\n_1 \geq 1, \ldots, n_i \geq 1}} \binom{\lambda}{n_1} \ldots \binom{\lambda}{n_i}.
\end{equation}
On pose \'egalement $q_{0,0} = 1$ et $q_{0,n} = 0$ si $n > 0$ et $q_{i,0} = 0$ si $i > 0$.
De mani\`ere \'equivalente, nous pouvons d\'efinir les $q_{i,n}(\lambda)$ par
\begin{equation} \label{eq:series_formelles}
(A^\lambda - 1)^i = \sum_{n=0}^\infty q_{i,n}(\lambda) X^n.
\end{equation}
On a par exemple:
\[
\begin{array}{| c | | c| c | c| c| c| c| c|}  
\hline  
i \backslash n & 0 & 1 & 2 & 3 & 4 & 5  \\
\hline \hline
0 & 1 & 0 & 0 & 0 & 0 & 0 \\
\hline
1 & 0 & \lambda & \binom{\lambda}{2} & \binom{\lambda}{3} & \binom{\lambda}{4} & \binom{\lambda}{5} \\
\hline
2 & 0 & 0 & \lambda^2 & \lambda^3 - \lambda^2 & \frac{7\lambda^4-18\lambda^3+11\lambda^2}{12} 
  & \frac{3\lambda^5 - 14\lambda^4 + 21\lambda^3 - 10 \lambda^2}{12}\\
\hline
3 & 0 & 0 & 0 & \lambda^3 & \frac{3}{2} (\lambda^4-\lambda^3) & \frac{5\lambda^5 - 12\lambda^4 + 7 \lambda^3}{4} \\
\hline
4 & 0 & 0 & 0 & 0 & \lambda^4 & 2(\lambda^5-\lambda^4) \\
\hline
5 & 0 & 0 & 0 & 0 & 0 & \lambda^5 \\
\hline
\end{array}
\]

\begin{proposition}  \label{prop:puiss_avec_qij}
Pour tout entier $n \geq 0$ on a
\[
\binom{\lambda x}{n} = \sum_{i=0}^{n} q_{i,n}(\lambda) \binom{x}{i}.
\]
Ou de mani\`ere \'equivalente
\[
B(\lambda)^{-1}\ B(\lambda x) = \sum_{i=1}^\infty \binom{x}{i} (A^\lambda - 1)^{i-1}.
\]
\end{proposition}
Avant d'entammer la preuve de cette proposition on en d\'eduit une forme plus explicite pour le morphisme $\cL$.
\begin{corollaire} \label{cor:morph}
Soit $\lambda \in \RR$ et $\beta,\delta \in V$. Soit $\cL$ l'endomorphisme de $G$ tel que 
	$\cL \big(\(1,0\)\big) = \(\lambda,\beta\)$ et $\cL \big( \(0,1\)\big) = \(0,\delta\)$. 
Alors pour un \'el\'ement $\(z,s\)$ de $G$, on a
\[
\cL \big( \(z,s\) \big) = \left\( \lambda z, \sum_{k=1}^\infty (A^\lambda - 1)^{k-1} \ \left(s_k \delta + \binom{z}{k} \beta \right) \right\).
\]
Aini, en \'ecrivant $\cL \(z,s\) = \(\lambda z, t\)$ avec $t = (t_1,t_2,\ldots)$ on a pour tout $k \geq 1$ l'expression
\[
t_k = \sum_{j=1}^k \sum_{i=1}^{k-j+1} q_{j-1,k-i}(\lambda) \left(\delta_i s_j + \beta_i \binom{z}{j}\right).
\]
\end{corollaire}

\begin{proof}[Preuve de la proposition~\ref{prop:puiss_avec_qij}]
Soit $x$ entier. Nous avons d'une part
\[
A^{\lambda x} = \sum_{k=0}^{+\infty} \binom{\lambda x}{k} X^k.
\]
D'autre part, en \'ecrivant $(A^\lambda)^x$ on a
\[
A^{\lambda x}
= \sum_{n=0}^\infty \sum_{\substack{k_1 + \ldots + k_x=k \\ k_1\geq 0,\ldots,k_x \geq 0}} \binom{m}{k_1} \ldots \binom{m}{k_x} X^k
= \sum_{k=0}^\infty \sum_{j=0}^x \binom{x}{j} q_{j,k}(\lambda) X^k.
\]
Comme $\binom{x}{j}$ s'annule pour $x$ entier et $j > x$ on a pour tout $x$ entier
\[
\binom{\lambda x}{k} = \sum_{j=0}^k \binom{x}{j} q_{j,k}(\lambda).
\]
Comme les deux membres de l'\'equation ci-dessus sont polynomiaux en $x$, l'\'egalit\'e est v\'erifi\'ee pour tout $x$ dans $\RR$.
Maintenant, en utilisant l'expression obtenue, nous pouvons r\'e\'ecrire
\begin{align*}
B(\lambda x)
= \sum_{k=0}^\infty \left(\sum _{i=0}^k q_{i,k}(\lambda) \binom{x}{i}\right) X^{k-1} 
= \sum_{i=0}^\infty \left(\sum_{k=0}^\infty q_{i,k}(\lambda) X^{k-1}\right) \binom{x}{i} 
= \sum_{i=0}^\infty \binom{m}{i} \frac{(A^\lambda-1)^i }{X}.
\end{align*}
Et donc
\[
B(\lambda)^{-1} B(\lambda x)
= \frac{X}{A^\lambda - I} \sum_{i=0}^\infty \binom{m}{i} \frac{(A^\lambda-1)^i }{X} 
= \sum_{i=0}^\infty \binom{m}{i} (A^\lambda -1)^{i-1}.
\qedhere
\]
\end{proof}

Nous pouvons d\'eduire imm\'ediatement de la proposition le corollaire suivant :
\begin{corollaire} \label{cor}
Soit $\(z,s\)$ un \'el\'ement de $G$ et soit $\cL: G \rightarrow G$ d\'efini par 
$\cL \big(\(1,0\)\big) = \(\lambda,\beta\)$ et $\cL \big( \(0,1\) \big)= \(0,\delta\)$.

Si $\delta_1=0$, la deuxi\`eme coordon\'ee du  $j$-\`eme coefficient de l'\'el\'ement $\cL^{(n)} \big(\(z,s\)\big)$ ne d\'epend que de 
$(m,s_1,\ldots,s_{j-n})$.
\end{corollaire}

\begin{proof}
Le corollaire d\'ecoule de la forme explicite de l'endomorphisme donn\'ee dans le corollaire~\ref{cor:morph}. La seconde partie s'obtient par r\'ecurrence sur $n$.
\end{proof}

Nous reviendrons dans la partie \ref{subse:qij} sur les propri\'et\'es asymptotiques des polyn{\^o}mes $q_{i,n}$.

\section{Cobords et approximation polynomiale des sommes it\'er\'ees} \label{se:approx}
Le but de cette partie est d'introduire et d'\'etudier le comportement des sommes it\'er\'ees qui nous permettrons de construire \`a la limite une solution \`a l'\'equation~\eqref{eq:E_lambda_delta}. Pour cela nous \'etudions le comportement asymptotique des it\'er\'ees des endomorphismes $\cL$ associ\'es aux paires $(\phi,\sigma)$ o\`u $\phi: \{a,b\} \rightarrow \RR$ est une fonction et $\sigma$ une substitution sur $\{a,b\}$. En effet, la puissance $k$-\`eme de $\cL$ est directement reli\'ee \`a la valeur des sommes de Birkhoff de $\phi$ au temps $n = \lambda^k$

\subsection{Approximation polynomiale.} \label{subsec:approximation_polynomiale}
Soit $\lambda > 0$ et $\beta,\delta \in V$. Comme pr\'ec\'edemment, on consid\`ere le morphisme $\cL$ tel que 
$\cL \big(\(1,0\)\big) = \(\lambda,\beta\)$ et $\cL \big((0,1\)\big) = \(0,\delta\)$. Afin d'\'etudier les puissances de l'endomorphisme $\cL$ on introduit des polyn\^omes $R_{i,n}(x)$ en posant pour $i \geq 0$ :
\begin{equation} \label{eq:R_def}
\cL^i \big(\(x,(0,\ldots,0)\) \big) = \left\( \lambda^i x,\, (R_{i,1}(x), R_{i,2}(x), R_{i,3}(x), \ldots) \right\).
\end{equation}
Les $R_{i,n}(x)$ sont bien des polyn\^omes en $x$ (d\'ependants des param\`etres $\beta$, $\delta$ et $\lambda$) d'apr\`es la formule pour les puissances et le fait que $\cL^i \big(\(x,0\)\big) = \Big(\cL^i \big( \(1,0\)\big)\Big)^x$.
De plus, en \'ecrivant $\cL^{i+1}\big(\(x,0\)\big) = \cL \big( \cL^i (\(x,0\))\big)$ on tire la relation de r\'ecurrence
\[
(R_{i+1,1}(x), R_{i+1,2}(x), \ldots) = \sum_{k=1}^\infty (A^\lambda - 1)^{k-1} \, \left(R_{i,k}(x) \delta + \binom{m}{k} \beta\right).
\]

Commen\c cons  par exprimer les polyn\^omes $R_{i,n}$ pour $i \geq n$ \`a l'aide des valeurs diagonales $R_{i+1,i+1}(1/\lambda^{i+1})$. Ces derni\`eres peuvent \^etre vues comme un vecteur propre de $\cL$ sur $G$.
\begin{proposition} \label{prop:R_et_c}
Soit $\lambda > 0 $ et $\beta \in V$ et $\delta \in V_2$ (i.e. $\delta = (0, \delta_2, \delta_3,\ldots)$). Soit $\cL$ l'endomorphisme de $G$ tel que
$\cL \(1,0\) = \(\lambda,\beta\)$ et $\cL \(0,1\) = \(0,\delta)$. On pose
\[
c_i = R_{i+1,i+1} \left(\frac{1}{\lambda^{i+1}} \right).
\]
Alors pour tout $i \geq 1$ et $n \geq 0$ on a
\[
R_{i+n,i}(x) = R_{i,i}(\lambda^n x)
\qquad \text{et} \qquad
R_{i,i}(x / \lambda^i) = \sum_{k=0}^{i-1} \binom{x}{i-k} c_k.
\]
De plus
\[
\cL \big(\(1, (c_0, c_1, c_2, \ldots) \) \big) = \(1, (c_0, c_1, c_2, \ldots)\)^{\lambda}.
\]
\end{proposition}

Pour la substitution de Prouet-Thue-Morse, les constantes $c_i$ sont toutes \'egales \`a z\'ero. On verra dans la proposition suivante, que ceci implique que les sommes de Birkhoff it\'er\'ees sont toutes born\'ees. Une autre cons\'equence est que la fonction limite $f_{\w}$ est nulle en z\'ero. Par contre pour les deux substitutions $a \mapsto aab, b \mapsto baa$ et $a \mapsto aab, b \mapsto aba$ ce n'est pas le cas.

\[
\begin{array}{|l|c|c|c|c|c|c|}
\hline
& c_0 & c_1 & c_2 & c_3 & c_4 & c_5 \\
\hline
a \mapsto aab, b \mapsto aba & 2/3 & 1/3 & 10/9& 11 & 8567/27 & 718435/27 \\
\hline
a \mapsto aab, b \mapsto baa &  1/3 & 1/3 & 23/9 & 440/9 & 74431/27 & 455949 \\
\hline
\end{array}
\]

\begin{proof}
On fixe un entier $i$ et nous faisons les calculs dans $G_i = G / V_{i+1}$.

Comme $\delta_1=0$ pour tout $s = (s_1, s_2, \ldots) \in V$ on a d'apr\`es le corollaire~\ref{cor}
\[
\cL^i \big( \(z, s\)\big) = \( \lambda^i z, (R_{i,1}(z), R_{i,2}(z), \ldots, R_{i,i}(z)) \).
\]
En particulier, $\cL^i \big(\(\lambda^{n} x, 0 \)\big)  = \cL^i \circ \cL^n \big( \(x,0\)\big)$ et donc
\[
\left\( \lambda ^{i+n} x, \left(R_{i,1}(\lambda^{n} x ),\ldots ,R_{i,i}(\lambda ^{n} x) \right) \right\)
=
\left\( \lambda ^{i+n}  x, \left( R_{i+n,1}(x),\ldots ,R_{i+n,i}(x) \right)\right\).
\]
Ce qui d\'emontre la premi\`ere \'equation.

Comme $\cL$ est un morphisme et $\(n,0\) = \(1/\lambda^i,0\)^{n\lambda^i}$ on a :
\begin{align*}
\cL^i \big( \(n,0n\)  \big)
&= \(n \lambda^{i},(R_{i,1}(n),\ldots,R_{i,i}(n)) \) \\
&= \left( \cL^i \big( \(1/\lambda^i,0\)\big) \right)^{n \lambda^i}  \\
&= \left\( 1, \left( R_{i,1}(1/\lambda^i),\ldots,R_{i,i}(1/\lambda^i) \right) \right\)^{n \lambda^i}.
\end{align*}
En utilisant la formule pour les puissances, on obtient que pour tout $1 \leq j \leq i$
\[
R_{i,j}(n) = \sum_{k=1}^j \binom{n \lambda^i}{j - k + 1} R_{i,k}(1 / \lambda^i).
\]
En utilisant la premi\`ere \'equation on peut remplacer $R_{i,k}(1 / \lambda^i)$ par $R_{k,k}(1/\lambda^k)$. On trouve alors la seconde relation en prenant $j=i$.

Pour d\'emontrer la derni\`ere partie, remarquons tout d'abord que
\begin{align*}
\cL^i \left( \left\( \frac{1}{\lambda^i}, 0 \right\)\right)
&= \(1, (R_{i,1}(1/\lambda^i), R_{i,2}(1/\lambda^i), \ldots, R_{i,i}(1/\lambda^i)) \) \\
&= \(1, (R_{1,1}(1/\lambda), R_{2,2}(1/\lambda^2), \ldots, R_{i,i}(1/\lambda^i)) \) \\
&= \(1, (c_0, c_1, \ldots, c_{i-1}) \).
\end{align*}
On a alors
\[
\cL \big( \(1, (c_0, c_1, \ldots, c_{i-1})\) \big)
= \cL^{i+1} \big( \( 1/\lambda^i, 0\) \big)
= \left(\cL^{i+1} \big( \( 1/\lambda^{i+1}, 0\) \big) \right)^\lambda 
= \(1, (c_0, \ldots, c_{i-1})\)^\lambda.
\qedhere
\]
\end{proof}

On se sert maintenant des polyn\^omes $R_{i,i}(x)$ pour redresser les sommes de Birkhoff.
Soit $\sigma$ une substitution de longeur constante $\lambda$ telle que $\delta_1(\sigma) = 0$.
Rappellons que cette condition signifie que $\sigma(a)$ et $\sigma(b)$ contiennent le m\^eme nombre de $a$.
Soit $\phi: \{a,b\} \rightarrow \RR$ une fonction non constante.
On associe \`a la fonction $\phi$ et la substitution $\sigma$ l'endomorphisme $\cL = \cL_{\phi,\sigma}$ de $G$ d\'efini dans la partie~\ref{subse:induit}.

Soit $R_{i,n}$ et $c_i$ les polyn\^omes et les constantes d\'efinis dans la proposition~\ref{prop:R_et_c}.
On pose pour tout entier $\ell \geq 1$ et tout r\'eel $x$ :
\begin{equation} \label{eq:def_p}
p_\ell (x)  = \sum_{i=0}^\ell c_i \binom{x}{\ell-i} = R_{\ell,\ell}(x/\lambda^\ell) + c_\ell.
\end{equation}

\begin{proposition} \label{prop:approximation_polynomiale}
Soit $\sigma$ une substitution de longueur constante $\lambda$ telle que $\delta_1(\sigma) = 0$.
Soient $p_\ell$ les polyn\^omes d\'efinis en~\eqref{eq:def_p}.
Alors, pour tout entier $\ell\geq 1$ et tout mot $\u \in \sigma^\ell(\{a,b\}^\NN)$ on a:
\begin{itemize}
\item pour tout r\'eel $x$, $p_{\ell}(x+1) - p_{\ell}(x) = p_{\ell-1} (x)$ ;
\item pour tout entier $n$,
\begin{equation} \label{eq:sommeetRbisbis}
\big{(} S^{(\ell)}_{n+1}(\phi,\u) - p_{\ell}(n+1) \big{)}  - \big{(} S^{(\ell)}_{n}(\phi,\u) - p_{\ell} (n) \big{)} = S^{(\ell-1)}_{(n)}(\phi,\u) - p_{\ell-1}(n).
\end{equation}
\item pour tout entier $n$,
\begin{equation} \label{eq:sommeetRbis}
S^{(\ell)}_{n \cdot \lambda^{\ell}}(\phi,\u)   = p_\ell (n \lambda^\ell) - c_\ell.
\end{equation}
\item la suite $\big{(}  S^{(\ell)}_n(\phi,\u) - p_\ell(n) \big{)}_{n\geq 0}$ est born\'ee.
\end{itemize}
\end{proposition}

\begin{proof}
Le premier point est imm\'ediat par la d\'efinition des polyn{\^o}mes $(p_\ell)_\ell$. En effet pour tout nombre r\'eel $x$ et tout entier $\ell \geq 1$ :
\[
p_\ell(x+1) - p_\ell (x) = \sum \limits_{i=0} ^{\ell} c_i \binom{x+1}{\ell-i} - \sum \limits_{i=0} ^{\ell} c_i \binom{x}{\ell-i}
= \sum \limits_{i=0} ^{\ell-1} c_i \binom{x}{\ell-1-i} = p_{\ell-1}(x).
\]
Le deuxi\`eme point d\'ecoule imm\'ediatement.

\bigskip

On note $\pi = \pi_\phi$ et soit $\u_0 \u_1 \ldots = \sigma^\ell(\v_0 \v_1 \ldots)$ un mot de $\sigma^\ell(\{a,b\}^\NN)$.
On consid\`ere l'endomorphisme $\cL = \cL_{\phi,\sigma}$ du lemme~\ref{lem:L_f_sigma}.
Fixons un entier $n \geq 0$. D'une part, $\cL \circ \pi = \pi \circ \sigma$, et donc
\[
\pi (\u_0\cdots \u_{n \lambda^\ell -1}) =\pi \circ \sigma^\ell (\v_0\cdots \v_{n -1}) =  \cL^\ell \circ \pi (\v_0\cdots \v_{n-1} ) .
\]
On a ainsi dans $G_\ell$
\[
\pi (\u_0\cdots \u_{\lambda^\ell n -1}) = \( n \lambda^\ell, (R_{\ell,1}(n),\ldots,R_{\ell,\ell}(n))\) 
\]
Par le lemme~\ref{lem:connexion}, nous obtenons
\[
S^{(\ell)}_{n \lambda^i}(\phi,\u) = R_{\ell,\ell}(n) = p_\ell(n \lambda^i) - c_\ell.
\]
Ce qui prouve~\eqref{eq:sommeetRbis}.

Le dernier point est facile car les incr\'ements de la suite $(S_k^{(\ell)}(\phi,\u) - p_\ell(k))_{k \geq 0}$ sont en nombre fini, donc born\'es.
\end{proof}

\begin{remarque}
Pour $\ell=2$, nous pouvons interpr\'eter le fait que $S^{(2)}_{\lambda^2}(\phi, \sigma^2(a)) = S^{(2)}_{\lambda^2}(\phi, \sigma^2(b))$ de mani\`ere visuelle. Prenons pour simplifier la fonction $\phi = \chi_a$.
Nous associons \`a tout mot $\v=\v_0\cdots \v_{n-1}$ une ligne bris\'ee $(X_k)_{k}$ dans $\RR^2$ d\'efinie par r\'ecurrence de la mani\`ere suivante. On pose $X_0=(0,0)$ et supposons avoir construit la ligne bris\'ee $(X_0,\ldots,X_{m_k})$ associ\'ee au mot $\v_0\cdots \v_{k-1}$. Si $\v_k = a$, on ajoute les deux termes $X_{m_k} + (0,1)$ et $X_{m_k} + (1,1)$ \`a notre suite. Si $\v_k = b$ on ajoute le terme $X_{m_k} + (1,0)$.
\begin{figure}[H]  
\begin{center}
\includegraphics[scale=0.7]{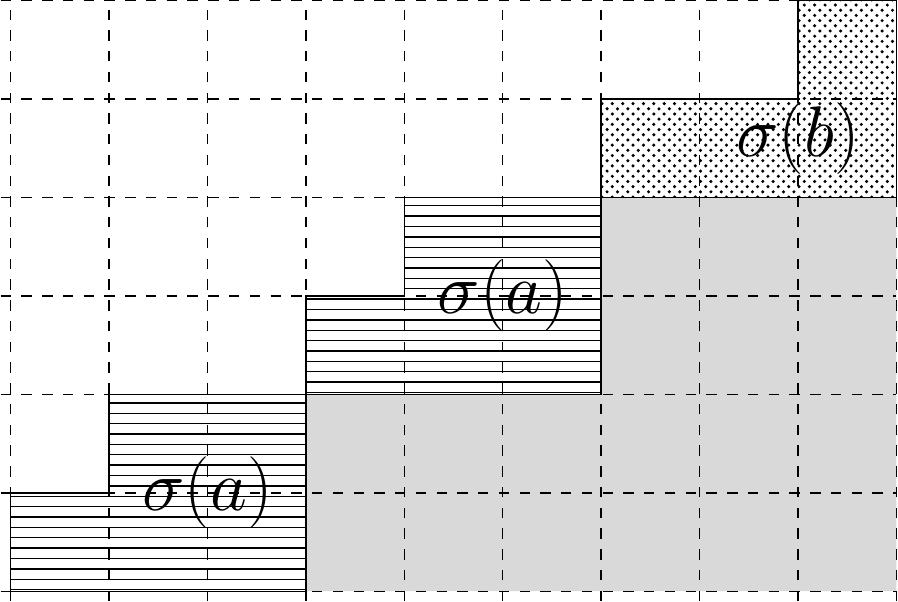} \hspace{1cm} \includegraphics[scale=0.7]{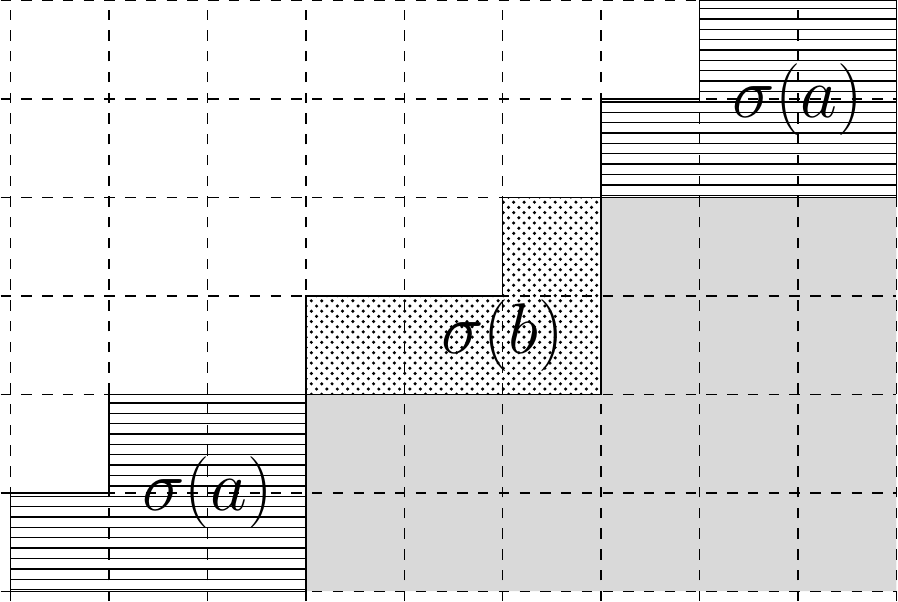}
\end{center}
\caption{Repr\'esentation de la ligne bris\'ee pour $\sigma: a \mapsto aab, b \mapsto aba$.}
\label{fig:S2_et_aire} 
\end{figure}
Comme on peut le voir sur la figure~\ref{fig:S2_et_aire}, la hauteur de la ligne bris\'ee au point $k$ est $S^{(1)}_k(\phi,\v)$ et $S^2_{k} (\chi_a,\v)$ est alors l'aire de la zone sous la ligne bris\'ee. Les aires situ\'ees entre sous les lignes bris\'ees jusqu'\`a l'abscisse $k=9=\lambda^2$ sont \'egales.
\end{remarque}

\subsection{Cobords.} \label{sec:cobords}
Dans cette section on revient sur les coefficients $c_i$ de la proposition~\ref{prop:R_et_c} et leur application au comportement moyen des sommes de Birkhoff dans la proposition~\ref{prop:approximation_polynomiale}. Nous montrons qu'un th\'eor\`eme plus g\'en\'eral est vrai (dans le cadre d'un alphabet quelconque). On montre \'egalement comment ce ph\'enom\`ene de sommes de Birkhoff born\'ees s'articule avec la notion de cobord.

On fixe un d\'ecalage $K \subset {\mathcal A}^\NN$ ou $K \subset {\mathcal{A}}^\ZZ$.

Une fonction continue $\phi: X \rightarrow \RR$ est un \emph{cobord continu} s'il existe une fonction continue $\psi: X \rightarrow \RR$ telle que $\phi = \psi \circ T - \psi$. Autrement dit, la fonction $\phi$ est dans l'image de l'op\'erateur $U_T - I: \psi \mapsto \psi \circ T - \psi$. Si $\phi$ est un cobord, alors la fonction $\psi$ telle que $\phi = \psi \circ T - \psi$ est d\'etermin\'ee \`a une constante pr\`es.
Si $\phi$ est un cobord, alors sa somme de Birkhoff se r\'eecrit simplement en fonction de $\psi$: $S^{(1)}_n(\phi,\u) = \psi (T^n \u) - \psi(\u)$. En particulier, la suite $S^{(1)}_n(\phi,\u)$ est born\'ee. La r\'eciproque est vraie, c'est le th\'eor\`eme de Morse-Hedlund : si $K$ est minimal, une fonction $\phi: K \rightarrow \RR$ dont la somme de Birkhoff est born\'ee est un cobord continu.

Nous dirons que la fonction $\phi$ est {\it un cobord infini} s'il existe une constante $\psi_0$ et une suite de fonctions continues $(\psi_n)_{n\geq 1}$ telle que $\phi + \psi_0 = \psi_1 \circ T - \psi_1$ et pour tout entier $n\geq 1$, 
\[
\psi_{n} = \psi_{n+1} \circ T- \psi_{n+1}.
\]
Autrement dit, la fonction $\phi$ est \`a une constant pr\`es dans $\bigcap_{n \geq 0} (U_T - I)^{n} \mathcal C(K)$.
La suite de fonctions $(\psi_n)_n$ est unique et pour tout entier $n \geq 1$, la moyenne de $\psi_n$ selon toute mesure $T$-invariante $\mu$ est nulle.
Avec ces notations, si $\phi$ est un cobord infini, alors pour tout $\u \in X$ on a
\[
\begin{array}{l}
S^{(1)} _n(\phi ,\u) = \psi_1 \circ T^n(\u) - \psi_1(\u) - n \psi_0 \ , \\
S^{(2)} _n(\phi ,\u)= S^{(1)} _n(\psi_1 ,\u) - n \psi_1(\u) = \psi_2(T^n \u) - \psi_2(\u) - n \psi_1(\u) - \binom{n}{2} \psi_0 \, \\
S^{(3)} _n(\phi ,\u)= S^{(1)} _n(\psi_2 ,\u) - n \psi_2(\u)-\binom{n}{2} \psi_1(\u)  = \psi_3 (T^n \u) - \psi_3(\u)- n \psi_2(\u) - \binom{n}{2} \psi_1(\u) - \binom{n}{3} \psi_0.
\end{array}
\]
On d\'eduit ais\'ement par r\'ecurrence que pour tout $\ell \geq 0$,
\begin{equation} \label{eq:cobord}
S^{(k)} _n (\phi ,\u)=  \psi_k(T^n \u) - \sum_{j=0}^k \binom{n}{k-j} \psi_j(\u). 
\end{equation}
Et plus g\'en\'eralement pour $i \geq 1$ et $k \geq 0$ on a
\begin{equation} \label{eq:cobord_bis}
S^{(k)}_n(\psi_i, \u) = \psi_{i+k}(T^n \u) - \sum_{j=1}^{k} \binom{n}{k-j} \psi_{i+j}(\u).
\end{equation}
En particulier, si on note $p_{\ell,\u}(n) = \sum_{k=0}^\ell \binom{n}{\ell-k} \psi_k(\u)$, alors $S^{(\ell)}_n(\phi, \u) - p_{\ell,\u}(n) = \psi_\ell(T^n \u)$ est une suite born\'ee. Autrement dit, les sommes de Birkhoff it\'er\'ees sont \`a distance born\'ee des valeurs d'un polyn\^ome. Et il est facile de voir que $p_{\ell, T\u}(n) = p_{\ell,\u}(n+1)$.

\bigskip

Dans la section pr\'ec\'edente, $\mathcal A=\{a,b\}$ et nous travaillions avec une fonction $\phi$ ne d\'ependant que de la premi\`ere lettre d'un mot, sur un d\'ecalage $K_\sigma$ provenant d'une substitution. Les constantes $c_i$ obtenues dans la proposition~\ref{prop:approximation_polynomiale} v\'erifient alors $c_i = - \psi_i(\w)$. Notons que dans ce cas, les $c_i$ d\'eterminent enti\`erement les fonctions $\psi_i$


\bigskip

Dans cette section, nous montrons que la proposition~\ref{prop:approximation_polynomiale} se g\'en\'eralise de la mani\`ere suivante.
\begin{proposition} \label{thm:cobord_infini}
Soit $\mathcal A$ un alphabet fini et $v=(v_\alpha)_{\alpha \in \cA}$ un vecteur d'entiers strictement positifs.
Soit $\sigma$ une substitution sur $\mathcal A$ telle que pour chaque couple de lettres $\alpha,\beta \in \mathcal A$, 
$\sigma(\alpha)$ contienne $v_\beta$ fois la lettre $\beta$. Autrement dit, l'image de chacune des lettres a le m\^eme ab\'elianis\'e $v$.
On note $\lambda = \sum_{\alpha \in \cA} v_\alpha$ la longueur de la substitution.
Alors toute fonction $\phi: \cA \rightarrow \RR$ de moyenne nulle est un cobord infini sur le d\'ecalage $(K_\sigma,T)$.
\end{proposition}
Pour une lettre $\alpha$, on note $[\alpha] = \{\u \in K_{\sigma,\ZZ} \ ; \  \u_0 = \alpha\}$.
Pour $0 \leq m < \lambda^k$, notons $\Cyl(k,m,\alpha) = T^m \sigma^k([\alpha])$.
Pour chaque entier $k \geq 0$, $K_\sigma$ se d\'ecompose de la mani\`ere suivante :
\begin{equation} \label{eq:decomposition_K_sigma}
K_\sigma = \bigcup_{\alpha \in \cA} \bigcup_{0 \leq m < \lambda^k} \Cyl(k,m,\alpha).
\end{equation}

Si $\sigma$ est une substitution positive, deux cas peuvent se produire: soit $K_{\sigma,\ZZ}$ est un ensemble fini (auquel cas tous les mots de $K_{\sigma,\ZZ}$ sont p\'eriodiques) soit $K_{\sigma,\ZZ}$ est infini et alors la d\'ecomposition~\eqref{eq:decomposition_K_sigma} est une r\'eunion disjointe (voir~\cite{Mosse}). Le premier cas appara\^it par exemple lorsque toutes les images de $\sigma$ sont identiques. Le th\'eor\`eme \'etant \'el\'ementaire pour les mots p\'eriodiques nous supposons que $K_{\sigma,\ZZ}$ est infini.

On note $F^k$ l'ensemble des fonctions de $K_\sigma$ dans $\RR$ qui sont constantes sur chaque $\Cyl(k,m,\alpha)$ et $F^k_0$ les fonctions de $F^k$ de moyenne nulle.

Le th\'eor\`eme ci-dessus se d\'eduit du lemme suivant.

\begin{lemme} \label{le:lemmecob}
Soit $\sigma$ une substitution sur $\cA$ v\'erifiant les m\^emes hypoth\`eses que celles du  th\'eor\`eme~\ref{thm:cobord_infini}.
Soit $\phi \in F^k_0$.
Alors $\phi$ est un cobord et si $\psi$ est telle que $\phi = \psi \circ T - \psi$ alors $\psi$ appartient \`a $F^{k+1}$.

Autrement dit, \[F^k_0 \subset (U_T - I)(F^{k+1}_0).\]
\end{lemme}

Ce lemme se d\'eduit lui-m\^eme du r\'esultat suivant 
\begin{lemme} \label{le:mesure}
Avec les m\^emes hypoth\`eses que dans le th\'eor\`eme~\ref{thm:cobord_infini}.
Pour tout $k,m,\alpha$, 
\[
\mu \left( \Cyl(k,m,\alpha) \right) = \frac{v_\alpha}{\lambda^{k+1}}.
\]
D'autre part, si $\phi \in F^k$ alors pour tout $\alpha$ et tout $\u \in \Cyl(k+1,0,\alpha)$ on a
\[
\sum_{m = 0}^{\lambda^{k+1}-1} \phi(T^m \u) = \int_{K_\sigma} \phi d\mu,
\]
o\`u on rappelle que $\mu$ est l'unique mesure de probabilit\'e sur $K_\sigma$ invariante par $T$.
\end{lemme}

\begin{proof}
Comme la mesure $\mu$ est invariante, $\mu(\Cyl(k,m,\alpha)) = \mu(\Cyl(0,m,\alpha))$.
On note $\mu_k$ la mesure d\'efinie pour tout ensemble mesurable $Y$ par
\[
\mu_k(Y) = \mu(T^k \sigma (Y)).
\]
Alors $\mu_0$, $\mu_1$, \ldots, $\mu_{\lambda-1}$ sont \`a supports disjoints et
$\mu_0 + \mu_1 + \ldots + \mu_{\lambda-1}$ est une mesure invariante dont il est facile de voir
que c'est une mesure de probabilit\'e. Par unique ergodicit\'e, elle est \'egale \`a $\mu$ et on
en d\'eduit que pour tout ensemble mesurable $Y \subset K_\sigma$
\[
\mu(\sigma(Y)) = \frac{\mu(Y)}{\lambda}.
\]
Ceci prouve la premi\`ere partie du lemme.

La seconde partie d\'ecoule du fait que chaque image de $\sigma$ a le m\^eme ab\'elianis\'e.
\end{proof}

\begin{proof}[Preuve du lemme \ref{le:lemmecob}]
Soit $\phi \in F^k_0$, prenons $\u \in \Cyl(k+1,0,\alpha)$ et posons pour $m < |\sigma|^{k+1}$, $\psi(T^m \u) = S^{(1)}_m(\phi,\u)$.
La fonction $\psi$ s'\'etend de mani\`ere unique en une fonction de $F^{k+1}$.
D'apr\`es la seconde partie du lemme~\ref{le:mesure}, on a pour tout $\u \in \Cyl(k+1,0,\alpha)$, et pour tout $m \geq 0$, $\psi(T^m \u) = S^{(1)}_m(f,\u)$. Et donc pour tout $\u \in K_\sigma$, $S^{(1)}_m(\phi,\u) = \psi(T^m \u) - \psi(\u)$.
\end{proof}

\begin{remarque}
Le th\'eor\`eme~\ref{thm:cobord_infini} ne se g\'en\'eralise pas simplement \`a d'autres substitutions. Les fonctions qui ne d\'ependent que de la premi\`ere lettre et qui sont des cobords se lisent sur la matrice d'incidence de la substitution (voir~\cite{adam1}).
Par exemple, pour les mots sturmiens sur $\{a,b\}$ (qui sont des codages de rotations et dont certains sont substitutifs), il est bien connu que les fonctions $\chi_a - \mu([a])$ et $\chi_b - \mu([b])$ sont des cobords: il existe $\alpha$ tel que
\[
\left\lvert S^{(1)}_n(\chi_a, \u) - k \alpha \right\rvert
\qquad \text{est born\'ee.}
\]
Cependant, il est montr\'e dans~\cite{pinner} que pour tout param\`etre $\beta$,
\[
\left\lvert S^{(2)}_n(\chi_a, \u) - \alpha \binom{n}{2} - \beta n \right\rvert
\qquad \text{n'est par born\'ee.}
\]
Plus g\'en\'eralement, l'\'etude des sommes de Birkhoff it\'er\'ees pour un syst\`eme substitutif est un sujet int\'eressant et assez peu d\'evelopp\'e mais sortant du cadre de cet article.
\end{remarque}

\subsection{Une remarque sur le mot de Prouet-Thue-Morse.}
La section pr\'ec\'edente pourrait laisser penser que toutes les fonctions continues sont des cobords.
Dans cette section, nous montrons que la plupart des fonctions sur le d\'ecalage $K_\sigma$ associ\'e au mot de Prouet-Thue-Morse et qui d\'ependent des deux premi\`eres lettres ne sont pas des cobords. Ce r\'esultat est un corollaire des travaux de B. Adamczewski~\cite{adam1}. 

Rappelons, que le mot de Prouet-Thue-Morse est un point fixe de la substitution $\sigma: a \mapsto ab, b \mapsto ba$.
Pour \'etudier les sommes de Birkhoff des fonctions sur deux lettres, il suffit de consid\'erer la substitution associ\'ee aux mots de deux lettres.  En notant$A = aa$, $B=ab$, $C=ba$ et $D=bb$ il s'agit de
\[
\sigma_2: A \mapsto BC \quad B \mapsto BD \quad C \mapsto CA \quad D \mapsto CB.
\]
On obtient le mot de Prouet-Thue-Morse en prenant le point fixe de $\sigma_2$ commen\c{c}ant par $B$ et en projettant sur l'alphabet $\{a,b\}$ via le morphisme $\tau: A,B \mapsto a$ et $C,D \mapsto b$.

Les sommes de Birkhoff de $\chi_{A}$ sur le d\'ecalage $K_{\sigma_2,\ZZ}$ sont \'egales \`a celles de la fonction $\chi_{aa}$ sur le d\'ecalage $K_{\sigma,\ZZ}$ o\`u
\[
\chi_{aa}(\u) = \left\{ \begin{array}{ll}
1 & \text{si $\u_0 \u_1 = aa$} \\
0 & \text{sinon.}
\end{array} \right.
\]
Plus pr\'ecis\'ement, pour tout mot $\u \in K_{\sigma_2,\ZZ}$ nous avons $S_n(\chi_A,\u) = S_n(\chi_{aa}, \tau(\u))$.

Les valeurs propres de la matrice d'incidence de $\sigma_2$ sont $2$, $0$, $1$ et $-1$. Chacune d'elles est associ\'ee \`a une fonction qui d\'epend d'au plus deux lettres

\[
\phi_2 = 1 \quad \phi_0 = \chi_a - \chi_b \quad \phi_1 = \chi_{ab} - \chi_{ba} \quad \phi_{-1} = 2(\chi_{aa}+\chi_{bb}) - (\chi_{ab}+\chi_{ba}).
\]
Ce sont des vecteurs propres \`a gauche pour la matrice d'incidence : pour $\alpha =2,0,1,-1$ et tout $n$ entier on a
\[
S_{2 n}(f_\alpha, \sigma(\u)) = \alpha S_n (f_\alpha, \u).
\]
En appliquant le crit\`ere de~\cite{adam1} on obtient.
\begin{proposition}
Soit $F = \RR \chi_{aa} \oplus \RR \chi_{ab} \oplus \RR \chi_{ba} \oplus \RR \chi_{ab} = \RR \phi_2 \oplus \RR \phi_0 \oplus \RR \phi_1 \oplus \RR \phi_{-1}$ les fonctions qui d\'ependent d'au plus deux lettres sur le d\'ecalage de Prouet-Thue-Morse.
Alors l'ensemble des cobords est $\RR \phi_0 \oplus \RR \phi_1$.
En particulier, pour tout mot de deux lettres $\v$ et toute constante $c$ la fonction $\chi_\v - c$ n'est pas un cobord.
\end{proposition}

Remarquons que la raison pour laquelle le th\'eor\`eme~\ref{thm:cobord_infini} ne s'applique pas est que les cylindre $[\v]$ ne se d\'ecompose pas simplement sur les cylindres $\Cyl(k,m,\alpha)$. Par exemple, pour $[aa]$, on peut \'ecrire
\[
[aa] = T \sigma^2(b) \cup T^3 \sigma^2(bb)
\]
dont on d\'eduit que
\[
[aa] = \Cyl(2,1,b) \cup \Cyl(4,9,a) \cup \Cyl(8,21,b) \cup \ldots
\]

\section{Asymptotiques} \label{se:asymptotique}
Dans cette section nous \'etudions en d\'etail les polyn\^omes $q_{i,n}(\lambda)$ intervenant dans les morphismes $\cL$, en particulier certaines valeurs asymptotiques. On en d\'eduit alors le comportement des coefficients $c_\ell$ et des polyn\^omes $R_{i,n}$ introduits dans la section~\ref{subsec:approximation_polynomiale}.

\subsection{Polyn{\^o}mes $q_{i,n}(\lambda)$.} \label{subse:qij}
Le principal outil de cette section va \^etre la r\'ecurrence suivante.
\begin{lemme} \label{lem:q_recurrence}
Les polyn\^omes $q_{i,n}$ v\'erifient
\[
q_{1,n} = \binom{\lambda}{n} \, , \qquad q_{n,n} = \lambda^{n} \, , \qquad q_{n,n+1} = \frac{\lambda-1}{2} n \lambda^n.
\]
et pour $2 \leq i \leq n$
\[
q_{i,n+1}  =  \frac{\lambda i}{n+1} \,  q_{i-1,n} + \frac{\lambda i - n }{n+1} \, q_{i,n}.
\]
\end{lemme}

\begin{proof}
En d\'erivant en $X$ l'\'egalit\'e de s\'eries formelles
\[
(A^\lambda-1)^i = \sum_{n \geq 0} q_{i,n} X^n
\]
on obtient
\begin{equation} \label{eq:relation_q_derive_serie}
n q_{i,n} = i \sum_{m=1}^{n-1} m \binom{\lambda}{m} q_{i-1,n-m}.
\end{equation}
Ainsi, en \'ecrivant $\binom{\lambda}{m+1} = \frac{\lambda - m + 1}{m} \binom{\lambda}{m}$, on a
\begin{align*}
n q_{i,n} &= i \sum_{m = 1}^{n-1} m \binom{\lambda}{m} q_{i-1,n-m} 
	= \lambda i q_{i-1,n-1} + i \sum_{m=2}^{n-1} m \binom{\lambda}{m} q_{i-1,n-m} \\
	&= \lambda i q_{i-1,n-1} + i \sum_{m=1}^{n-2} (m+1) \binom{\lambda}{m+1} q_{i-1,n-m-1} \\
	&= \lambda i q_{i-1,n-1} + i \sum_{m=1}^{n-2} (\lambda - m) \binom{\lambda}{m} q_{i-1,n-m-1} \\
	&= \lambda i q_{i-1,n-1} + i \lambda \sum_{m=1}^{n-2}  \binom{\lambda}{m} q_{i-1,n-m-1} 
		-i  \sum_{m=1}^{n-2} m \binom{\lambda}{m} q_{i-1,n-m-1} \\
	&= \lambda i q_{i-1,n-1} + \lambda i q_{i,n-1} - i\sum_{\alpha_1 + \ldots + \alpha_i = n-1} \alpha_1 \binom{\lambda}{\alpha_1} \ldots \binom{\lambda}{\alpha_i}.
\end{align*}
On conclut en utilisant \`a nouveau~\eqref{eq:relation_q_derive_serie}.
\end{proof}

\begin{proposition} \label{pr:q}
Pour tout $i \geq 1$, on a lorsque $n \to \infty$
\begin{equation} \label{eq:q}
q_{n,n+i}(\lambda) \underset{n}{\sim} \frac{(\lambda-1)^i}{2^i i!} \  n^i \lambda^n .
\end{equation}
\end{proposition}

\begin{remarque}
Pour $i \geq 1$, on peut montrer par r\'ecurrence qu'il existe des polyn\^omes $P_i \in \QQ[\lambda,n]$ tel que
\[
q_{n,n+i}(\lambda) = \frac{\lambda^n (\lambda-1)}{i! 2^i} P_i(\lambda,n).
\]
De plus
\[
P_i(\lambda, n/(\lambda+1)) = n^i + a^{(i)}_{1}(\lambda) \ n^{i-1} + \ldots + a^{(i)}_{i-2}(\lambda)  n^2 + a^{(i)}_{i-1}(\lambda) n.
\]
o\`u $a^{(i)}_k \in \QQ[\lambda]$ est de degr\'e $k$.

Par exemple :
\begin{align*}
P_1\left(\lambda,\frac{n}{\lambda+1}\right) &= n \\
P_2\left(\lambda,\frac{n}{\lambda+1}\right) &= n^2 + 1/3(\lambda - 3) n \\
P_3\left(\lambda,\frac{n}{\lambda+1}\right) &= n^3 + (\lambda -3)\ n^2 + 2(-\lambda+1)n \\
P_4\left(\lambda,\frac{n}{\lambda+1}\right) &= n^4 + 2(\lambda-3)\ n^3 + \left(\frac13 \lambda^2 - 10 \lambda + 11\right)\ n^2 
	+\left(- \frac 2 {15} \lambda^3 - \frac{14}{15} \lambda^2 + \frac{62}{5} \lambda - 6 \right)n.
\end{align*}
\end{remarque}

\begin{proof}
Posons
\[
Q_{i,n} = \frac{n (n+1) \ldots (n+i-1)}{\lambda^n} \ q_{n,n+i}.
\]
Alors, $Q_{0,n} = q_{n,n} = \lambda^n$ et la r\'ecurrence du lemme~\ref{lem:q_recurrence} nous donne
\[
Q_{i+1,n} = Q_{i+1,n-1} + (n(\lambda-1) - i) Q_{i,n}.
\]
On en extrait la formule
\begin{equation} \label{eq:Q_recurrence}
Q_{i+1,n} = \sum_{m=0}^n (m(\lambda-1) - i) Q_{i,m}.
\end{equation}

On d\'emontre maintenant par r\'ecurrence que $Q_{i,n} \sim (\lambda-1)^i n^{2i}/ (2^i \ i!)$.
Cette formule est vraie pour $i=0$, supposons la vraie jusqu'au rang $n$. Alors, en utilisant la formule~\eqref{eq:Q_recurrence} on trouve
\begin{align*}
Q_{i+1,n} \sim (\lambda-1) \sum_{m=0}^n m Q_{i,m} 
\sim \frac{(\lambda-1)^{i+1}}{2^i\ i!} \sum_{m=0}^n m^{2i+1} 
\sim \frac{(\lambda-1)^{i+1}}{2^i\ i!}\ \frac{n^{2i+2}}{2i+2}
\end{align*}
En revenant \`a $q_{n,n+i}$, on obtient le r\'esultat.
\end{proof}

\begin{lemme} \label{le:techniquepourq}
Soit $\lambda$ un entier sup\'erieur ou \'egal \`a $1$. Alors
\[
\sup \{ q_{i,j}({\lambda}) \ ; \  1\leq i \leq j \leq n \}  \leq (2\lambda-1)^n.
\]
\end{lemme}

\begin{proof}
Notons $A_n = \sup \{ q_{i,j}({\lambda}) ; 1\leq i \leq j \leq n \}$.
Nous montrons ce r\'esultat par r\'ecurrence.
On a $A_1 = |q_{1,1}({\lambda})|={\lambda}$. Fixons un entier $n$ et $i\in\{2,\ldots,n\}$, alors :
\begin{flalign*}
q_{i,n+1}({\lambda})  
	&=  \frac{\lambda i }{n+1}  q_{i-1,n}({\lambda}) +\frac{ \lambda  i -n }{n+1}  q_{i,n}({\lambda}) 
		\leq   \frac{\lambda n }{n+1}  A_n +\lambda \frac{ \lambda   -1 }{n+1}  A_n\\
	&\leq \frac{n}{n+1} (2\lambda-1)A_n 
		=   (2\lambda-1)(2\lambda-1)^{n} 
		= (2\lambda-1)^{n+1}. 
\end{flalign*}
D'autre part le r\'esultat est aussi vrai pour $q_{n+1,n+1}({\lambda}) = {\lambda}^{n+1}$ ce qui termine la preuve.
\end{proof}

\subsection{Comportement asymptotique diagonal.} \label{subse:cad}
Soit $\lambda$ un r\'eel positif, $\beta \in V$ et $\delta \in (V_2 \backslash V_3)$ (i.e. $\delta_1 = 0$ et $\delta_2 \not= 0$). 
Soit $\cL$ le morphisme de $G$ tel que $\cL \big(\(1,0\)\big) = \(\lambda,\beta\)$ et $\cL \big((\(0,1\) \big)= \(0,\delta\)$. 
On a alors $\cL (V) \subset V_2$.

Dans cette section nous \'etudions le comportement de $\cL$ restreint \`a $V$. 
Pour cette raison, les \'el\'ements de $V$ seront directements \'ecrits $(s_1,s_2,\ldots)$ plut\^ot que $\(0,(s_1,s_2,\ldots)\)$.

Si $s = (s_1,\ldots,s_n,\ldots)\in V$, nous notons $s^{(k)} = X^{-k} \cL^k s$ autrement dit
\[
\cL^k (s) = ( \  \underbrace{0\ ,\ \ldots \ ,\ 0}_{k \mbox{ fois}} \ , s_1^{(k)},s^{(k)}_2, \ldots ).
\]
Comme $\delta_1=0$, ces \'el\'ement sont bien d\'efinis et $s^{(k)}_n$ ne d\'epend que de $s_1,s_2,\ldots,s_n$ (voir corollaire~\ref{cor}).

\begin{proposition} \label{prop:Phi}
Soit $\lambda > 0$, $\beta \in V$ et $\delta \in V_2 \backslash V_3$.
Soit $\cL$ le morphisme de $G$ tel que $\cL \big(\(1,0\) \big)= \( \lambda, \beta \)$ et $\cL \big(\(0,1\)\big) = \(0,\delta\)$.

Il existe une fonction $\Phi_\ell: V \rightarrow \RR$ lin\'eaire et qui ne d\'epend que des $\ell$ premi\`eres coordonn\'ees telle que pour tout $s \in V$ on ait
\[
s_\ell^{(n)} = \Phi_\ell(s) \ \delta_2^n \ \lambda^\frac{(n+\ell-1)(n+\ell-2)}{2} + o_n\left( \delta_2^n \lambda^\frac{(n+\ell-1)(n+\ell-2)}{2}  \right).
\]
De plus, $\Phi_\ell(\bold{b}_\ell) = 1$ et
\[
s_1^{(n)} = s_1\ \delta_2^n\ \lambda^{n(n-1)/2}
\quad \text{et} \quad
s_2^{(n)} = s_2 + \sum_{m=0}^{n-1} \lambda^{-n-1} \left( \frac{n (\lambda-1)}{2} + \frac{\delta_3}{\delta_2} \right).
\]
En particulier
\[
\Phi_1(s) = s_1 \quad \text{et} \quad 
\Phi_2(s) = s_2 + \frac{1}{\lambda-1} \left( \frac{1}{2} + \frac{\delta_3}{\delta_2} \right) s_1.
\]
\end{proposition}

\begin{proof}
Comme nous n'utilisons que la restriction de l'endomorphisme $\cL$ sur $V \subset G$ on \'ecrira plus simplement $\cL (v)$ pour $\cL\big( \(0,v\)\big)$.

Nous prouvons le r\'esultat g\'en\'eral par r\'ecurrence sur $\ell$. Nous commen\c{c}ons par \'etudier $s_1^{(n)}$ et $s_2^{(n)}$ pour lesquels tout est explicite.

\bigskip

\noindent
\textbf{Cas $\ell=1$ et $\ell=2$:}

Par la forme explicite de $\cL$ (voir corollaire~\ref{cor:morph}), on a
\[
\cL (s)  = \left( 0,\, \delta_2 s_1,\, \delta_3 s_1 + \lambda \delta_2 s_2,\, s^{(1)}_3, s^{(1)}_4,\, \ldots \right)
\]
Donc $s_1^{(1)} =\delta_2 s_1$ et $s_2^{(1)} = \delta_3 s_1 + \lambda \delta_2 s_2$.
En appliquant \`a nouveau $\cL$ on trouve
\[
\cL^{2}( s) = \left( 0,\, 0 ,\, \lambda \delta_2 s_1^{(1)},\, \left(\binom{\lambda}{2} \delta_2 + \lambda \delta_3\right) s_1^{(1)} + \lambda^2 \delta_2 s_2^{(1)},\, s_3^{(2)},\, s_4^{(2)},\, \ldots \right)
\]
Nous trouvons donc $s_1^{(2)} = \lambda \delta_2^2 s_1$ et $s_2^{(2)} = \left( \binom{\lambda}{2} \delta_2 + \lambda \delta_3\right) s_1^{(1)} + \lambda^2 \delta_2 s_2^{(1)}$.
On montre alors simplement par r\'ecurrence que
\begin{equation} \label{eq:rec_si_1_et_2}
s_1^{(n+1)} = \delta_2\, q_{n,n}\, s_1^{(n)} = \delta_2 \lambda^n s_1^{(n)}
\quad \text{et} \quad
s_2^{(n+1)} = (q_{n,n+1}\, \delta_2 + q_{n,n}\, \delta_3) s_1^{(n)} + q_{n+1,n+1}\, \delta_2\, s_2^{(n)}.
\end{equation}
Ainsi $s^{(n)}_1 =s_1\, \delta_2^n\, \lambda^\frac{n(n-1)}{2}$. Pour l'\'etude de $s_2^{(n)}$, posons
\[
y_n = \frac{s_2^{(n)}}{\delta_2^n \lambda^{n(n+1)/2}}.
\]
En divisant la relation de r\'ecurrence pour $s_2^{(n)}$ on obtient
\[
y_{n+1}
= y_n + \frac{1}{\delta_2^{n+1} \lambda^{(n+1)(n+2)/2}}
\big{(} q_{n,n+1} \delta_2+q_{n,n} \delta_3     \big{)} s^{(n)}_1
= y_n + \frac{q_{n,n+1} \delta_2 + q_{n,n} \delta_3}{\delta_2 \lambda^{2n+1}} s_1 
\]
o\`u on a utilis\'e $\frac{(n+1)(n+2)}{2}-(n+1)=\frac{n(n+1)}{2}$ et $\frac{(n+1)(n+2)}{2}-\frac{n(n-1)}{2}=2n+1$.
Maintenant, nous savons que $q_{n,n}(\lambda) = \lambda^n$ et $q_{n,n+1}(\lambda) = \frac{n}{2}  (\lambda-1) \lambda^n$
Nous trouvons alors
\[
\frac{q_{n,n+1} \delta_2+q_{n,n} \delta_3  }{\delta_2 \lambda^{2n+1}}
=
\frac{1}{\delta_2 \lambda^{2n+1}}  \left( \frac{n}{2}  (\lambda-1) \lambda^n \delta_2+ \lambda^n \delta_3 \right)
=
\frac{1}{\lambda^{n+1}} \left( \frac{n}{2}  (\lambda-1) +  \frac{\delta_3}{\delta_2}  \right).
\]
La s\'erie de terme g\'en\'eral $\displaystyle \lambda^{-n-1}\left( \frac{n}{2}  (\lambda-1) +  \frac{\delta_3}{\delta_2}  \right)$ est sommable et donc la suite $(y_n)_n$ converge vers 
\[
\Phi_2(\bold{s}) = s_2 + s_1 \sum_{m=0}^\infty \lambda^{-n-1} \left(\frac{n}{2} (\lambda-1) + \frac{\delta_3}{\delta_2} \right)
= s_2 + \frac{1}{\lambda-1} \left( \frac{1}{2} + \frac{\delta_3}{\delta_2} \right) s_1.
\]

\noindent
\textbf{Passage de $\ell-1$ \`a $\ell$.}  \label{subse:ellmoinsun}
Pour le cas g\'en\'eral, on \'etablit d'abord une relation de r\'ecurrence similaire \`a~\eqref{eq:rec_si_1_et_2} qui exprime $s_\ell^{(n)}$ sous la forme d'une s\'erie. Nous utilisons ensuite l'estimation asymptotique des $q_{n,n+i}$ pour conclure qu'elle converge.

D'apr\`es la forme explicite de $\cL$ (voir corollaire~\ref{cor:morph}), on trouve
\begin{equation} \label{eq:moritell}
s_\ell^{(n+1)} =
\sum \limits_{j=1} ^{\ell} \sum \limits_{i=1}^{\ell-j+1} q_{j+n-1,n+\ell-i}  \delta_i  \ s_j^{(n)} .
\end{equation}
Posons
\[
y_n = \frac{s_\ell^{(n)}}{\delta_2^n \lambda^{(n+\ell-1)(n+\ell-2)/2}}.
\]
Le coefficient de $s_\ell^{(n)}$ dans la relation \eqref{eq:moritell} est $q_{n+\ell-1,n+\ell-1}=\lambda^{n+\ell-1}$.
En remarquant que $\binom{n+\ell}{2} -\binom{n+\ell-1}{2} = n+\ell-1$
et en divisant la relation \eqref{eq:moritell} par $\delta_2^{n+1} \lambda^{ \binom{n+\ell}{2} }$, nous trouvons :
\[
y_{n+1}
= y_n + \frac{1}{\delta_2^{n+1}    \lambda^{ \binom{n+\ell}{2} }   }
\sum \limits_{j=1} ^{\ell-1} \left( \sum \limits_{i=1}^{\ell-j+1} q_{j+n-1,n+\ell-i} \delta_i  \right) s_j^{(n)} .
\]
Nous voulons montrer que la s\'erie de terme g\'en\'eral $y_{n+1} - y_n $ converge.
Cette s\'erie est elle-m{\^e}me d'une somme d'un nombre fini de termes.
Nous allons \'etudier chacun de ces termes et montrer qu'ils sont sommables.
Fixons  $j\in \{1,\ldots,\ell-1\}$ et $i\in\{1,\ldots,\ell-j+1\}$ et \'etudions le comportement asymptotique de la suite
\[
z_n =  \frac{1}{\delta_2^{n+1} \lambda^{(n+\ell-1+1)(n+\ell-2+1)/2}} q_{j+n-1,n+\ell-i} \delta_i  s_j^{(n)}.
\]
D'apr\`es l'hypoth\`ese de r\'ecurrence, $z_n$ est \'egal \`a :
\[
q_{j+n-1,n+\ell-i} \delta_i 
\left( \Phi_j(s)
  \frac{\delta_2^n \lambda^{(n+j-1)(n+j-2)} }{\delta_2^{n+1} \lambda^{(n+\ell-1+1)(n+\ell-2+1)/2}}
  +  o\left( \frac{ \delta_2^n \lambda^\frac{(n+j-1)(n+j-2)}{2} }{\delta_2^{n+1} \lambda^{(n+\ell-1+1)(n+\ell-2+1)/2}} \right) \right) .
\]
Nous pouvons imm\'ediatement simplifier cette expression en remarquant que
\[
\binom{n+\ell}{2} - \binom{n+j-1}{2} = n(\ell-j+1)+j' \mbox{ avec }j' = (\ell^2 -\ell -j^2+2j+j-2 )/2.
\]
Donc
\[
z_n = 
q_{j+n-1,n+\ell-i} \delta_i   \left(  \Phi_j(s)  \frac{  1  }{\delta_2  \lambda^{n(\ell-j+1)+j'} }
  +  o\left(  \frac{  1  }{\delta_2  \lambda^{n(\ell-j+1)+j'} } \right) \right).
\]
Nous avons vu dans la section \ref{subse:qij} le comportement asymptotique de la suite $ q_{j+n-1,n+\ell-i}(\lambda) $.
D'apr\`es la proposition~\ref{pr:q}, nous trouvons :
\[
q_{j+n-1,n+\ell-i}(\lambda)
=q_{n+\ell-i-(\ell-j+1),n+\ell-i}(\lambda)
\underset{n\to+\infty}{\sim}
\frac{(n+\ell-i)^{\ell-j+1}}{2^{\ell-j+1} (\ell-j+1)!} \left( \frac{\lambda-1}{\lambda}\right) ^{\ell-j+1}  \lambda^{n+\ell-i} .
\]
Donc, nous obtenons l'\'equivalent suivant :
\[
q_{j+n-1,n+\ell-i}(\lambda)  \underset{n}{\sim}  A  n^{\ell-j+1} \lambda^{n}
 \quad \mbox{ avec } \quad A = \frac{ \lambda^{\ell-i}  }{2^{\ell-j+1} (\ell-j+1)!} \left( \frac{\lambda-1}{\lambda}\right) ^{ \ell-j+1}  .
\]
Soit encore,
\[
z_n
=\Phi_j(s)  \frac{  A  n^{\ell-j+1} \lambda^{n}  }{\delta_2  \lambda^{n(\ell-j+1)+j'} }
+ o\left(  \frac{  A  n^{\ell-j+1} \lambda^{n}  }{\delta_2  \lambda^{n(\ell-j+1)+j'} } \right) \
=\Phi_j(s)  \frac{  A  n^{\ell-j+1}  }{\delta_2  \lambda^{n(\ell-j)+j'} }
+  o\left(  \frac{  A  n^{\ell-j+1}   }{\delta_2  \lambda^{n(\ell-j)+j'} } \right) .
\]
La suite $(y_n)_n$ converge donc vers un r\'eel  $\Phi_{\ell}(s)$.
Puisque $y_n$ est la somme de $s_\ell$ et d'une combinaison lin\'eaire de $\Phi_j(s)$, la fonction $\Phi_\ell(s)$ est bien lin\'eaire en $s_1,\ldots,s_\ell$, et le coefficient de $s_\ell$ est $1$.

Ce qui finit la d\'emonstration du r\'esultat. 
\end{proof}

\subsection{Asymptotiques des sommes de Birkhoff et des polyn\^omes d'approximation.} \label{se:po}
On montre dans cette section deux applications de la proposition \ref{prop:Phi}
qui nous permettront de contr\^oler la taille des sommes de Birkhoff et des polyn\^omes d'approximation $R_{i,j}$.

\begin{corollaire} \label{cor:limite_diff_birkhoff} \label{cor:limite_difference_Rij}
Soit $\sigma$ une substitution $\lambda$-uniforme telle que $\delta_1(\sigma) = 0$ et $\delta_2 = \delta_2(\sigma) \not= 0$.
Soit $\phi: \{a,b\} \rightarrow \RR$ telle que $\phi(a) \not= \phi(b)$.
Soit $R_{i,j}$ les polyn\^omes d'approximation de la section~\ref{subsec:approximation_polynomiale}. 
Alors pour tout mot fini $\u_0 \u_1 \ldots \u_{m-1}$ et tout entier $k \geq 0$,
\[
\lim_{n \to \infty} \frac{S^{(n+k)}_{m \lambda^n}(\phi, \sigma^n(\u)) - R_{n,n+k}(m)}{\delta_2^n\  \lambda^{(n+k-1)(n+k-2)/2}}
\]
converge vers $\Phi_k(S_{m}(\phi,\u))$ lorsque $n$ tend vers $+\infty$.

Pour tout entier $k$ et tout r\'eel $x$,
\[
\frac{R_{n+k,n+k} (x)  - R_{n,n+k}(\lambda ^kx)  }{\delta_2^n\ \lambda^{(n+k-1)(n+k-2)/2}}
\]
converge vers $\Phi_k \left( (R_{k,1}(x), R_{k,2}(x), \ldots) \right)$ lorsque $n$ tend vers $+\infty$.
\end{corollaire}

\begin{proof}[Preuve du corollaire~\ref{cor:limite_difference_Rij}]
Par construction, on a
\[
\cL^n \circ \pi(\u_0 \ldots \u_{m-1}) = \big\(m \lambda^n, S^{(n)}_{m \lambda^n}(\phi,\sigma^n(\u))\big\)
\]
et
\[
\cL^n\big( \(m,0\)\big) = \big\(m \lambda^n, (R_{1,n}(m), R_{2,n}(m), \ldots) \big\).
\]
Maintenant, remarquons que dans $G$ on a $\(z,s\)^{-1}\ \(z,t\) = \(0,s-t\)$. En particulier
\begin{align*}
\cL^n \big(\(m,0\)^{-1} \pi(\u)\big)
&= \big\(m \lambda^n, (R_{n,1}(m), R_{n,2}(m), \ldots)\big\)^{-1} \big\(m \lambda^n, S_{m \lambda^n}(\phi, \sigma^n(\u)) \big\) \\
&= \big\(0, S_{m \lambda^n}(\phi, \sigma^n(\u)) - (R_{n,1}(m), R_{n,2}(m), \ldots) \big\).
\end{align*}
Il suffit alors, d'appliquer la proposition~\ref{prop:Phi}.

On montre la seconde partie de mani\`ere similaire.
Soit $x$ un r\'eel.
Alors :
\begin{align*}
\cL^n \circ  \cL^{k} \big( \(x,0\) \big) &= \big\(\lambda^{n+k} x,(R_{n+k,1}(x),\ldots,R_{n+k,j}(x),\ldots) \big\)\\
\cL^n  \big( \(\lambda^k x,0\)\big) &=\big\(\lambda^{n+k} x,(R_{n,1}(\lambda^k x),\ldots,R_{n,j}(\lambda^k x),\ldots)\big\).
\end{align*}
On peut alors \'ecrire
\begin{align*}
\cL^n \big(\(0, (R_{k,1}(x), R_{k,2}(x), \ldots) \) \big)&= 
\cL^n \circ \cL^k \big(\(x,0\)\ \(\lambda^k x, 0\)^{-1}\big) \\
&= \big\(0, (R_{n+k,1}(x) - R_{n,1}(\lambda^kx),\,R_{n+k,2}(x) - R_{n,2}(\lambda^k x), \ldots)\big\).
\end{align*}
Il suffit d'appliquer la proposition~\ref{prop:Phi} pour conclure.
\qedhere
\end{proof}

\subsection{Coefficients $c_\ell$.} \label{se:convcl}
Soient $\lambda$ un entier, $\beta \in V$ et $\delta \in V_2 \backslash V_3$. Comme auparavant, on leur associe un endomorphisme $\cL$. On suppose ici que $\beta$ et $\delta$ ont des coefficients nuls \`a partir du rang $\lambda+1$. C'est le cas lorsque $\cL$ est associ\'e \`a une paire $(\phi,\sigma)$ de fonctions $\phi: \{a,b\}^* \rightarrow \RR$ et d'une substitution $\sigma$ de longueur constante $\lambda$.
On notera $||\delta||=\sup\{|\delta_i|;1\leq i\}$.

Nous reprenons les polyn\^omes $R_{i,j}$ et les nombres $c_\ell$ introduit dans la proposition~\ref{prop:R_et_c}.
Rappelons qu'ils sont d\'efinis 
\[
\cL^i \big(\(x,0\) \big)= \left\( \lambda^i x,\, (R_{i,1}(x), R_{i,2}(x), R_{i,3}(x), \ldots) \right\).
\]
et qu'on pose alors
\[
c_i = R_{i+1,i+1} \left( \frac{1}{\lambda^{i+1}} \right).
\]
Nous allons montrer le r\'esultat suivant :
\begin{proposition} \label{prop:cellconv}
Soit $c_i$ comme ci-dessus. Alors, la suite $\left( c_\ell \cdot \delta_2^{-\ell} \lambda^{-(\ell-1)(\ell-2)/2} \right)_{\ell \geq 0}$ converge.
\end{proposition}

\begin{proof}
Fixons un entier $\ell \geq \lambda$.
Commen{\c c}ons par rappeler la relation de la proposition~\ref{prop:R_et_c} :
\[
\cL \big( \(1, (c_0, c_1, c_2, \ldots) \) \big)= \(1, (c_0, c_1, c_2, \ldots)\)^{\lambda}.
\]
D'apr\`es la forme explicite de l'endomorphisme $\cL$ (corollaire~\ref{cor:morph}), nous trouvons pour $\ell > \lambda$ :
\begin{align*}
\lambda c_\ell+\binom{\lambda}{2} c_{\ell-1}+\cdots+\binom{\lambda}{\ell+1} c_0 &=
\lambda^{\ell-1}\delta_2 c_{\ell-1}  + \sum \limits_{j=1} ^{\ell-1}   \sum \limits_{i=1}^{\ell-j+2} q_{j-1,\ell+1-i}  \delta_i c_{j-1}  \\
& = \lambda^{\ell-1}\delta_2 c_{\ell-1}  + \sum \limits_{j=2} ^{\ell-1}   \sum \limits_{i=2}^{\ell-j+2} q_{j-1,\ell+1-i}  \delta_i c_{j-1}.
\end{align*}
La seconde \'egalit\'e est d\^ue au fait que $q_{0,i}=0$ si $i > 0$.
La double somme se majore avec le lemme~\ref{le:techniquepourq} :
\[
\begin{array}{ll}
\left| \sum \limits_{j=1} ^{\ell-1}   \sum \limits_{i=2}^{\ell-j+2} q_{j-1,\ell+1-i}  \delta_i c_{j-1}\right|
	&\leq ||\delta|| \cdot \sup \limits_{1\leq j \leq \ell-2} |c_j| \cdot  \sum \limits_{j=1} ^{\ell-1}   \sum \limits_{i=2}^{\ell-j+2} q_{j-1,\ell+1-i} \\
	& \leq ||\delta|| \cdot \sup \limits_{1\leq j \leq \ell-2} |c_j| \cdot  \left( \sum \limits_{k=0} ^{\ell-2} q_{k,\ell-1}+ \sum \limits_{r=0} ^{\ell-2} \sum \limits_{k=0} ^r q_{k,r} \right) \\
	& \leq  ||\delta|| \cdot \sup \limits_{1\leq j \leq \ell-2} |c_j| \cdot  \left( \ell-1+\binom{\ell}{2} \right)  (2\lambda-1)^{\ell-1} \\
	& \leq  ||\delta|| \cdot \sup \limits_{1\leq j \leq \ell-2} |c_j| \cdot   \binom{\ell+1}{2}   (2\lambda-1)^{\ell-1} .
\end{array}
\]

Pour l'autre partie, on utilise simplement :
\[
\left| \binom{\lambda}{2} c_{\ell-1}+\cdots+\binom{\lambda}{\ell+1} c_0 \right|
\leq
2^\lambda\ \sup_{0 \leq j \leq \ell-1} |c_j|.
\]
Nous trouvons donc :
\begin{equation} \label{eq:majoration_convergence_cell_1}
| \lambda c_\ell - \lambda^{\ell-1}\delta_2 c_{\ell-1}  | \leq   ||\delta|| \, \binom{\ell+1}{2} \, (2 \lambda-1)^{\ell-1}  \sup \limits_{1\leq j \leq \ell-2} |c_j| + \sup \limits_{1\leq j \leq \ell-1} |c_j| \cdot 2^\lambda.
\end{equation}
Notons $\widetilde{c_j} = c_j \cdot \delta_2^{-j} \lambda^{-(j-1)(j-2)/2}$. Nous avons
\[
\left\lvert \frac{c_j}{\delta_2^\ell \lambda^{(\ell-1)(\ell-2)/2} } \right\rvert =
\frac{|\widetilde{c_j}|}{|\delta_2|^{\ell-j} \lambda^{(\ell-1)(\ell-2)/2 - (j-1)(j-2)/2}  }
\leq \frac{|\widetilde{c_j}|}{\lambda^{(\ell-1)(\ell-2)/2 - (j-1)(j-2)/2}}.
\]
En divisant la relation \eqref{eq:majoration_convergence_cell_1} par $\delta_2^\ell \lambda^{(\ell-1)(\ell-2)/2}$, on obtient alors
\begin{align*}
|\widetilde{c_\ell} - \widetilde{c_{\ell-1}}| &\leq
\|\delta\| \binom{\ell+1}{2} \frac{(2\lambda-1)^{\ell-1}}{\lambda^{2\ell-1}} \sup_{0 \leq j \leq \ell-2} \widetilde{|c_j|} + \frac{2^\lambda}{\lambda^{\ell-1}} \sup_{0 \leq j \leq \ell - 1} |\widetilde{c_j}|
\end{align*}
car pour $1\leq j \leq \ell-2$ : $\binom{\ell-1}{2}-\binom{j-1}{2}= \ell + \ell-1 +\binom{\ell-3}{2}-\binom{j-1}{2} \geq 2 \ell-1$ et $\binom{\ell-1}{2} - \binom{\ell-2}{2} = \ell-1$.

Maintenant, si $\theta$ est tel que $\frac{2\lambda-1}{\lambda^2} < \theta < 1$ alors pour $\ell$ assez grand on a
\[
| \tilde c_\ell - \tilde c_{\ell-1}  | \leq  \theta^\ell \sup \limits_{1\leq j \leq \ell-1} |\tilde c_j|.
\]
On peut alors conclure en utilisant le lemme~\ref{lem:convergence_suite}.
\end{proof}

\section{Preuve du th\'eor\`eme} \label{se:preuve}
On se donne une substitution $\sigma$ de longueur constante $\lambda$ telle que $\delta_1(\sigma) = 0$, $\delta_2(\sigma) \not= 0$ et admettant un point fixe $\w=\w_0 \cdots\w_n\cdots$.
Soit $\phi: \{a,b\} \rightarrow \RR$ une fonction non constante.
Nous allons d\'efinir une suite de fonctions en escalier $(f^{(\ell)})_{\ell\geq 1}$ qui convergera vers la fonction $f_{\w,\phi}=f_\w$ du th\'eor\`eme.

On reprend les constantes $c_i$ de la section~\ref{subsec:approximation_polynomiale} et la suite de fonctions $\psi_i: K_\sigma \rightarrow \RR$ de la section~\ref{sec:cobords} telles que
\[
\psi_0 = - \int \phi\,, \quad \phi + \psi_0 = \psi_1 \circ T - \psi_1
\quad\text{et} \quad
\quad \psi_n = \psi_{n+1} \circ T - \psi_{n+1}.
\]
Et $c_i$ v\'erifient $c_i = - \psi_i(\w)$.
Rappelons qu'on a un lien direct entre ces fonctions et les sommes de Birkhoff via
\[
\psi_i(T^n \w) = S^{(i)}_n(\phi,\w) - p_i(n) \quad \text{o\`u} \quad p_i(n) = \sum_{k=0}^i \binom{n}{i-k} c_k.
\]

On d\'efinit une suite de fonctions $f^{(i)}: \RR \rightarrow \RR$ par
\[
f_\w^{(i)}(x) = \frac{\psi_i(T^m \w)}{\delta_2^{i-1} \lambda^{(i-1)(i-2)/2}} \qquad \text{avec $m = \left\lfloor \lambda^{i-1} x \right\rfloor$}.
\]

Pour des entiers $i \geq 0$ et $0 \leq m < i$ et une lettre $\alpha$, on rappelle la notation $\Cyl(i,m,\alpha) = T^m \sigma^i([\alpha])$ (voir section~\ref{sec:cobords}).
On rappelle que la fonction $\psi_i$ est constante sur chaque $\Cyl(i,m,\alpha)$ et on notera par $\psi(i,m,\alpha)$ cette valeur.
On d\'efinit deux fonctions $f^{(i)}_a, f^{(i)}_b: [0,\lambda] \rightarrow \RR$ par
\[
f^{(i)}_a(x) = \frac{\psi(i, \lfloor \lambda^{i-1} x \rfloor, a)}{\delta_2^{i-1} \lambda^{(i-1)(i-2)/2}}
\qquad \text{et} \qquad
f^{(i)}_b(x) = \frac{\psi(i, \lfloor \lambda^{i-1} x \rfloor, b)}{\delta_2^{i-1} \lambda^{(i-1)(i-2)/2}}.
\]
On a alors pour $x \in [0,\lambda]$ et $n$ entier positif la relation:
\[
f_\w^{(i)}(x + \lambda n) = f^{(i)}_{\w_n}(x).
\]
Comme pour tout $\u \in \Cyl(i, 0, \alpha)$ on a $\psi_i(T^{\lambda^i} \u) = \psi_i(\u)$ on en d\'eduit que
\[
f_\w^{(i)}(0) = f^{(i)}_a(0) = f^{(i)}_b(0) = \frac{-c_i}{\delta_2^i \lambda^{(i-1)(i-2)/2}}.
\] 
D'autre part, comme pour tout mot $\u$ on a $\psi_i(T \u) - \psi_i(\u) = \psi_{i-1}(\u)$ on obtient l'\'equation
\begin{equation} \label{eqiet}
 f_\w^{(i)} \left( x + \frac{1}{\lambda^{i-1}}  \right) - f_\w^{(i)}(x) = \frac{1}{\delta_2 \lambda^{i-2}}\ f_\w^{(i-1)}(\lambda x).
\end{equation}

Tout le reste de cette section sera d\'edi\'ee \`a la preuve du r\'esultat suivant qui pr\'ecise l'\'enonc\'e du th\'eor\`eme~\ref{thm:construction} de l'introduction.
\begin{theoreme} \label{thm:thm_principal_precis}
Soit $\sigma$, $f_\w^{(i)}$, $f^{(i)}_a$ et $f^{(i)}_b$ comme ci-dessus. Les suites $f_\w^{(i)}$, $f^{(i)}_a$ et $f^{(i)}_b$ convergent uniform\'ement vers des fonctions continues $f_\w$, $f_a$ et $f_b$ sur respectivement $\RR^+$, $[0,\lambda]$ et $[0,\lambda]$. De plus:
\begin{enumerate}
\item $f_a$ et $f_b$ ne sont pas nulles, \label{item:non_nullite}
\item $f_a(0) = f_b(0) = \lim -c_i / \delta_2^{i-1} \lambda^{(i-1)(i-2)/2}$, \label{item:valeur_en_zero}
\item pour $x \in [0,\lambda]$ et $n$ entier: $f_\w(x + \lambda n) = f_{\w_n}(x)$, \label{item:lin_rec}
\item $f_\w$ est solution de l'\'equation $(E_{\lambda,\delta_2(\sigma)})$. \label{item:equation_fonctionnelle}
\end{enumerate}
\end{theoreme}
Une fois la convergence \'etablie, les items~\ref{item:valeur_en_zero}, \ref{item:lin_rec} et~\ref{item:equation_fonctionnelle} proviennent directement de la construction.

\subsection{Convergence aux points $\lambda$-adiques.}
On appelle un point \emph{$\lambda$-adique} un r\'eel $x$ tel qu'il existe un entier $n$ tel que $\lambda^n x$ soit entier.

\begin{lemme} \label{lem:convergence_points_rat}
Les fonctions $f_\w^{(i)}$, $f_a^{(i)}$ et $f_b^{(i)}$ convergent en tout point $\lambda$-adique. 

De plus, si on note respectivement $\lambda_a$ et $\lambda_b$ le nombre de $a$ et $b$ dans $\sigma(a)$ alors pour tout entier $0 \leq k < \lambda$ on a pour tout $i \geq 1$
\[
f_a^{(i)}(k+1) - f_a^{(i)}(k) = \left\{ \begin{array}{ll}
\frac{\lambda_b}{\lambda} (\phi(a) - \phi(b)) & \text{si $(\sigma(a))_k = a$} \\
\frac{\lambda_a}{\lambda} (\phi(b) - \phi(a)) & \text{si $(\sigma(a))_k = b$}.
\end{array} \right.
\]
et
\[
f_b^{(i)}(k+1) - f_b^{(i)}(k) = \left\{ \begin{array}{ll}
\frac{\lambda_b}{\lambda} (\phi(a) - \phi(b)) & \text{si $(\sigma(b))_k = a$} \\
\frac{\lambda_a}{\lambda} (\phi(b) - \phi(a)) & \text{si $(\sigma(b))_k = b$}.
\end{array} \right.
\]
\end{lemme}

\begin{proof}
Soit $x$ un r\'eel tel que $\lambda^{i_0-1} x = m \in \NN$.
On note $k_j = \delta_2^{j-1} \lambda^{(j-1)(j-2)/2}$ les coefficients de renomalisation.
Alors pour tout $i \geq 0$ :
\[
f^{(i+i_0)}(x) = \dfrac{\psi_{i+i_0}(T^{m\lambda^i} \w)}{k_{i+i_0}} \ ; \quad 
f^{(i+i_0)}_a(x) = \dfrac{\psi(i+i_0,m\lambda^i, a)}{k_{i+i_0}} \quad \mbox{ et } \quad
f^{(i+i_0)}_b(x) = \dfrac{\psi(i+i_0, m \lambda ^i, b)}{k_{i+i_0}} .
\]
Nous allons montrer que la suite $\big(f^{(i)}(x)\big)_i $ converge, ce qui prouvera le r\'esultat.
\[
\begin{array}{ll}
f^{(i+i_0)}(x) 
&= \dfrac{\psi_{i+i_0}(T^{m\lambda^i} \w)}{k_{i+i_0}}	
	= \dfrac{S^{(i+i_0)}_{m\lambda^i} (\phi,\w) - p_{i+i_0}(m\lambda^i)}{k_{i+i_0}} \\
& =\dfrac{1}{\delta_2^{i_0-1}} 
 \dfrac{	 S^{(i+i_0)}_{ m {\lambda} ^i}(\phi,\w) - p_{i+i_0}( m \lambda^i)  }{ \displaystyle \delta_2^n  {\lambda}^{\binom{i+i_0}{2}}}  \\
&= 
 \dfrac{ \displaystyle	 S^{(i+i_0)}_{m {\lambda} ^i}(\phi,\w) - R_{i+i_0, i+i_0} (m) - c_{i+i_0}}
{k_{i+i_0}}    \mbox{ d'apr\`es l'\'equation \eqref{eq:def_p}} \\
&=
\dfrac{ S^{(i+i_0)}_{m{\lambda} ^i}(\phi,\w) - R_{i,i+i_0}\left( m  \right)  }{  k_{i+i_0}  }
+ \dfrac{ R_{i,i+i_0}\left( x  \right)   - R_{i+i_0,i+i_0}\left( x \right)} {k_{i+i_0}}
- \dfrac{c_{i+i_0} } {k_{i+i_0}}
\end{array}
\]
Les trois termes ci-dessus convergent d'apr\`es le corollaire~\ref{cor:limite_difference_Rij} et la proposition~\ref{prop:cellconv}.

Dans le cas particulier $i_0=1$, on obtient une formule explicite et qui ne d\'epend pas de $i$. Ceci est directement reli\'e au fait que, dans la Proposition~\ref{prop:Phi}, les valeurs de $s_1^{(n)}$ renormalis\'ees sont constantes.
Comme dans la section~\ref{sec:cobords}, on note $\Cyl(0,0,a)$ le cylindre de $K_{\sigma,\ZZ}$ des mots $\cdots \u_{-1} \cdot \u_0 \u_1 \cdots$ tels que $\u_0 = a$. D'apr\`es les formules pour les sommes de Birkhoff~\eqref{eq:cobord} et~\eqref{eq:cobord_bis}, on sait que pour $\u \in \sigma^i(\Cyl(0,0,a)) = \Cyl(i,0,a)$ on a
\[
\psi_{i+1}(T^{\lambda^i} \u) - \psi_{i+1}(\u) = S^{(1)}(\psi_i, \sigma^i(a)).
\]
De la m\^eme fa\c{c}on pour $\u \in \Cyl(i,0,b)$ en rempla\c{c}ant $a$ par $b$ ci-dessus. En particulier, cette quantit\'e ne d\'epend pas de $\u$. D'autre part, nous savons par la proposition~\ref{prop:Phi} que pour tout $i \geq 1$ que
\[
S^{(1)}_{\lambda^i}(\psi_i,\sigma^i(a)) - S^{(1)}_{\lambda^i}(\psi_i,\sigma^i(b)) = S^{(i+1)}_{\lambda^i}(\phi,\sigma^i(a)) - S^{(i+1)}_{\lambda^i}(\phi,\sigma^i(b))
= (\phi(a) - \phi(b)) \delta_2^i \lambda^{(i-1)(i-2)/2}.
\]
De plus, $\psi_{i+1}$ est de moyenne nulle et donc
\[
\lambda_a S^{(1)}(\psi_i, \sigma^i(a)) + \lambda_b S^{(1)}(\psi_i, \sigma^i(b)) = 0.
\]
En recollant ces deux \'egalit\'es, on obtient le r\'esultat annonc\'e.
\end{proof}

\begin{figure}[H]
\begin{center}
\includegraphics{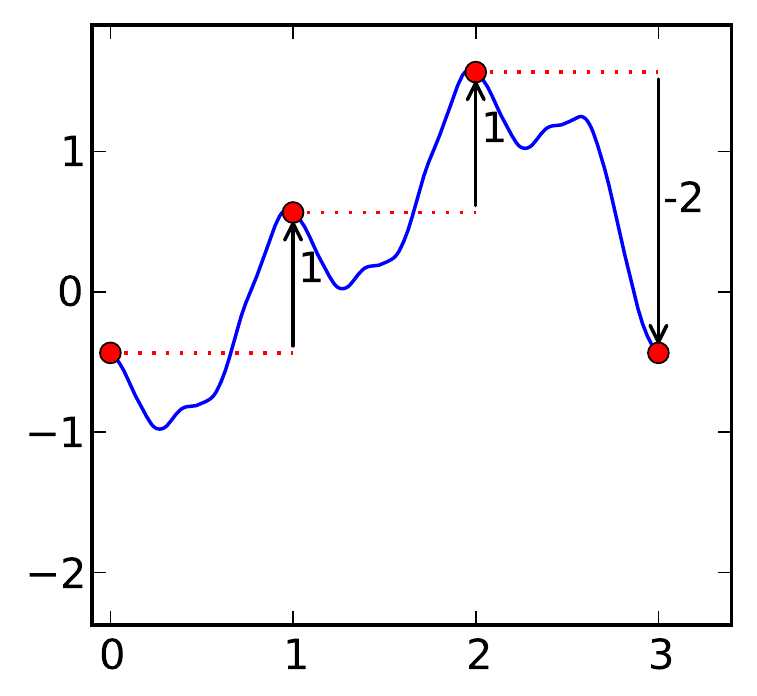} \vspace{.5cm} \includegraphics{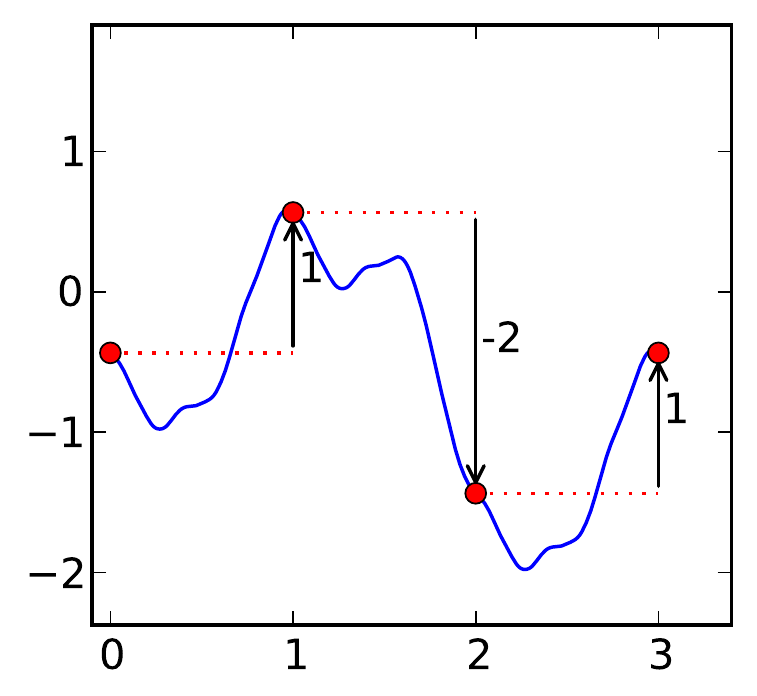}
\end{center}
\caption{Les fonctions $f_a$ et $f_b$ pour $a \mapsto aab, b \mapsto aba$ et pour $\phi = \chi_a - 2 \chi_b$.}
\label{fig:fa_fb_aab_aba}
\end{figure}

\begin{figure}[H]
\begin{center}
\includegraphics{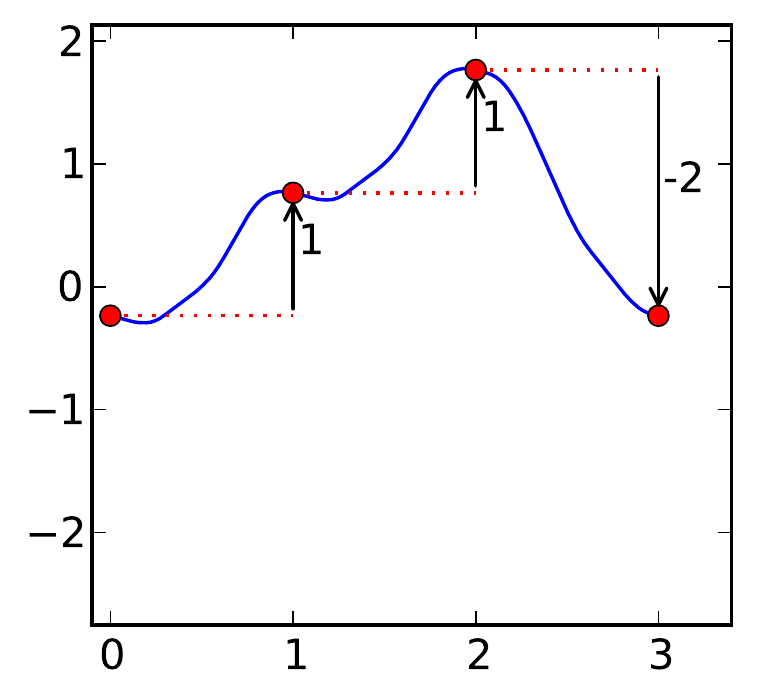} \vspace{.5cm} \includegraphics{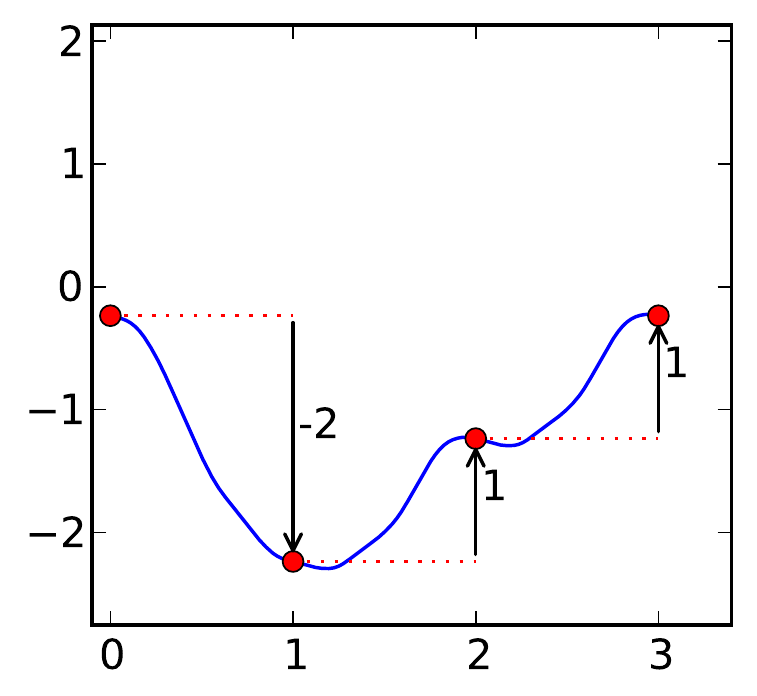}
\end{center}
\caption{Les fonctions $f_a$ et $f_b$ pour $a \mapsto aab, b \mapsto baa$ et pour $\phi = \chi_a - 2 \chi_b$.}
\label{fig:fa_fb_aab_baa}
\end{figure}

\begin{figure}[H]
\begin{center}
\includegraphics{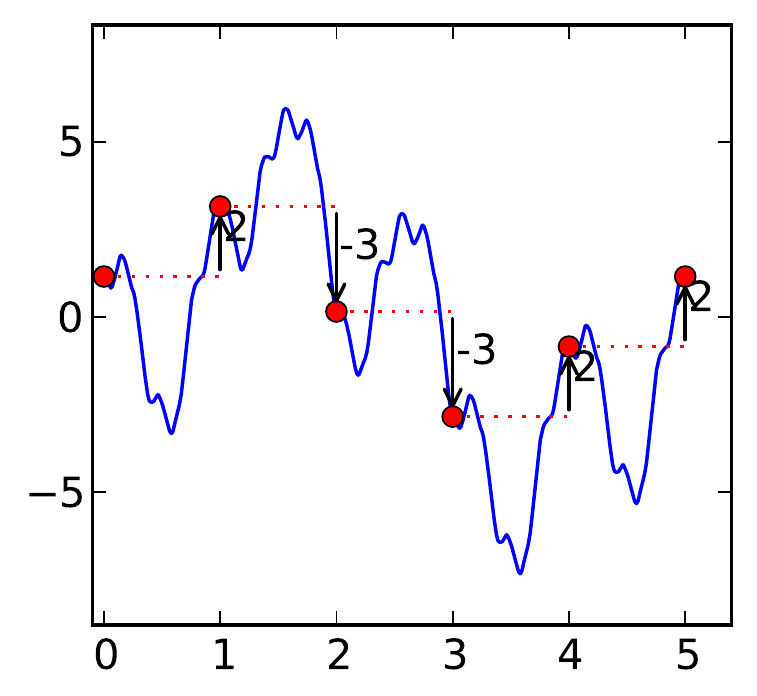} \vspace{.5cm} \includegraphics{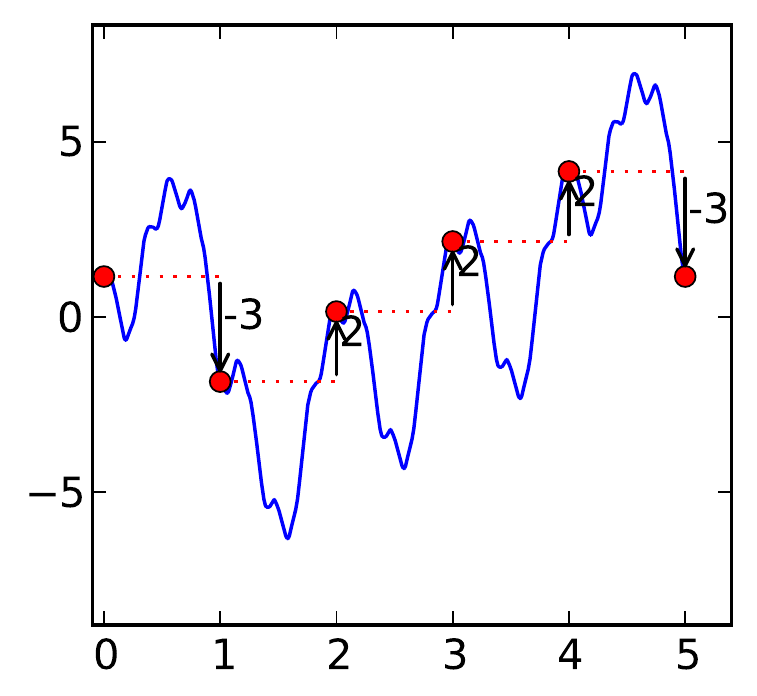}
\end{center}
\caption{Les fonctions $f_a$ et $f_b$ pour $a \mapsto abbaa, b \mapsto baaab$ et pour $\phi = 2\chi_a - 3\chi_b$.}
\label{fig:fa_fb_abbaa_baaab}
\end{figure}

Ainsi les valeurs aux points entiers sont constantes le long de la suite $f^{(i)}$. Ce n'est plus vrai pour les points de la forme $m/\lambda^2$ du fait de la formule pour $s_2^{(n)}$ dans la proposition~\ref{prop:Phi}. Lorsqu'on choisit $\phi = \lambda_b \chi_a - \lambda_a \chi_b$ alors la diff\'erence $f_a(k+1)-f_a(k)$ vaut soit $\lambda_b$ si la $k$-\`eme lettre de $\sigma(a)$ est un $a$ et vaut $-\lambda_a$ sinon. On peut voir ce fait sur les figures~\ref{fig:fa_fb_aab_aba}, \ref{fig:fa_fb_aab_baa} et \ref{fig:fa_fb_abbaa_baaab}.

\subsection{Borne uniforme.} \label{subse:bu}
Nous attaquons la partie la plus d\'elicate de la preuve: montrer que les fonctions $f^{(i)}$ sont born\'ees. La preuve est tr\`es similaire \`a celle de la proposition~\ref{prop:cellconv} pour la convergence des coefficients $c_\ell$.
\begin{proposition} \label{prop:borne_uniforme}
La suite de fonctions $f^{(i)}$ est uniform\'ement born\'ee.
\end{proposition}
En cons\'equence, puisque $f^{(i)}$ est une "concat\'enation" de fonctions $f_a^{(i)}$ et $f_b^{(i)}$, ces suites de fonctions sont \'egalement uniform\'ement born\'ees.

\begin{proof}
Fixons un entier $\ell > \lambda+1$. Nous noterons $q_{i,n}=q_{i,n}(\lambda) $ et $k_j = \delta_2^j \lambda^{(j-1)(j-2)/2}$.

Nous commen{\c c}ons par relier les valeurs prises par la fonction $f^{(\ell)}$, \`a celles prises par les fonctions $f^{(1)},\ldots,f^{(\ell-1)}$.

D'apr\`es le corollaire~\ref{cor:morph} et le travail fait dans la partie \ref{subse:induit}, pour tout entier $r$ :
\begin{equation} \label{eq:pout1}
S_{r\lambda} ^{(\ell)} =   \sum \limits_{j=2} ^{\ell-1}   \sum \limits_{i=2}^{\ell-j+1} q_{j-1,\ell-i} \left( \delta_i   \cdot  S_r^{(j)} +\beta_i \binom{r}{j}\right).
\end{equation}

Puisque pour tout entier $r$, et tout entier $j\leq \ell$, $S_{i \lambda^\ell} ^{(j)} = p_j(i \lambda^\ell)-c_j$, nous en d\'eduisons que
\begin{equation} \label{eq:pout2}
p_\ell( r\lambda )-c_\ell =  \sum \limits_{j=2} ^{\ell-1}   \sum \limits_{i=2}^{\ell-j+1} q_{j-1,\ell-i} 
 	\left( \delta_i   \cdot  \left( p_j( r)-c_j \right) +\beta_i \binom{r}{j}\right).
\end{equation}

Rappelons maintenant la relation v\'erifi\'ee par les coefficient $c_i$ :
\[
c_\ell = 
	\frac{1}{\lambda} \sum \limits_{j=2} ^\ell   \sum \limits_{i=1}^{\ell-j+2} q_{j-1,\ell+1-i}  \delta_i c_{j-1}
	-\frac{1}{\lambda}\binom{\lambda}{2} c_{\ell-1}-\cdots-\frac{1}{\lambda}\binom{\lambda}{\lambda+1} c_{\ell-\lambda}.
\]

Nous posons alors

\[
\mathcal C_\ell =\sum \limits_{j=2} ^\ell   \sum \limits_{i=1}^{\ell-j+2}  \left( q_{j-1,\ell+1-i}  \frac{\delta_i}{\lambda}   c_{j-1} -q_{j-1,\ell-i} \delta_i  c_j \right)
-\frac{1}{\lambda}\sum \limits_{k=1} ^\lambda \binom{\lambda}{k+1} c_{\ell-k}
\mbox{ et }
\tilde{\mathcal C}_\ell =\frac{1} { k_\ell }   \mathcal C_\ell .
\]

En rassemblant les \'equations \eqref{eq:pout1} et \eqref{eq:pout2},  nous trouvons :
\[
f^{(\ell)} \left( \frac{ r {\lambda} }{{\lambda}^{\ell-1}} \right) 
	= \frac{1} { k_\ell}  
	 \sum \limits_{j=2} ^{\ell-1}   \sum \limits_{i=2}^{\ell-j+1} q_{j-1,\ell -i} \delta_i   \cdot   \left( S^{(j)}_r - p_j(r) \right)  + \tilde{\mathcal C}_\ell.   
\]

Nous trouvons donc :
\[
f^{(\ell)} \left( \frac{ r {\lambda} }{{\lambda}^{\ell-1}} \right) 
	= \sum \limits_{j=2} ^{\ell-1}   \sum \limits_{i=2}^{\ell-j+1} 
 	q_{j-1,\ell -i} \delta_i   \cdot  f^{(j)}   \left( \frac{r }{{\lambda}^{\ell-1}} \right)   \cdot \frac{k_j}{k_\ell} + \tilde{\mathcal C}_\ell.   
\]

D'apr\`es l'\'equation \eqref{eqiet}, entre deux valeurs de $ \frac{i {\lambda} }{{\lambda}^{\ell-1}} $, la fonction $f_{\ell}$ est une somme renormalis\'ee d'au plus ${\lambda}$ termes de $f^{(\ell-1)}(x)$. 

Remarquons que pour $1\leq k \leq \ell-2$ : $\binom{\ell-1}{2}- \binom{\ell-k-1}{2} \geq \binom{\ell-1}{2}- \binom{\ell-3}{2}=2\ell-5$.
Rappelons que les poyn\^omes $q_{i,j}$ sont \`a valeurs dans $\NN$.
En notant $\|\delta\| = \sup\{|\delta_i|;i\in \NN^*\} = \max \{|\delta_i|;1\leq i \leq \lambda\}$, $M = \| \delta \| \cdot \lambda^5$
et puisque $|\delta_2|\geq 1$ :
\[
||f^{(\ell)}||_\infty \leq ||f^{(\ell-1)}||_\infty+ M \frac{1}{\lambda^{2\ell}}
 	\left( \sum \limits_{j=2} ^{\ell-2}   \sum \limits_{i=2}^{\ell-j+1} q_{j-1,\ell -i}  \right) \cdot  ||f^{(j)}||_\infty
 	+{\lambda} \frac{\| \delta \|}{\lambda^{\ell-2}}    ||f^{(\ell-1)}||_\infty+| \tilde{\mathcal C}_\ell|.
\]

Notons alors $F_r = \sup \{ \| f^{(j)}\|;1\leq j \leq r \}$, en utilisant la majoration du lemme~\ref{le:techniquepourq} et $\theta = \frac{2\lambda-1}{\lambda^2}\in ]0,1[$ :
\begin{equation} \label{eq:groscalculencore}
||f^{(\ell)}||_\infty \leq \left(1+\frac{M}{\lambda^{\ell}} \right) ||f^{(\ell-1)}||_\infty+ M \cdot \theta ^\ell \cdot  \binom{\ell}{2}  \cdot F_{\ell-2} 	+| \tilde{\mathcal C}_\ell|.
\end{equation}

On \'etudie maintenant le comportement de la suite $ |\tilde{\mathcal C}_\ell|$.
\[
\mathcal C_\ell =  \sum \limits_{i=2}^{\lambda} q_{1,\ell+1-i}  \delta_i c_{1} 
	+\sum \limits_{j=2} ^{\ell-1}   \sum \limits_{i=2}^{\ell-j+1} \left( \frac{1}{\lambda}  q_{j,\ell+1-i}-q_{j-1,\ell-i} \right)  \delta_i c_{j}
	-\frac{1}{\lambda}\sum \limits_{k=1} ^\lambda \binom{\lambda}{k+1} c_{\ell-k}.
\]
Le coefficient de $c_{\ell-1}$ est nul dans la double somme.
En notant $\| \lambda \| = \sup\{ \binom{\lambda}{j}  ; j\geq 1\}$ et $\tilde c_j = \sup \{|c_k|;1\leq k \leq j\}$ :
\[
|\mathcal C_\ell | \leq 	
	2 \lambda (2\lambda-1)^\ell \| \delta \| \tilde c_{\ell-2}
	+||\lambda|| \cdot \tilde c_{\ell-1}
\]

Puisque $\frac{1}{k_\ell} c_\ell$ converge, cette suite est born\'ee et il existe une constante $D$ telle que pour tout entier $\ell$ et
pour $k\in \{1,\ldots, \ell-2\}$, $\frac{1}{k_\ell} |c_k| \leq \frac{1}{\lambda^{2\ell}} D$ et 
\[
\frac{\tilde c_{\ell-2}}{k_\ell} \leq \frac{1}{\lambda^{2\ell}} D \mbox{ et } \frac{\tilde c_{\ell-1}}{k_\ell} \leq \frac{1}{\lambda^{\ell}} D .
\]
Finalement il existe une constante $\Delta$ tel que
\begin{equation} \label{lequation}
| \tilde C_\ell | = \left| \frac{C_\ell}{k_\ell} \right|  \leq \Delta \left( \theta^\ell+\frac{1}{\lambda^\ell} \right).
\end{equation}

En rassemblant \eqref{eq:groscalculencore} et \eqref{lequation} :
\[
\| f^{(\ell)} \|_\infty \leq \left(1+\frac{M}{\lambda^{\ell}}  + M \theta^\ell \binom{\ell}{2} \right) F_{\ell-1} + \Delta \left( \theta^\ell+\frac{1}{\lambda^\ell} \right).
\]
On peut alors appliquer le lemme~\ref{lem:convergence_suite} qui permet  de conclure que la suite $(\| f^{(\ell)}\|)_\ell$ est convergente donc born\'ee.
\end{proof}

\subsection{Fin de la preuve.}
\begin{lemme} \label{lem:fi_lipshitz}
Pour tout $i \geq 1$, tout $x,y \in [0,\lambda]$ et $\u \in \{a,b\}$ :
\[
\left|f_\u^{(i)}(y) - f_\u^{(i)}(x) \right| \leq \frac{\lambda}{\delta_2} \left\|f^{(i-1)}_\u \right\| \left( |y - x| + \frac{2}{\lambda^{i-1}} \right).
\]
De plus, pour tout $i \geq 1$ et $x,y \in \RR^+$ :
\[
\left|f_\w ^{(i)}(y) - f_\w ^{(i)}(x) \right| \leq \frac{\lambda}{\delta_2} \left\|f_\w^{(i-1)} \right\| \left( |y - x| + \frac{2}{\lambda^{i-1}} \right).
\]
\end{lemme}

\begin{proof}
Nous faisons la preuve avec $\u=a$.
Soit $x_i = \lfloor \lambda^{i-1} x \rfloor / \lambda^{i-1}$ et $y_i = \lfloor \lambda^{i-1} y \rfloor / \lambda^{i-1}$.
Alors $f_a^{(i)}(x_i) = f_a^{(i)}(x)$ et $f_a^{(i)}(y_i) = f_a^{(i)}(y)$. Maintenant pour tout $\u \in K_\sigma$
\[
|\psi_i(T^{\lambda^{i-1} y_i} \u) - \psi_i(T^{\lambda^{i-1} x_i} \u) | 
= | S^{(1)}_{\lambda^{(i-1}(y_i - x_i)}(\psi_{i-1}, T^{\lambda^{i-1} x_i} \u)|
\leq \lambda^{i-1} |y_i - x_i|  \cdot \|\psi_{i-1}\|.
\]
Et donc
\[
|f^{(i)}_a(y_i) - f^{(i)}_a(x_i)| \leq
\frac{\lambda^{i-1} \delta_2^{i-1} \lambda^{(i-2)(i-3)/2}}{\delta_2^i \lambda^{(i-1)(i-2)/2}} |y_i - x_i| \cdot \|f^{(i-1)}_a\| = \frac{\lambda}{\delta} \|f^{(i-1)}_a\|  \cdot |y_i - x_i|.
\]
Il suffit alors de remarquer que $|y_i - x_i| \leq |y-x| + 2/\lambda^{i-1}$ pour conclure.

La m\^eme preuve montre la deuxi\`eme relation de l'\'enonc\'e.
\end{proof}

\begin{proposition}
Les fonctions $f_\w$, $f_a$ et $f_b$ sont continues et la fonction $f$ est solution de l'\'equation int\'egrale $(E_{\lambda,\delta_2})$.
\end{proposition}

\begin{proof}
	La continuit\'e des fonctions est une cons\'equence directe du lemme~\ref{lem:fi_lipshitz} et du fait que les fonctions $f_a^{(i)}$, $f_b^{(i)}$ et $f^{(i)}$ sont uniform\'ement born\'ees (proposition~\ref{prop:borne_uniforme}).

En reprenant l'\'equation \eqref{eqiet} et puisque $\binom{i-1}{2} +i-1=\binom{i-2}{2}$, nous trouvons pour tout entier $n$ : 
\[
\int_0 ^{n/{\lambda}^{i-1}} f^{(i)} (t) \dt = \sum \limits_{j=0} ^{n-1} f^{(i)} \left( \frac{j}{{\lambda}^{i-1}} \right) \frac{1}{{\lambda}^{i-1}}
=
\delta_2 f^{(i+1)} \left( \frac{n}{ \lambda^{i} } \right)-\delta_2f^{(i+1)}(0).
\]
En passant \`a la limite, la fonction $f$ v\'erifie l'\'equation int\'egrale en tout point $\lambda$-adique.
Par densit\'e des points $\lambda$-adiques et par continuit\'e de la fonction, $f$ est solution de l'\'equation int\'egrale.
\end{proof}


\appendix

\section{Le cas $\delta_2 = 0$}
L'hypoth\`ese $\delta_2\neq 0$ dans le th\'eor\`eme~\ref{thm:construction} est essentielle. Nous expliquons quels r\'esultats restent valables dans ce cas et pr\'esentons le r\'esultat de simulations.

Soit $\sigma$ une substitution sur $\{a,b\}$ de longueur constante $\lambda$ et telle que $|\sigma(a)|_a = |\sigma(b)|_a$ On suppose \'egalement que $\sigma(a) \not= \sigma(b)$ (car sinon le point fixe de $\sigma$ est un mot p\'eriodique). On sait par la proposition~\ref{prop:injection_monoide_dans_G} qu'il existe un entier $2 \leq k < \lambda$ tel que $\delta_2 = \delta_3 = \ldots = \delta_{k-1} = 0$ et $\delta_k \not= 0$. Dans cet article nous avons construit une solution \`a l'\'equation~\eqref{eq:E_lambda_delta} lorsque $k=2$.

De mani\`ere g\'en\'erale, le th\'eor\`eme~\ref{thm:cobord_infini} s'applique. En particulier, nous savons que l'on peut construire une solution $\psi = (\psi_1,\psi_2,\ldots)$ \`a l'\'equation cohomologique $\phi = (U_T - A) \psi$ d\`es que $\phi: \{a,b\} \rightarrow \RR$ est de moyenne nulle. En particulier, toutes les sommes de Birkhoff (recentr\'ees) sont born\'ees. 
Il semblerait que les $2n$-i\`emes sommes it\'er\'ees (respectivement les $2n+1$-i\`emes sommes it\'er\'ees) le long d'un point fixe de $\sigma$ convergent lorsque le temps est renormalis\'e par $\lambda^n$ (resp. $\lambda^{n+2}$).

\begin{figure}[H]  
\begin{center}
\includegraphics[width=0.95\textwidth]{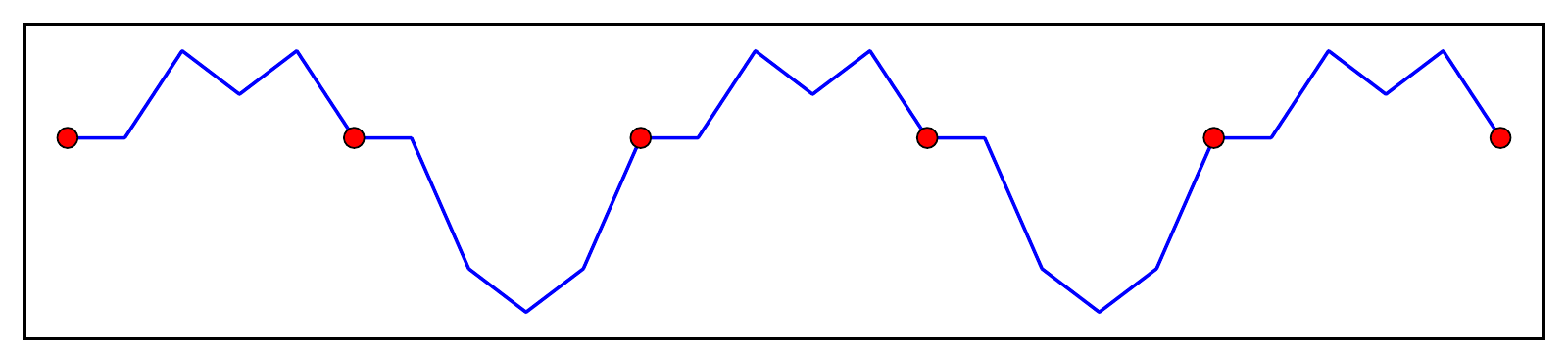} \\
\includegraphics[width=0.95\textwidth]{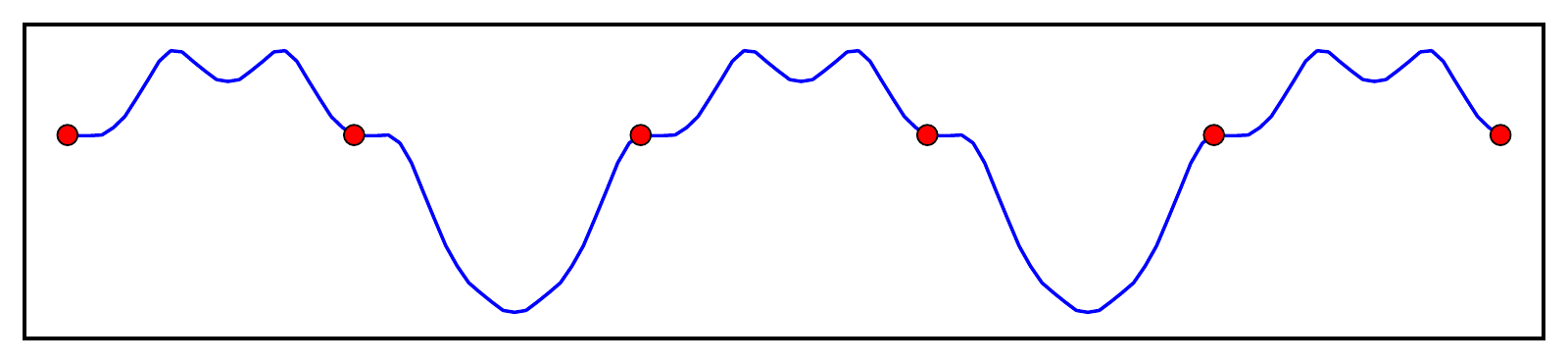} \\
\includegraphics[width=0.95\textwidth]{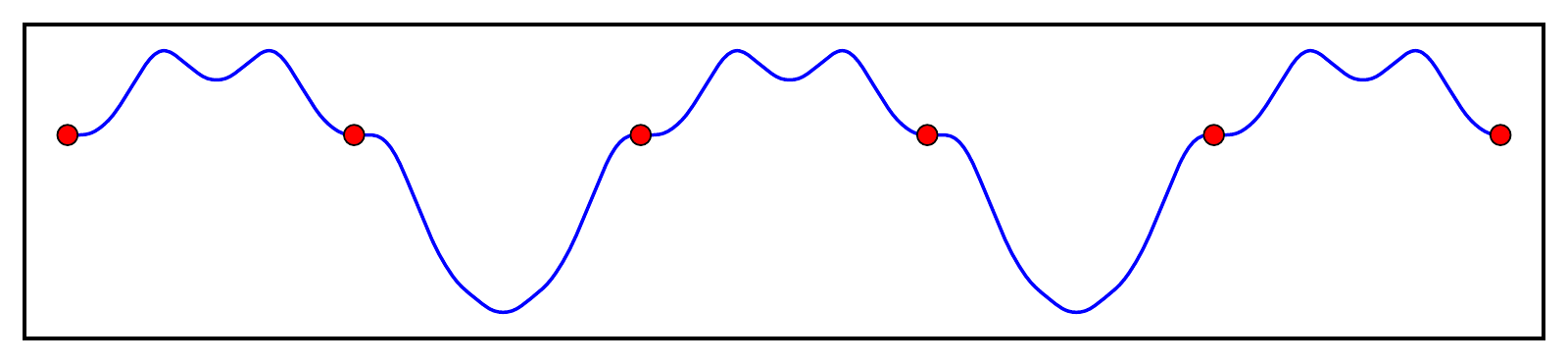} \\
\end{center}
\caption{Sommes $S^{(2n)}_{\lambda^{n}}(\phi,\sigma^{2n}(a))$ pour $n=1,2,3$ et la substitution $a \mapsto ababa, b \mapsto baaab$.}
\label{fig:fa_pair_delta2_nul}
\end{figure}

\begin{figure}[H]
\begin{center}
\includegraphics[width=0.95\textwidth]{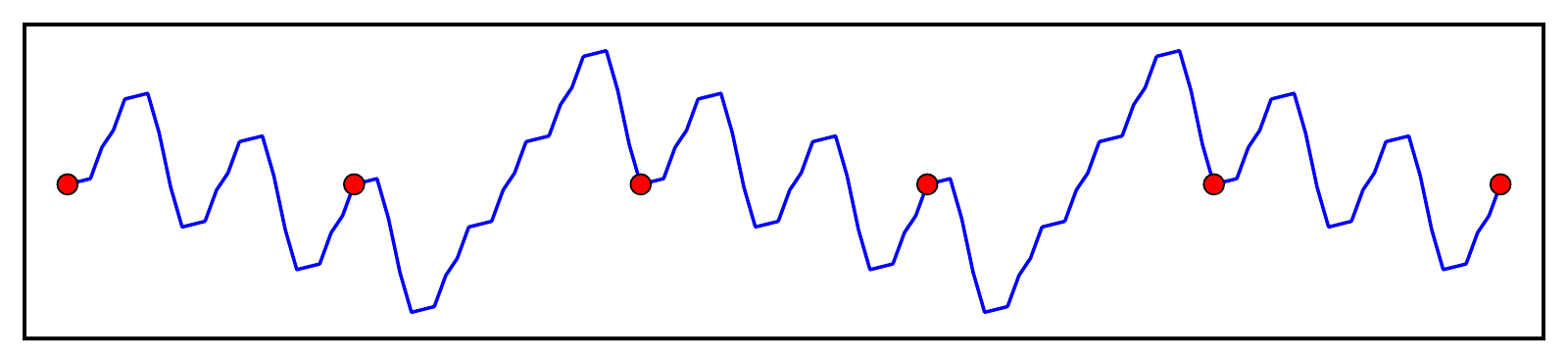} \\
\includegraphics[width=0.95\textwidth]{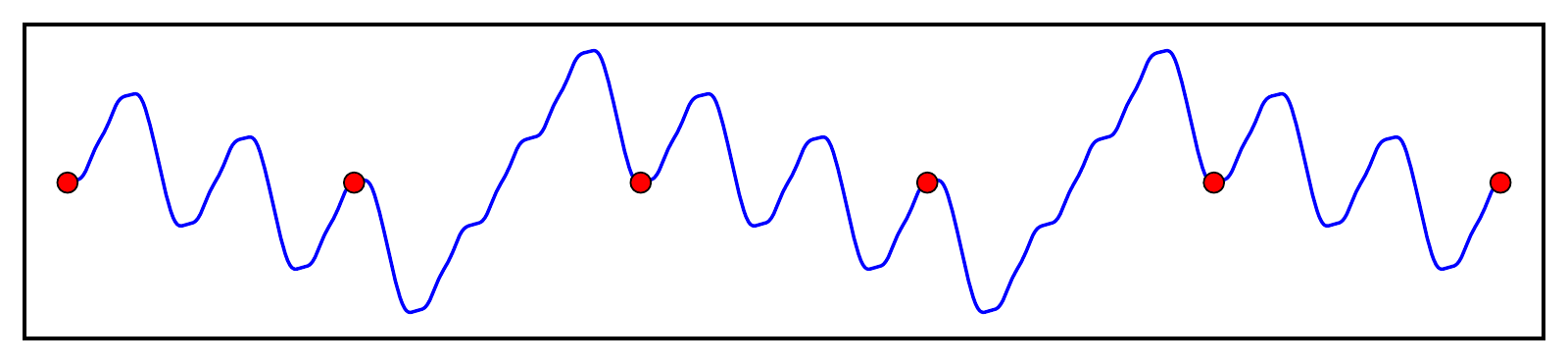} \\
\includegraphics[width=0.95\textwidth]{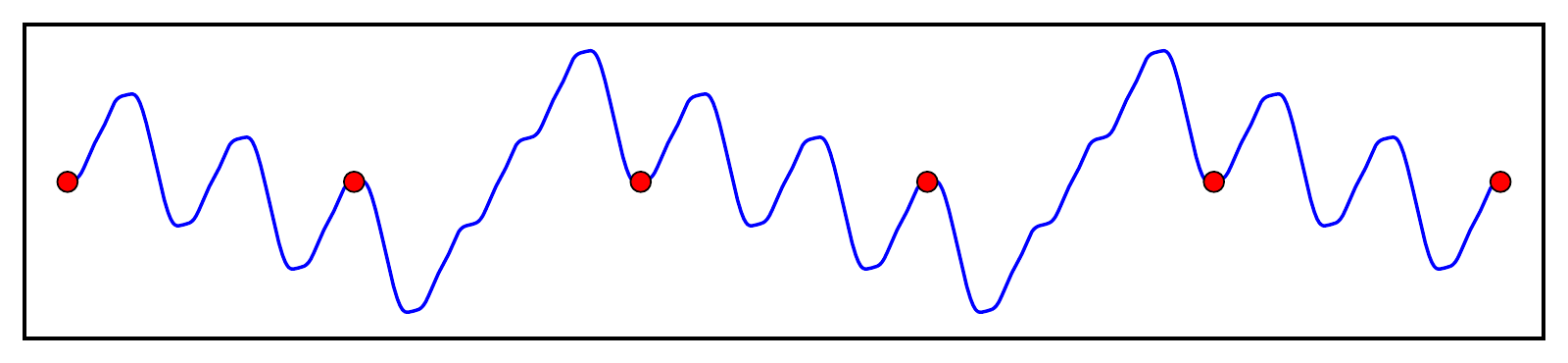}
\end{center}
\caption{Sommes $S^{(2n+1)}_{\lambda^{n+2}}(\phi,\sigma^{2n+1}(a))$ pour $n=1,2,3$ et la substitution $a \mapsto ababa, b \mapsto baaab$.}
\label{fig:fa_impair_delta2_nul}
\end{figure}

Les deux exemples les plus simples sont donn\'es par les substitutions
\[
\sigma_1: \left\{ \begin{array}{l@{\mapsto}l}
a & ababa \\
b & baaab 
\end{array} \right.
\quad \text{et} \quad
\sigma_2: \left \{ \begin{array}{l@{\mapsto}l}
a & abbaa \\
b & baaba
\end{array} \right.
\]
Elles v\'erifient $\delta(\sigma_1) = (0,0,3,3,1)$ et $\delta(\sigma_2) = (0,0,2,3,1)$.
Il semblerait que pour ces deux substitutions, les sommes $S^{(2n)}_{\lambda^n}(\phi,\sigma^{2n}(a))$ et $S^{(2n+1)}_{\lambda^{n+2}}(\phi,\sigma^{2n+1}(a))$ convergent et vers des limites distinctes. On a trac\'e ces sommes dans les figures~\ref{fig:fa_pair_delta2_nul} et~\ref{fig:fa_impair_delta2_nul}. Les points sur les dessins correspondent aux "valeurs enti\`eres" de la limite, c'est-\`a-dire aux temps $k \lambda^{n-1}$ pour $S^{(2n)}$ et aux temps $k \lambda^{2n}$ pour $S^{(2n+1)}$. On voit bien que ces valeurs sont toutes identiques et ceci peut-\^etre prouv\'e en suivant simplement la d\'emonstration de~\ref{thm:cobord_infini} ou bien en calculant les it\'er\'es de $\cL$.

\section{Une construction S-adique}
Soit $\sigma$ une substitution de longueur constante $\lambda$ telle que $\delta_1(\sigma) = 0$ et $\delta_2(\sigma) \not= 0$.
Dans le th\'eor\`eme~\ref{thm:thm_principal_precis} nous avons construit deux fonctions $f_a:[0,\lambda] \rightarrow \RR$ et $f_b:[0,\lambda] \rightarrow \RR$ comme limite des $k$-i\`emes sommes de Birkhoff sur respectivement $\sigma^k(a)$ et $\sigma^k(b)$. \'Etant donn\'e un mot infini $\u \in \{a,b\}^\NN$ on peut construire une fonction $f_\u: \RR_+ \rightarrow \RR$ en posant
\[
f_\u(x) = f_{\u_m}(x - m \lambda) \qquad \text{o\`u $m = \lfloor x /\lambda \rfloor$.}
\]
En adaptant la preuve du th\'eor\`eme~\ref{thm:thm_principal_precis}, il est facile de voir que
\begin{equation}
\label{eq:eq_fonc_sub}
\int_0^{\lambda x} f_{\sigma(\u)} \ds = \delta_2\ (f_{\u}(x) - f_{\u}(0)).
\end{equation}
Lorsque $\u = \sigma(\u)$ on retrouve l'\'equation~\eqref{eq:E_lambda_delta}.

\bigskip

Plus g\'en\'eralement, il semble possible de construire des fonctions limites lorsque les substitutions varient. Les constructions de ce type sont appel\'ees S-adiques (voir~\cite{DLR} ou bien \cite{BD-survol_S_adique}). Nous pr\'esentons un exemple associ\'e aux deux substitutions suivantes
\[
\sigma_1: a \mapsto aaabb, b \mapsto aabab
\quad \text{et} \quad
\sigma_2: a \mapsto ababb, b \mapsto abbab.
\]
Nous avons choisi ces deux substitutions car elles ont la m\^eme longueur $\lambda = 5$ et le m\^eme vecteur $\delta = (0,1,1,0,0)$.

Pour tout mot infini $\omega= \omega_0 \omega_1 \ldots \in \{1,2\}^\NN$ (on dira \emph{suite directrice}) on peut d\'efinir un mot infini comme limite
\[
\u(\omega) = \lim_{n \to \infty} \sigma_{\omega_0} \sigma_{\omega_1} \ldots \sigma_{\omega_{n-1}} (a).
\]
Par exemple pour $\omega = 1111\ldots$ on retrouve le point fixe de $\sigma_1$ et pour $\omega = 2222\ldots$ on retrouve le point fixe de $\sigma_2$. Plus g\'en\'eralement tout suite directrice p\'eriodique correspond \`a un point fixe de substitution. \'Etant donn\'e cette suite $\omega$ on peut \'egalement lui associer des d\'ecalages $K_{\omega,\NN}$ et $K_{\omega,\ZZ}$.

La preuve du th\'eor\`eme~\ref{thm:cobord_infini} s'adapte directement : toute fonction $\phi:\{a,b\} \rightarrow \RR$ est un cobord infini sur $K_{\omega,\ZZ}$. Cependant, si nous prenons na\"ivement la suite des cobords successifs $\psi_0$, $\psi_1$, \ldots il n'y a pas de renormalisation convenable qui fasse converger les fonctions.

On consid\`ere l'ensemble $\{1,2\}^\ZZ$ des mots bi-infinis sur $S$ et les fonctions $\phi_1 = 2 \chi_a - 3 \chi_b$ et $\phi_2 = 3 \chi_a - 2 \chi_b$. Si $\omega_0 = 1$ alors $\phi_1$ est de moyenne nulle sur $K_{\omega,\ZZ}$ alors que si $\omega_0 = 2$ alors $\phi_2$ est de moyenne nulle.

\begin{figure}[H]
\includegraphics{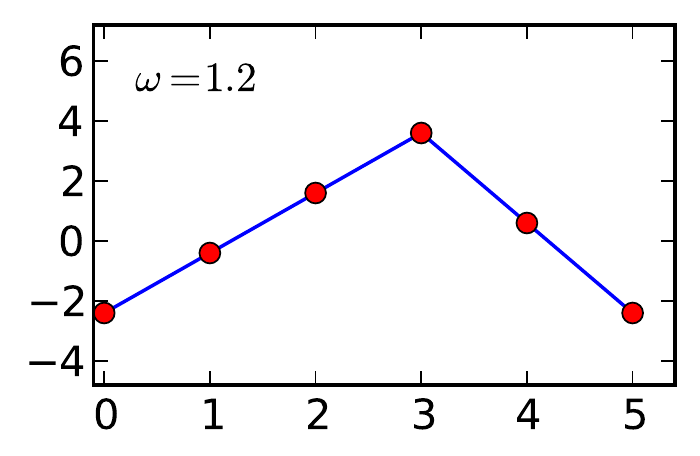}   \hspace{1cm}   \includegraphics{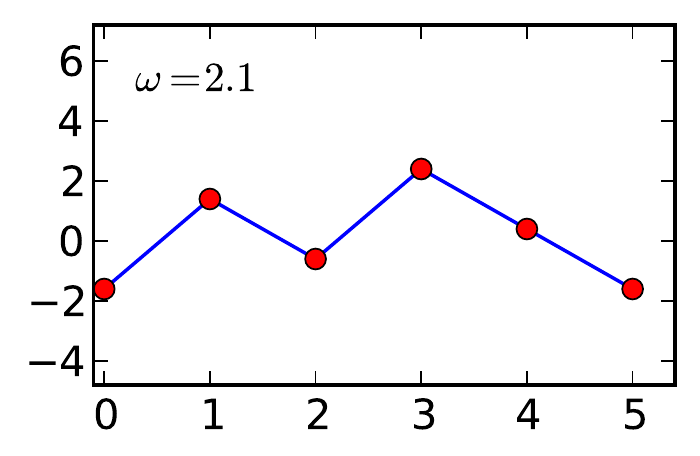}      \\
\includegraphics{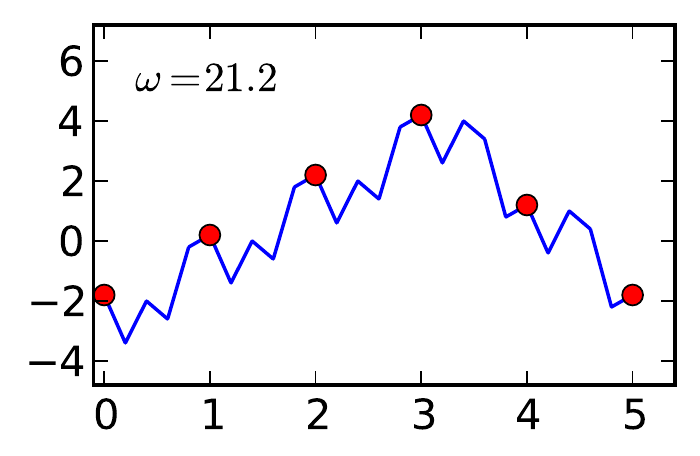}  \hspace{1cm}   \includegraphics{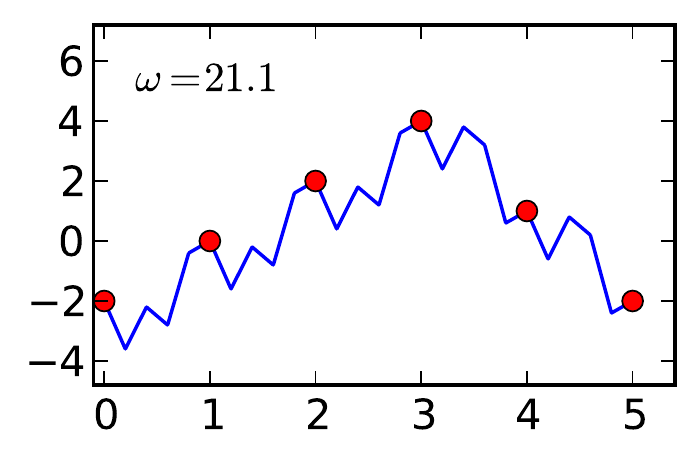}     \\
\includegraphics{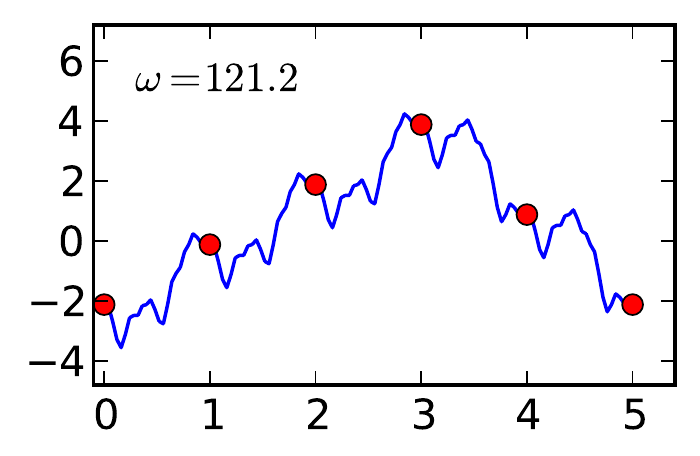} \hspace{1cm}   \includegraphics{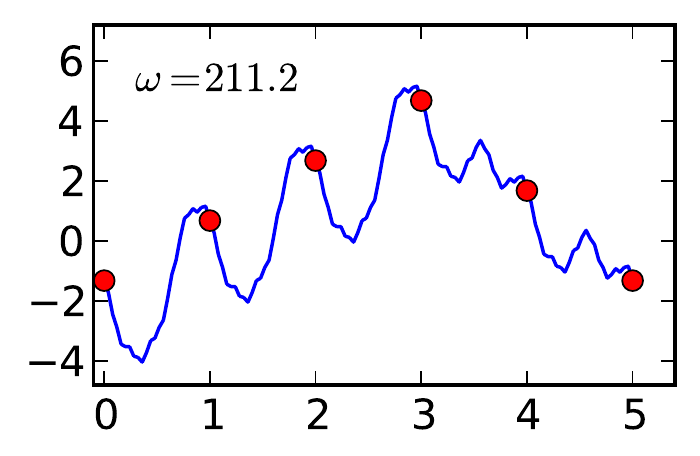}    \\
\includegraphics{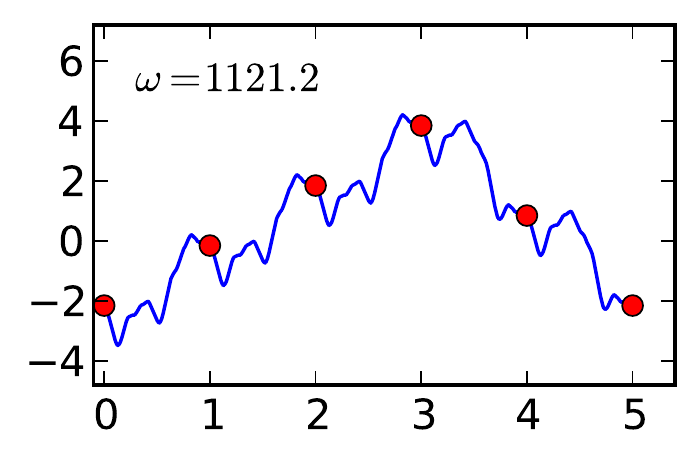} \hspace{1cm}  \includegraphics{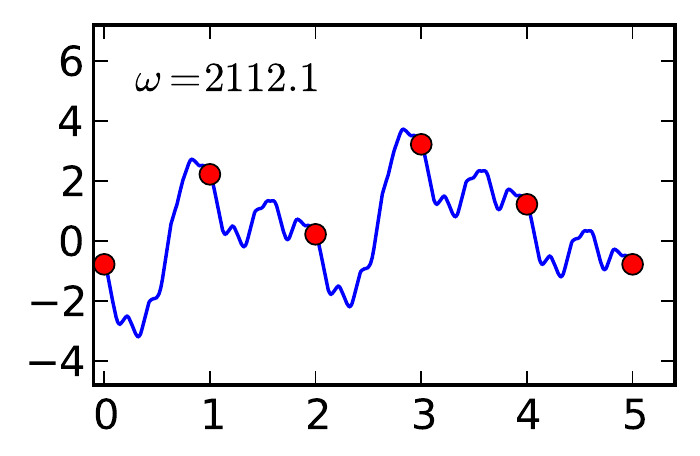} \\
\includegraphics{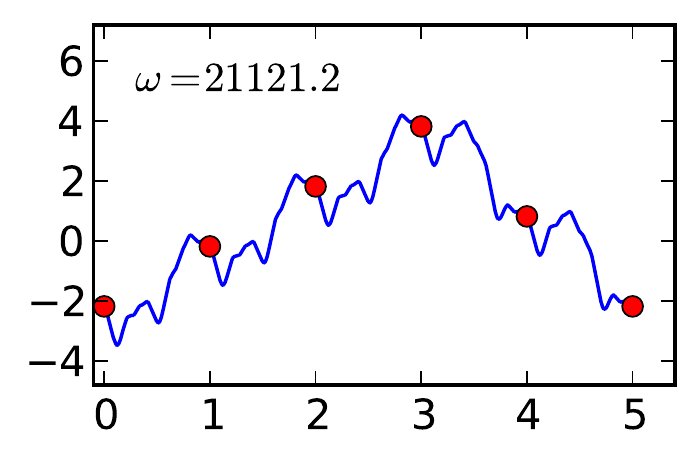}
\caption{Dessins des premi\`eres fonctions $f^{(i)}_{\omega,a}$. Les graphiques \`a gauche sont obtenus en prolongeant $\omega$ vers la gauche tandis que les graphiques sur la droite sont obtenus en prolongeant $\omega$ vers la droite. Remarquer que les diff\'erences entre les valeurs aux points entiers ne d\'ependent que de $\omega_{-1}$.}
\label{fig:fa_adique}
\end{figure}

On note $\ZZ_+ = \{0,1,2,\ldots\}$ et $\ZZ_- = \{-1,-2,\ldots\}$. \'Etant donn\'e un mot infini, on marquera la position de l'indice $0$ avec un point ainsi : $\omega = \ldots \omega_{-2} \omega_{-1} . \omega_0 \omega_1 \ldots$.
\`A chaque mot $\omega \in S^{\ZZ_- \cup \{0\}}$ nous allons construire un couple de fonctions $f_a: [0,\lambda] \rightarrow \RR$ et $f_b: [0,\lambda] \rightarrow \RR$. On d\'efinit $\psi_{i,\omega}: K_{T^{-i} \omega, \ZZ} \to \RR$ comme le $i$-\`eme cobord de la fonction $\phi_{\omega_{-i}}$. Bien s\^ur, $\psi_i$ ne d\'epend que de $\omega_{-i}, \omega_{-i+1} \ldots \omega_{-1}$ et $\omega_0$. De plus, il est constant sur les cylindres $\Cyl_\omega(i,m,\alpha) := T^m \sigma_{\omega_{-i}} \sigma_{\omega_{-i+1}} \ldots \sigma_{\omega_{-1}} ([\alpha])$. On pose
\[
f_{\omega,a}^{(i)}(x) = \frac{\psi_{i,\omega}(\Cyl_\omega(i,m,a))}{5^{(i-1)(i-2)/2}}
\quad \text{et} \quad
f_{\omega,b}^{(i)}(x) = \frac{\psi_{i,\omega}(\Cyl_\omega(i,m,b))}{5^{(i-1)(i-2)/2}}
\]
o\`u $m = \lfloor 5^{i-1} x \rfloor$. En suivant la preuve de la proposition~\ref{lem:convergence_points_rat}, on peut montrer que les diff\'erences entre deux points entiers successifs de $f_{\omega,a}$ et $f_{\omega,b}$ sont constantes et ne d\'ependent que de $\sigma_{\omega_{-1}}$. Nous ne le d\'emontrons pas formellement, mais $f^{(i)}_{\omega,a}$ et $f^{(i)}_{\omega,b}$ ainsi d\'efinies convergent. Voir les graphiques de la figure~\ref{fig:fa_adique}.

En admettant que la construction pr\'ec\'edente est valide, nous pourrions associer \`a chaque mot bi-infini $\omega \in S^\ZZ$ une fonction $f_\omega$ en concat\'enant $f_{a,\omega}$ et $f_{b,\omega}$ (qui ne d\'epend que de $\omega_-$ et $\omega_0$) selon le motif de $\u(\omega)$ (qui ne d\'epend que de $\omega_+$). La famille de fonctions $f_\omega$ v\'erifirait alors l'\'equation fonctionnelle
\[
\int_{\lambda x}^{\lambda y} f_\omega(s) \ds = \delta (\sigma_{\omega_0}) \Big(f_{T \omega}(y) - f_{T \omega}(x)\Big).
\]
o\`u $T: S^\ZZ \rightarrow S^\ZZ$ est le d\'ecalage sur les suites directives.


\begin{thebibliography}{9}

\bibitem[Ad04]{adam1}
Adamczewski, B.
\textit{Symbolic discrepancy and self-similar dynamics},
Ann. Inst. Fourier
\textbf{54}
(2004),
2201--2234.

\bibitem[BD14]{BD-survol_S_adique}
Berth\'e, V. et Delecroix, V.
\textit{Beyond substitutive dynamical systems: S-adic expansions},
\`a para\^itre dans RIMS Lecture note `Kokyuroku Bessatu'.

\bibitem[DLR13]{DLR}
Durand F., Leroy J. et Richomme G.
\textit{Do the properties of an S-adic representation determine factor complexity?},
J. of Integer sequences
\textbf{16} (2013)

\bibitem[BDMO08]{BogachevaDerfelbMolchanovcOckendond}
Bogachev L., Derfel G., Molchanov S. and Ockendon J.,
\textit{On bounded solutions of the balanced generalized pantograph equation},
Topics in stochastic analysis and nonparametric estimation,
IMA Vol. Math. Appl.,
Springer, New York,
\textbf{145} (2008) 24--49.

\bibitem[Be12-a]{jeff1}
Bertazzon, J.-F.
\textit{Symbolic approach and induction in the Heisenberg group}, 
Discret and Cont. Dyn. Syst
\textbf{32} (2012) 4, 1209--1229.

\bibitem[Be12-b]{jeff2}
Bertazzon, J.-F.
\textit{Resolution of an integral equation with the Thue-Morse sequence},
Indagationes Mathematicae,
\textbf{23} (2012) 4, 327--336.

\bibitem[Fa66]{Fabius66}
Fabius, J.
\textit{A probabilistic example of a nowhere analytic $C^\infty$-function},
Zeitschrift für Wahrscheinlichkeitstheorie und Verwandte Gebiete
\textbf{5} (1966), Issue 2, pp 173--174.

\bibitem[Ge94]{MR1326950}
Gelbrich, G.
\textit{Self-similar periodic tilings on the {H}eisenberg group},
Journal of Lie theory
\textbf{4}
(1994)
31--37.

\bibitem[Mo96]{Mosse}
Moss\'e, B.
\textit{Reconnaissabilit\'e des substitutions et complexit\'e des suites automatiques},
Bull. Soc. Math. France
\textbf{124}
(1996)
2
329--346.

\bibitem[Pi00]{pinner} 
Pinner, C.G. 
\textit{On the one-sided boundedness of sums of fractional parts {$(\{n\alpha+\gamma\}-\frac 12)$}},
Journal of Number Theory, \textbf{81}, (2000) 1, 170--204.	

\bibitem[Pr11]{prunescu}
Prunescu, M.
\textit{The \mbox{T}hue-Morse-Pascal double sequence and similar structures},
C. R. Acad. Sci.
\textbf{349} (2011) 939--942.

\bibitem[Py00]{Py}	
Fogg, N. Pytheas
\textit{Substitutions in dynamics, arithmetics and combinatorics},
Edited by V. Berth\'e, S. Ferenczi, C. Mauduit and A. Siegel.
Lecture Notes in Mathematics, Springer-Verlag, Berlin, 2002.

\bibitem[Yo06]{Yoneda2006}
Yoneda, T.
\textit{On the functional-differential equation of advanced type {$f'(x)=af(2x)$} with {$f(0)=0$}}
J. Math. Anal. Appl.,  \textbf{37} n. 1 (2006), 320--330.

\bibitem[Yo07]{Yoneda2007}
Yoneda, T.
{On the functional-differential equation of advanced type {$f'(x)=af(\lambda x)$}, {$\lambda>1$}, with {$f(0)=0$}},
J. Math. Anal. Appl., \textbf{332}, no. 1 (2007), 487--496.

\end{thebibliography}
\end{document}